\documentstyle[12pt]{article}
\begin{document}
\title{ Estimates  at or beyond
endpoint in harmonic analysis: Bochner-Riesz means and spherical
means
}
\author{{    Shunchao Long }}
\date{}
\maketitle
\begin{center}
\begin{minipage}{120mm}
\vskip 0.1in {
\begin{center}
{\bf Abstract}
\end{center}
\par   We  introduce  some new functions spaces
 to investigate some problems at or beyond endpoint. First,
  we prove that Bochner-Riesz means $B_R^\lambda$ are bounded
  from some subspaces of $L^p_{|x|^\alpha}$ to $L^p_{|x|^\alpha}$
  for   $ \frac{n-1}{2(n+1)}<\lambda \leq \frac{n-1}{2}, 0 < p\leq
p'_\lambda=\frac{2n}{n+1+2\lambda}, n(\frac{p}{p_\lambda}-1)<
\alpha<n(\frac{p}{p'_\lambda}-1)$, and $0<R<\infty,$ and so are  the
maximal Bochner-Riesz means $B_*^\lambda$ for $ \frac{n-1}{2}\leq
\lambda < \infty, 0 < p\leq 1$ and $-n< \alpha<n(p-1)$. From these
we obtain the $L^p_{|x|^\alpha}$-norm convergent property of
$B_R^\lambda $ for these $\lambda,p,$ and $\alpha$. Second, let
$n\geq 3,$  we prove that the maximal spherical means are bounded
from some subspaces of $L^p_{|x|^\alpha}$ to $L^p_{|x|^\alpha}$ for
$0<p\leq \frac{n}{n-1}$ and $ -n(1-\frac{p}{2})<\alpha<n(p-1)-n$. We
also obtain  a $L^p_{|x|^\alpha}$-norm convergent property of the
spherical means for such $p$ and $\alpha$. Finally, we   prove that
some new types of $|x|^\alpha$-weighted estimates hold at or beyond
endpoint for many operators, such as Hardy-Littlewood maximal
operator, some maximal and truncated singular integral operators,
the maximal Carleson operator, etc. The new estimates can be
regarded as some substitutes for the $(H^p,H^p)$ and $(H^p,L^p)$
estimates for the operators which fail to be  of types $(H^p,H^p)$
and $(H^p,L^p)$.

}
\end{minipage}
\end{center}
\vskip 0.1in
\baselineskip 16pt
\par

\par \begin{center}
{\bf Table of contents}
\end{center}
\par \textrm{1.}  Introduction

\par \textrm{2.}  $  BL^{p,s}_{|x|^{\alpha}}$     and their  basic properties

\par \textrm{3.}  A  molecular theorem

\par \textrm{4.}   Hardy-Littlewood maximal operator

\par \textrm{5.}  Bochner-Riesz means

 \par \textrm{6.}   Spherical means

\par \textrm{7.} Some   of other sublinear operators

\par \textrm{8.} Remarks
\par \section*
 {\bf  1.~~ Introduction }

\par ~~~~It is well-known that the Lebesgue spaces  $L^{p}=L^p({\bf R}^n)$ play an
important role in harmonic analysis. A lot of problems are solvable
on $L^{p}$, but usually only for $p$ in certain
subinterval $(A,A')$ of $(0,\infty)$, not outside of $(A,A')$ (i.e.
not at or beyond endpoint). For example, $L^p$-inequalities hold for
Carleson operator for $p$ in $E=(1, \infty)$, Bochner-Riesz means in
$E=(\frac{2n}{n+1+2\lambda}, \frac{2n}{n-1-2\lambda})$ (conjecture,
$0<\lambda\leq \frac{n-1}{2}$), the maximal spherical means in
$E=(\frac{n}{n-1}, \infty)$,   Hardy-Littlewood maximal operator in
$E=(1, \infty)$, etc, but do not for $p$ in ${\bf R}\setminus E$
(see \cite{S3}). The boundedness properties of these operators are
related to the corresponding convergence problems, respectively (see
\cite{S2,St61970}).

\par  Generally, $A$ and $A'$ are called endpoints.

When $p$ is at or  beyond endpoint, many operators are unbounded on
$L^p$,  so $L^p$ are too large to be the domains of these operators,
yet  too small to be the range since the target spaces of such
operators exceed $L^p$. Therefore,
 one generally  looks for  some larger spaces than
   $ L^p$ as the range spaces of these operators. The weak-$ L^p$ are suitable
   substitutes for  $ L^p$ for some operators, for example, weak-$
   L^1$ for $
   L^1$
    for Hardy-Littlewood maximal function $Mf$ which is of weak-type
    $(L^1,L^1)$. On the other hand,
 it is also natural  to look for some
subspaces of $L^{p}$ from which the operator considered is bounded.
 A famous example of such subspaces
 is
  the classical Hardy spaces $H^{p}$ with $0<p\leq 1$ which are
  some subspaces  of $L^{p}$ .
They are suitable
 substitutes for  $L^p $ for many  questions when $0<p \leq 1$, for
 example, for the boundedness estimates of certain singular integral
 operators.
  These operators are of type $(H^p,L^p)$ when $0<p \leq 1$ (see
\cite{S3}).

 \par But, $H^p, 0<p\leq 1,$ are  still not suitable for a lot
of problems in harmonic analysis, such as the convergent of Fourier
series and Bochner-Riesz means at the critical index,  the
boundedness property of Hardy-Littlewood maximal operator, etc. At
the same time, $H^p=L^p$ when $p>1$,   naturally, $H^p$ are not
suitable for the problems beyond endpoint which is larger than 1,
such as the convergent of Bochner-Riesz mean  when the index
$\lambda$ is less than the critical index $\frac{n-1}{2}$  and $ 1<
p\leq p'_\lambda=\frac{2n}{n+1+2\lambda}$,  the spherical means when
$1<p\leq \frac{n}{n-1}$, etc.

\par  In fact,
 when $p$ is   at or beyond endpoint, both weak-$L^p$ and $H^p$ are  not suitable for a lot
of problems such as boundedness properties of the maximal
Bochner-Riesz mean, the maximal spherical means and   the Carleson
operators, since the maximal Bochner-Riesz mean as $ 0< p\leq
p'_\lambda$ (see \cite{Stei, SG}), the maximal spherical means as
$0<p\leq \frac{n}{n-1}$ (see \cite{S3}) and   Carleson operators as
$0<p\leq 1$ (see \cite {Ko1}) fail to be of both weak-type
$(L^p,L^p)$ and type $(H^p,L^p).$

\par A natural question is then  how to extend the known solutions of these
problems   to the cases of at or beyond endpoints.

\par Surprisingly to me, we see here that a lot of problems are
solvable at or beyond endpoints if we consider the weighted Lebesgue
spaces $L^p_{|x|^\alpha}$, i.e.  the problems are solvable on some
subspaces of $L^p_{|x|^\alpha}$.  $L^p_{|x|^\alpha}$ may be more
suitable functions spaces to study certain problems beyond endpoint
than the non-weighted ones.


\par In this paper  we will  study some problems  at or beyond
endpoint  such as  boundedness properties of  Bochner-Riesz means,
the maximal spherical
 means and many other operators arising  in harmonic analysis, and norm convergent.

\par To this end, we will introduce some new   functions spaces
$BL^{p,s}_{|x|^{\alpha}}$ (in fact, they were first introduced by
the  author himself in \cite{Longs}).  One maybe initially want to
use these spaces as some subspaces of $L^p$, liking Hardy spaces
$H^p$ (which are subspaces of $L^p$, as we know) and the block
spaces $B_s$ introduced in \cite{TW1} (which are subspaces of
$L^1$), since they are analogous
 to the atom-$H^p$ and  $B_s$. However, we see here that these
spaces are subspaces of $L^{p}_{|x|^{\alpha}}$ for certain ranges of
$p, \alpha $ and $ s $,
 and are very useful when they are used as subspaces of
$L^{p}_{|x|^{\alpha}}$.

 \par The estimates of Carleson operators and a.e. convergent
 problems will be considered in our next paper.

\par I am greatly indebted  to
J.Wang   and D.S.Fan  for helpful discussions and constant support,
and to X.C.Li for useful discussions and suggestions when I visited
him in 2006.
\\{\bf  1.1  Bochner-Riesz means  }

\par The Bochner-Riesz means in ${\bf R}^n$ are defined as
$$(B_R^\lambda f)\hat{}(\xi)=(1-|\xi /R|^2)^\lambda_+\hat{f}(\xi), ~~~\lambda\geq 0, 0<R<\infty,$$
where $\hat{}$ denotes the Fourier transform. By a simple argument,
 that $B_R^\lambda f$ converges to $f$ in $L^p$-norm
 as $R$ tends to $\infty$ is equivalent to the boundedness in $L^p$ of $B_1^\lambda
 $ (see \cite{S3}),
 i.e.
$$
 \|B_1^\lambda f\|_{L^{p}}\leq
 C\|f\|_{L^{p}}.\eqno (1.1)
 $$

\par A great deal has been already known about (1.1).
Bochner \cite{B} proved that (1.2) holds for $\lambda>(n-1)/2$ and
$1\leq p \leq \infty $.
  For $0<\lambda \leq (n-1)/2 $, Herz  \cite{H} showed that, if (1.1)
holds, then
$$
 \frac{2n}{n+1+2\lambda}=p'_\lambda
 <p<p_\lambda=\frac{2n}{n-1-2\lambda}.\eqno (1.2)
 $$
\par Since then, a famous conjecture states that

\par {\bf Bochner-Riesz conjecture} For $0<\lambda\leq (n-1)/2$,  (1.1) holds under condition (1.2).

\par This conjecture was completely  proved in ${\bf R}^2$ by Carleson and Sj\"{o}lin \cite{CaS}, see also
\cite{F1}, \cite{Cord1} and \cite{Hor}.  In the higher dimensions,
it is still open, but there are several results. Stein \cite{S4}
showed that, for $0<\lambda \leq (n-1)/2 $, if $
 \frac{2(n-1)}{n-1+2\lambda}<p<\frac{2(n-1)}{n-1-2\lambda},
$ then (1.1) holds, and this includes that  the conjecture is true
for $\lambda = \frac{n-1}{2}. $ The argument of Fefferman and Stein
\cite{F2} and the restriction theorem of Tomas and Stein \cite{To}
imply the case of $\frac{n-1}{2(n+1)}<\lambda  , $ and Lee
\cite{Lee} extended this to the case of $\frac{n-2}{2(n+2)}<\lambda
. $
\par See also  \cite{B3}, \cite{B4},  \cite{W},
\cite{TV1}, \cite{TV2} and \cite{TVV} for some specific progress in
${\bf R}^3$.

\par   A natural endpoint version of the Bochner-Riesz conjecture is
the statement that
\par {\bf Endpoint Bochner-Riesz conjecture} For $0<\lambda\leq (n-1)/2$ and $p'_\lambda=\frac{2n}{n+1+2\lambda}$,
$B_1^\lambda$ is of weak-type $(L^{p'_\lambda},L^{p'_\lambda}).$

\par The endpoint conjecture was proved  in two  dimensions by Seeger \cite{Seeg},
 for
$\frac{n-1}{2(n+1)}<\lambda\leq \frac{n-1}{2}$  by Chist
\cite{Christ1, Christ2} and for  $\lambda=\frac{n-1}{2(n+1)} $ by
Tao \cite{T1} in the higher dimensions.

\par And Tao \cite{T2} showed that the weak-type endpoint
Bochner-Riesz conjecture is equivalent to the  standard
Bochner-Riesz conjecture.

\par But, beyond endpoint, i.e. $p<p'_\lambda$,  it seems that
$B_1^\lambda$ fails to be of weak-type $(L^p,L^p)$ for $\lambda\leq
(n-1)/2$, for example, it is known that $B_1^\lambda$ fails to be of
weak-type $(L^1,L^1)$ for all $\lambda < (n-1)/2$.
\par
{\bf Question } A natural question  is then what happens beyond the
endpoint, i.e. for $0<p<p'_\lambda=\frac{2n}{n+1+2\lambda},0<\lambda
\leq \frac{n-1}{2 } $.

\par
For the question above, it is  natural to consider the Hardy spaces
$H^p$ when $p\leq 1$. Stein,   Taibleson and   Weiss \cite{STW}
considered this and obtained some $H^p$ estimates of $B_1^\lambda$,
but only for $\lambda\geq \frac{n-1}{2}$. Unfortunately, $H^p$ is
not suitable  for the norm convergence of $B_R^\lambda$, even for
$\lambda= \frac{n-1}{2}$. In fact, by a result of Stein \cite{Stei}
(see also \cite{Sato} for the weighted case), there exists a
function $f\in H^1$ such that
$$\limsup_{R\rightarrow \infty}B^{(n-1)/2}_Rf(x)=\infty, ~~{\rm a.e.}.$$
This implies that the convergene of $B^{(n-1)/2}_Rf(x)$ to $f(x)$ in
$L^1$ fails as $R\rightarrow \infty$. At the same time, $H^p=L^p$
when $p>1$,  so $H^p$ is not a suitable space  for the norm
convergence of $B_R^\lambda$ at this time.

\par Although we did not find the subspaces of $L^p$ we are looking forward
to here,  but to our surprise, we found that there exist some
subspaces of $L^{p} _{|x|^{\alpha}}$ so that $B^{(n-1)/2}_Rf(x)$
converges to $f(x)$ in $L^{p} _{|x|^{\alpha}}$ for the functions $f$
in these subspaces for $p\leq p'_\lambda,$ i.e. at or beyond
endpoint.

\par The question above has also been studied for the radial functions and $p>1$
in \cite{Ande}, and  a  weighted inequality has been obtained for
$B_1^\lambda$.

\par The weighted estimate of $B_1^\lambda$ for $p$ in the range
$(p'_\lambda,p_\lambda)$ has been investigated by several authors.
For $p=2$,  Rubio \cite{Rubio}, Hirschman \cite{Hir1} showed that
$B_1^\lambda$ is bounded on
 $L^2_{|x|^\alpha}$ provided  $\lambda < \frac{n-1}{2}$ and
$n(\frac{2}{p_\lambda}-1)=-1-2\lambda< \alpha<1+2\lambda
=n(\frac{2}{p'_\lambda}-1)$.  Carbery, Rubio, and Vega \cite{CRV}
 proved that the maximal Bochner-Riesz mean  $B_*^\lambda=\sup_{R>0}|B_R^\lambda|$ is bounded on
 $L^2_{|x|^\alpha}$ for $  \lambda>0,p=2$ and
$n(\frac{2}{p_\lambda}-1)< \alpha\leq 0$.  For the
 critical  index $\lambda=\frac{n-1}{2}$ and $1<p<\infty$ Shi and Sun \cite{SS} proved that
$B_*^{\frac{n-1}{2}}$ is bounded on
 $L^p_{|x|^\alpha}$ provided
$n(\frac{p}{p_{(n-1)/2}}-1)=-n< \alpha<n(p-1)
=n(\frac{p}{p'_{(n-1)/2}}-1)$.

\par Naturally, it is conjectured that:
let $ 0 < \lambda\leq \frac{n-1}{2}, p'_\lambda <p<p_\lambda,
n(\frac{p}{p_\lambda}-1)< \alpha<n(\frac{p}{p'_\lambda}-1)$, then
$B_1^\lambda$ (even $B_*^\lambda$) is of type $(L^{p}
_{|x|^{\alpha}},L^{p} _{|x|^{\alpha}})$.

\par Our first aim of    this paper
 is to prove that an extension of this conjecture holds also
 to the case of at of
 beyond endpoint
 $p'_\lambda$, i.e. to $0<p\leq p'_\lambda$.
 That is
\par let $0<R<\infty,$ and let $ \frac{n-1}{2(n+1)}< \lambda\leq
\frac{n-1}{2},    0<p< p_\lambda' , n(\frac{p}{p_\lambda}-1)<
\alpha<n(\frac{p}{p'_\lambda}-1)$, then
$$ B_R^\lambda ~{\rm is ~ of ~type ~}
(BL^{p, p_\lambda}_{|x|^{\alpha}},L^{p} _{|x|^{\alpha}})   $$
 and (1.5) below holds for all $ f\in BL^{p, p_\lambda}_{|x|^{\alpha}}$.

 \par In fact, we have more:

\par let $
\frac{n-1}{2(n+1)}< \lambda\leq \frac{n-1}{2}, p'_\lambda
<s<p_\lambda, 0<p\leq s $ and $ n(\frac{p}{s}-1)\leq
\alpha<n(\frac{p}{p'_\lambda}-1)$. Then
  $$
B_1^\lambda ~{\rm is ~ of ~type ~} (BL^{p,s}_{|x|^{\alpha}},BL^{p,s}
_{|x|^{\alpha}}) , \eqno (1.3)
$$
and  let $0<R<\infty,$ and let $ \frac{n-1}{2(n+1)}< \lambda\leq
\frac{n-1}{2}, p'_\lambda <s<p_\lambda, 0<p< s , n(\frac{p}{s}-1)<
\alpha<n(\frac{p}{p'_\lambda}-1)$ and $\alpha\leq 0$, then
$$
B_R^\lambda ~{\rm is ~ of ~type ~} (BL^{p,s}_{|x|^{\alpha}},L^{p}
_{|x|^{\alpha}}) \eqno (1.4)
$$
 for all $0<R<\infty, $ and
$$
\|B_R^\lambda f-f\|_{L^p_{|x|^\alpha}} \rightarrow 0,~~~~
R\rightarrow \infty\eqno (1.5)
$$
for all $ f\in BL^{p,s}_{|x|^{\alpha}}$.


\par When $\lambda = \frac{n-1}{2}$,  these results above hold for the maximal Bochner-Riesz
operators. Noticing that $n(\frac{p}{p_{(n-1)/2}}-1)=-n$ and $
n(\frac{p}{p'_{(n-1)/2}}-1)=n(p-1)$, we will prove that:

\par  let $ \lambda = \frac{n-1}{2},
 0<p\leq 1,
-n< \alpha<n(p-1)$ and $ s\geq \frac{pn}{n+\alpha}$, then
  $$
B_*^\lambda ~{\rm is ~ of ~type ~}
(BL^{p,s}_{|x|^{\alpha}},BL^{p,s}_{|x|^{\alpha}})  , \eqno (1.6)
$$
$$ B_*^\lambda ~{\rm is ~ of ~type ~} (BL^{p,s}_{|x|^{\alpha}},L^{p} _{|x|^{\alpha}}) ,\eqno (1.7)
$$
 and (1.5) holds.
\par And when $\lambda > \frac{n-1}{2}$, these results above, i.e. (1.6), (1.7) and (1.5) hold also.
\\
{\bf 1.2 Spherical means }

\par For a  locally integrable function $f$
 on ${\bf R}^n $
 we define the spherical means
$$A_t (f)(x)= \int_{{\bf S}^{n-1}}f(x-t\theta)d\sigma(\theta)$$
for $t\in E\subseteq (0,\infty)$,  where $d\sigma$ is the
 normalized Lebesgue measure on the unit sphere ${\bf S}^{n-1}$.
\par Stein \cite{S1} showed that
$$\lim_{t\rightarrow0,t\in E}A_t (f)(x)=f(x)~~~~{\rm a.e.}\eqno(1.8)$$
when $E=(0,\infty)$, provided $f\in L^p, p>n/(n-1)$ and $n\geq 3.$
Bourgain \cite{B1} extended this result to the case of $n=2.$
Unfortunately, (1.8) fails for $p\leq n/(n-1),$  see \cite {S3}.
\par This
caused a type of problems:
   what happens  at or beyond the
 endpoint, i.e.
  $   p\leq
 n/(n-1)$?

\par It is showed  that for certain subset $E$ of $(0,\infty)$,  (1.8)   holds for
  $1<p<\infty$, furthermore, similar results are also proved for  the endpoint  $p=1$.
 In fact,   when
 $E=\{2^k:k\in {\bf Z}\}$, i.e. for the lacunary spherical means,   (1.8) was proved by Calder\'{o}n \cite{Ca} and Coifman and Weiss
 \cite{CW2} for $f$  in   $L^p $ with $1<p<\infty$,   by Christ \cite{Christ1} for $f$    in Hardy
space $H^1 $, and by
 Seeger, Tao and Wright  \cite{STW1} and \cite{STW2} for $f$
locally in $L\log \log L $. When
 $E$  satisfies some regularity conditions, (1.8) was proved for $f$  in   $L^p $ with
  $1<p\leq n/(n-1)$ by
 Seeger, Tao and Wright in \cite{STW3}, and Seeger, Wainger and Wright in
 \cite{SWW1} and \cite{SWW2},  etc.

\par According to a theorem by Stein \cite{S2}, pointwise convergence (1.8) for $f\in L^p$ is
equivalent with a weak type $(p,p)$  boundedness properties of the
maximal function
$${\mathcal{M}}_Ef(x)=\sup_{t\in E}|A_t (f)(x)|.$$

\par At the
same time, the boundedness properties of
${\mathcal{M}}={\mathcal{M}}_{(0,\infty)} $ on $ L^p_{|x|^\alpha}$
were obtained in \cite{DV}, that is, if $n\geq 2, n/(n-1)<p$, then
 $$
 {\mathcal{M}} \textrm {~is~ of~ type}~ ({L^p_{|x|^\alpha}},
 {L^p_{|x|^\alpha}}) \eqno(1.9)
 $$
  for
 $ 1-n<\alpha< n(p-1)- p $ and (1.9) does  not hold  for $\alpha <1-n
$. While for general   ${\mathcal{M}}_E$ this was considered in
\cite{DS}.

\par In this paper, the second work is to extend (1.9) to $p\leq n/(n-1)$,  as follows:
 let $  n\geq 3, 0<p\leq n/(n-1)$ and $ -n(1-p/2)\leq \alpha< n(p-1)-n $, then
$$
{\mathcal{M}} ~{\rm is ~ of ~type ~}
(BL^{p,2}_{|x|^{\alpha}},BL^{p,2}_{|x|^{\alpha}}), \eqno (1.10)
$$
and let $  n\geq 3, 0<p\leq n/(n-1)$ and $ -n(1-p/2)<\alpha<
n(p-1)-n $, then
$$
 {\mathcal{M}} \textrm {~is~ of~ type}~ (BL^{p,2}_{|x|^{\alpha}},
 {L^p_{|x|^\alpha}}). \eqno(1.11)
 $$

 \par It is worth to point out that at the endpoint ${\mathcal{M}}
$
 is a bounded operator from $L^{n/(n-1),1}$ to  $L^{n/(n-1),\infty}$
 if $n\geq 3$, where $L^{p,q}$ is the usual Lorentz space, and this was
 proved by Bourgain \cite{B2}. For general set $E$, these types of  estimates
 can be found in \cite{STW3} and \cite{SWW1}.
\\{\bf  1.3 Some  sublinear operators    }
\par For $p \in (1,\infty)$, it is known that type $(L^p,L^p)$ holds
 for many important operators arising in harmonic analysis, while
 for $p \in (0,1]$, there are several cases: the operator

 \par {\bf A}. fails  (or not know) to be of both weak-type $(L^1,L^1)$, and
 types
 $(H^p,H^p)$ and
 $(H^p,L^p)$ for $p \in (0,1]$, such as   Carleson
 operator \cite{Ko1,Zy}, the maximal Bochner-Riesz means at the critical index \cite{S2,H},
  certain oscillatory
 singular integral operators \cite{Pan6}, etc.,

\par {\bf B}. is of  weak-type $(L^1,L^1)$ but fails  (or not know) to be of type
 $(H^p,H^p)$ and  type
 $(H^p,L^p)$ for $p \in (0,1]$,
 such as Bochner-Riesz means at
the critical index, see \cite{Christ1,H,Zy}, some maximal operators
including
 Hardy-Littlewood maximal
 operator, the maximal  Calder\'{o}n-Zygmund operators,
 the maximal  Hilbert transform and  the maximal  Riesz
 transform,   some truncated operators including   the truncated Calder\'{o}n-Zygmund operators,
  the truncated Hilbert transform and the truncated Riesz
 transform,
 etc.,

\par {\bf C}. is of  weak-type $(L^1,L^1)$, and
 types
 $(H^p,L^p)$ and   $(H^p,H^p)$
 for $p $ in part of  $(0,1]$ but not  (or not know)
 for  $p $ in the remain part of  $(0,1]$. For example, it is not known whether  type
 $(H^p,L^p)$ and  type
 $(H^p,H^p)$ hold for $p$ in $(0,1/2]$
for Hilbert transform and the partial summation operators of Fourier
series, for $p$ in $(0,n/(n+\delta)]$ for Calder\'{o}n-Zygmund
operators. Furthermore,  it is known that
 type $(H^p,H^p)$ does not hold for $p$ in $(0,1)$ for strongly singular integral
 operators, see \cite {FS72,Sjolin76},
 and type $(H^p,L^p)$ does not hold for $p$ in $(0,1)$
 for certain oscillatory
 singular integral operators, see \cite {Pan3,Pan6},

\par {\bf D}. is of  both weak-type $(L^1,L^1)$, and  types
 $(H^p,L^p)$ and
 $(H^p,H^p)$  for full $p \in (0,1]$ such as Riesz transform.

\par
A natural question  is then what happens for these operators when
$p$ is  in the remain ranges.

 \par In this paper, our third aim is to consider the boundedness properties of
some  operators mentioned above when $p$ is in the remain ranges. We
will show that some new types of estimates hold when $p$ is in the
remain ranges of $(0, 1]$ for many of these operators, such as
Hardy-Littlewood maximal operator, some singular integral operators
as well as the corresponding maximal operators and truncated
operators, Carleson
 operator, the maximal Bochner-Riesz means at the critical index, Calder\'{o}n-Zygmund
 operator as well as the corresponding maximal operators and truncated
operators, the strongly singular integral operator, some oscillatory
singular integrals operators,
 etc..
 That is, we shall prove that
 they  are of
 $$ {\rm type~}(BL^{p,s}_{|x|^{\alpha}}, BL^{p,s}_{|x|^{\alpha}})  \eqno(1.12)$$
 for
$1<s<\infty, 0<p\leq s$ and $ -n(1-p/s)\leq \alpha<n(p-1)$,
 and
$$ {\rm type~}(BL^{p,s}_{|x|^{\alpha}}, {L}^{p} _{|x|^{\alpha}})\eqno(1.13)$$
 for
$1<s<\infty, 0<p< s, -n(1-p/s)<\alpha<n(p-1)$ and  $\alpha\leq 0$.
  These  can be regarded
as some substitutes for the type $(H^p,H^p)$ and $(H^p,L^p)$  for
$p$ in the remain range of $(0,1]$. These  estimates are also new
for the operators which are of type $(H^p,H^p)$ and type
$(H^p,L^p)$. And these   can also be regarded as an extension of the
type $({L}^{p} _{|x|^{\alpha}}, {L}^{p} _{|x|^{\alpha}})$ from
$1<p<\infty$ to $0<p\leq 1$.

\par In fact, we will prove that (1.12) and (1.13) hold
for the  sublinear  operator  $H$
 which satisfies the size
condition
$$
|Hf(x)| \leq C \|f\|_{L^{1}}/|x-x_0|^n,\eqno(1.14)
$$
when supp $f\subseteq B(x_0,2^k)$ and $|x-x_0| \geq 2^{k+1}$ with
$k\in {\bf Z} $, and the countable sublinear property
$$
|Hf|\leq \sum|\lambda_j||Ha_j| \eqno(1.15)
$$
 for   $f=\sum\lambda_ja_j \in \dot{B}L^{p,s}
_{|x|^{\alpha}}$, where each $a_j$ is a   $(  p,s, \alpha)-$block.

\par  (1.14) and (1.15) are satisfied by a lot of operators mentioned above, see Section 7.

\par The paper is organized as follows. In section 2, we introduce
the spaces $BL^{p,s}_{|x|^{\alpha}} $, and  list some of their
properties, including the relationship $BL^{p,s}_{|x|^{\alpha}}
\subset L^p_{|x|^\alpha}$. Such relationship and the boundedness of
operators on $BL^{p,s}_{|x|^{\alpha}} $  imply the boundedness from
$BL^{p,s}_{|x|^{\alpha}} $  into $L^p_{|x|^\alpha}$. In section 3 we
establish a molecular theorem for $BL^{p,s}_{|x|^{\alpha}}$, which
is a foundation for the estimates of operators on
$BL^{p,s}_{|x|^{\alpha}} $. In section 4, we prove the boundedness
of
 Hardy-Littlewood maximal operator on $BL^{p,s}_{|x|^{\alpha}} $,
 and this is also used in the estimates of
Bochner-Riesz means and the maximal spherical means. In sections 5
and    6, we establish the estimates of Bochner-Riesz means and the
maximal spherical means respectively. In section 7 we obtain some
estimates of some other sublinear operators, including some maximal
and truncated singular integral operators,  the maximal Carleson
operators, etc. In section 8 we prove two sharp results relating to
the quasinorm of $BL^{p,s}_{|x|^{\alpha}}$ and Carleson operator.

\section*
 {\bf  2.~~  $
 BL^{p,s}_{|x|^{\alpha}}$    and their  basic properties   }
{\bf 2.1 Definition}

\par Throughout this paper,
 we define
 $B(x_0,r)=\{x\in {\bf R}^n: |x-x_0| \leq r \},$ and denote
$\|g\|_{L_w^q}=\left( \int_{{\bf R}^n} |g(x)|^qw(x)dx
\right)^{1/q},\|g\|_{L^q }=\|g\|_{L_1^{q}}, q'=q/(q-1){\rm ~for~ }
q\geq 1, \bar{p} = {\rm min} \{p, 1\} {\rm ~ for ~} 0<p<\infty,
\mu_\alpha(x)=|x|^\alpha $ and $ \mu_\alpha(B)=\int_B
\mu_\alpha(x)dx$. And $\chi _{E}$ is the characteristic function of
set $E$.

\par {\bf Definition 2.1} ~~Let $- \infty < \alpha < \infty ,  0< p < \infty $ and $
0< s \leq \infty.$
\par A. A function $a(x)$ is said to be a
 $ (p, s, \alpha)$-block on ${\bf R}^n$ (centered at  $x_0$), if
\par (i)~~~~supp $a\subseteq B(x_0,r)\subset {\bf R}^n, r>0,$
\par (ii)~~~~$\|a\|_{L^{s}   }\leq |B(x_0,r)|^{-\alpha/pn-1/p +1/s};$

\par B. A function $a(x)$ is said to be a
 $ (p, s, \alpha)$-block of restrict I-type on ${\bf R}^n$ (centered at  $x_0$), if
(i) is replaced by

\par (i)$'$~~~~supp $a\subseteq B(x_0,r)\subset {\bf R}^n, r>1;$

\par C. A function $a(x)$ is said to be a
 $ (p, s, \alpha)$-block of restrict II-type on ${\bf R}^n$ (centered at  $x_0$), if
(i) is replaced by

\par (i)$''$~~~~supp $a\subseteq B(x_0,r)\subset {\bf R}^n, 0<r\leq 1.$

\par {\bf Definition 2.2}~~ Let $- \infty < \alpha < \infty , 0<p< \infty
$ and $0< s \leq \infty.$
  The   function spaces
$ BL^{p,s}_{|x|^{\alpha}}({\bf R}^n)$ are defined as
\begin{eqnarray*}
 BL^{p,s}_{|x|^{\alpha}}({\bf R}^n)=\{  f&: &
 f=\sum\limits_{k=-\infty}^{\infty} \lambda _ka_k ,
 \\&&\textrm { where
 each $a_k$ is a $(p,s, \alpha)$-block on ${\bf R}^n$,}
 \sum\limits_{k=-\infty}^{\infty} |\lambda _k|^{\bar{p}} <+ \infty \},
 \end{eqnarray*}
 here,  the
 "convergence"  means  a.e. convergence.
 Moreover, we define a quasinorm on $ BL^{p,s}_{|x|^{\alpha}}({\bf
 R}^n)$ by
$$\|f\|_{BL^{p,s}_{|x|^{\alpha}}({\bf R}^n)}=
 \inf \left(\sum\limits_{k=-\infty}^{\infty}|\lambda _k|^{\bar{p}}\right)^{1/{\bar{p}}},$$
where the infimum is taken over all the decompositions of $f$ as
above.
\par Similarly, we can also define the non-homogeneous function spaces $ \dot {B}L^{p,s}_{|x|^{\alpha}}({\bf R}^n)$
  with $ (p, s, \alpha)$-blocks of restrict I-type, and the non-homogeneous function spaces $\ddot{ B}L^{p,s}_{|x|^{\alpha}}({\bf R}^n)$
  with $ (p, s, \alpha)$-blocks of restrict II-type.

\par For simplicity, we will omit the notes $({\bf R}^n)$ below.
\par The balls in Definition 2.1 and Definition 2.2 can be replaced by cubes, and
$$
{\rm ball}- BL^{p,s}_{|x|^{\alpha}} = {\rm
cube}-BL^{p,s}_{|x|^{\alpha}}.
$$

\par We will use these two definitions in this paper simultaneously, even in the
 proof of a theorem.

\par {\bf Comments }

 \par $ BL^{p,s}_{|x|^{\alpha}}$ are similar to the classical   Hardy spaces $H^p$,
the weighted Hardy spaces $H^p_{|x|^{\alpha}}$,  and the block
spaces $B_{s}(I)$, $I=(-1/2,1/2)$,  but there are some differences
between them.


\par (1) As we see, $ BL^{p,s}_{|x|^{\alpha}}$ is the
space generated by blocks.   $ BL^{p,s}_{|x|^{\alpha}}$ is an
extension of the classical   Hardy spaces $H^p$ with $0<p\leq 1$,
since
 $  {H}^{p} $ has an
atom decomposition ( see
 \cite{CW} for  details).
 We see that $ BL^{p,s}_{|x|^{\alpha}}$ is different
from $H^{p}$, since $H^p$- atom has the cancelation properties while
$(p,s, \alpha)$ - block has not.

Generally,  $H^p$ is regarded as a subspace of $L^p$, but here we
use $ BL^{p,s}_{|x|^{\alpha}}$ as a subspace of the weighted
Lebesgue space $L^p_{|x|^\alpha}$ for $0<p\leq 1$.
\par (2) It is well-known that the weighted Hardy space $  {H}^{p}_{|x|^{\alpha}}
$ is a subspace of $L^p_{|x|^\alpha}$, and it plays an analogous
role as $H^p$.
 $  {H}^{p}_{|x|^{\alpha}} $ has an
atom decomposition  similar to $H^p$
 (see   \cite{Gc} for   details).
       We see that the size condition and cancelation condition of ${H}^{p}_{|x|^{\alpha}} $ atom in \cite{Gc}
   are different
from  $(p,s, \alpha)$ - blocks.
\par The weighted Hardy space $  {H}^{p}_{|x|^{\alpha}} $ shares
an alternative  atomic decomposition  (see \cite{ST} for details).
Its atom has the cancellation properties but $(p,s, \alpha)$ - block
has not.

\par (3)
 $
BL^{p,s}_{|x|^{\alpha}}$ is an extension of the block space
$B_{s}(I)$ with $1<s\leq \infty$ and  $I=(-1/2,1/2)$ introduced by
Taibleson and Weiss in the study of the convergence of the Fourier
series in \cite{TW1}.
 The norm of the function $f$ in $B_{s}(I)$ is defined as
  $$ \|f\|_{B_{s}(I)}=\inf
\sum\limits_{k=-\infty}^{\infty} |\lambda _k|
 \left(1+{\rm log } \left(
\frac{1}{|\lambda _k|}\sum\limits_{i=-\infty}^{\infty} |\lambda
_i|\right)\right)    ,  $$ where the infimum is taken over all the
decompositions of $f$ in $B_{s}(I)$. This is different from that of
$ BL^{p,s}_{|x|^{\alpha}}$.
 Some basic properties of
$B_{s}(I)$ and its applications to harmonic analysis have been
studied by many authors (see \cite{TW1}, \cite{LTW1},
etc), and  $B_{s}(I)$ was used as a subspace of $L^1$ there.

\par  $
B^{\phi}_{s}(S^{n-1}) $, where  $1\leq s \leq \infty $ and $ \phi$
is a nonnegative function on ${\bf R}^+$,  is an extension of
$B_{s}(I)$  to higher dimensions. Like in $B_{s}(I)$, The norm of
the function $f$ in $ B^{\phi}_{s}(S^{n-1}) $ is different from that
of $ BL^{p,s}_{|x|^{\alpha}}$ (see \cite{TW1} for details).
   $
B^{\phi}_{s}(S^{n-1}) $ is useful in the study of boundedness of
some generalized singular integrals with rough kernels
(see, for example,
\cite{LTW1})
.

\par {\bf   Remarks}
\par (1) The space $ BL^{p,s}_{|x|^{0}} $
for $0<p<s\leq \infty $ and $\alpha=0$ was introduced in \cite{Lo}
and the space $ BL^{1,s}_{|x|^{\alpha}} $ for $1\leq s<\infty$ and $
-n<\alpha$ in \cite{BRV}. They are $h_{p,s}$ in \cite {Lo} and
${\mathcal{B}}_{s,\alpha+n}$ in \cite {BRV}, respectively.
\par (2) $ BL^{1,1}_{|x|^{0}} =L^1 $. See
\cite{S3}. In fact,   $ L^{p}= BL^{p,s}_{|x|^{0}}  $ for $0<p\leq 1$
and $0< p< s \leq \infty$ (see
 Corollary 2.2 below), and   because of this we use the notes $BL^{p,s}_{|x|^{\alpha}}$.
\\{\bf 2.2 Properties}
\par Next, let us give a number of properties of $
BL^{p,s}_{|x|^{\alpha}} $. Among these properties, Theorem 2.1 will
be used many times in the estimates of operators, but the others
will seldom or even never be used in the following sections.
However, in order to help us better understand the spaces $
BL^{p,s}_{|x|^{\alpha}} $, we show all these properties and prove
them. We also proved in \cite{Longs} that the dual of
$BL^{p,s}_{|x|^{\alpha}} $ are the classical Morrey spaces
$M^{s'}_{1-1/p-\alpha/np}$ for $1<s<\infty, 0<p\leq s$ and
$n(p/s-1)<\alpha<n(p-1)$.
\\ {\bf  2.2.1.~~ Some relationships between
$BL^{p,s}_{|x|^{\alpha}} $  and $L^{p}_{|x|^{\alpha}},  L^p$ and
$H^p$  }
\\ {\bf A. Between $BL^{p,s}_{|x|^{\alpha}}$   and
 $L^{p}_{|x|^{\alpha}}$}

\par  {\bf Theorem 2.1 }  Let $ 0<p< s \leq \infty$ and $ -n(1-p/s) < \alpha \leq 0 $.
Then
$$BL^{p,s}_{|x|^{\alpha}} \subset L^{p}_{|x|^{\alpha}}.$$
\par {\bf Proof }   Let $b$ be a $(p,s,\alpha)$-block with supp $b \subset Q$.
\par When $0<s<+\infty,$ using H\"{o}lder inequality for the index $p/s$,
and noticing $\alpha \leq 0$ and $-n(1-p/s) < \alpha$, we have
\begin{eqnarray*}
\|b\|^p_{L^{p}_{|x|^{\alpha}}}
  &\leq & \|b\|^p_{L^{s}} \left(\int_{Q
}|x|^{\alpha s/(s-p)}dx\right)^{(s-p)/s}
 \\ &\leq & \|b\|^p_{L^{s}}
\left(\int_{\widetilde{Q} }|x|^{\alpha s/(s-p)}dx\right)^{(s-p)/s}
 \\ &\leq & C|Q|^{(-\frac{\alpha}{p n}-\frac{1}{p}+\frac{1}{s})p}
 |Q|^{ (\frac{\alpha}{p n}+\frac{1}{p}-\frac{1}{s})p}
=C,
\end {eqnarray*}
 where $\widetilde{Q}$ denotes the balls
obtained by translation of $Q$ and  centered at $0$.
\par When $s=+\infty,$ we have
\begin{eqnarray*}
\|b\|^p_{L^{p}_{|x|^{\alpha}}}
  \leq  \|b\|^p_{L^{\infty}}  \int_{Q
}|x|^{\alpha }dx \leq  \|b\|^p_{L^{\infty}}
 \int_{\widetilde{Q} }|x|^{\alpha  }dx
 \leq  C|Q|^{(-\frac{\alpha}{p n}-\frac{1}{p} )p} |Q|^{ \frac{\alpha}{  n}+1}
=C.
\end {eqnarray*}

\par $C$ are independent of $b$.

Let $f  \in BL^{p,s}_{|x|^{\alpha}}$. For any $\varepsilon>0,$ there
exists a decomposition of   $f=\sum_{j=1}^{\infty} \lambda _i a_j ,$
such that $\sum_{j }^{ }|\lambda_j|^{\bar{p}}\leq
\|f\|^{\bar{p}}_{BL^{p,s}_{|x|^{\alpha}}}+\varepsilon <\infty$.
Then, we have
$$
\|f\|^p_{L^{p}_{|x|^{\alpha}}}
 \leq \sum_{j }^{
}|\lambda_j|^p \|a_j\|^p_{L^{p}_{|x|^{\alpha}}}
 \leq C\sum_{j }^{
}|\lambda_j|^p<\infty
 $$
for $0<p\leq 1$, and
$$
\|f\| _{L^{p}_{|x|^{\alpha}}}
 \leq \sum_{j }^{
}|\lambda_j|  \|a_j\| _{L^{p}_{|x|^{\alpha}}}
 \leq C\sum_{j }^{
}|\lambda_j|<\infty
 $$
for $1<p<\infty$. Thus, $f\in L^{p}_{|x|^{\alpha}}.$

\par  {\bf Theorem 2.2 } Let $  0<s \leq \infty ,0<p\leq 1 $  and $ 0\leq \alpha < \infty $.
Then
$$BL^{p,s}_{|x|^{\alpha}} \supset L^{p}_{|x|^{\alpha}}.$$

\par {\bf Proof } We use the definition of cube-$BL^{p,s}_{|x|^{\alpha}}$
which equals to ball-$BL^{p,s}_{|x|^{\alpha}}$. Let $f\in
L^{p}_{|x|^{\alpha}}$. We first consider the simple functions
$f=\sum_{j }^{ }c_j\chi_{Q_j}$, where $\{Q_j\} $ is an at most
countable set  of cubes whose interiors are mutually disjoint. Then
$$f=\sum_{j }^{ }c_j|Q_j|^{1/p+\alpha/np}|Q_j|^{-1/p-\alpha/np}\chi_{Q_j}=\sum_{j }^{ }m_jb_j,$$
where $m_j=c_j|Q_j|^{1/p+\alpha/np}$ and
$b_j=|Q_j|^{-1/p-\alpha/np}\chi_{Q_j}$. And noticing $\alpha \geq 0,
$
$$
\|f\|^p_{L^{p}_{|x|^{\alpha}}}=\sum_{j }^{
}|c_j|^p\int_{Q_j}|x|^{\alpha}dx \geq \sum_{j }^{
}|c_j|^p\int_{\widetilde{Q_j}}|x|^{\alpha}dx =\sum_{j }^{
}|c_j|^p|Q_j|^{1 +\alpha/n } =\sum_{j }^{ }|m_j|^p,\eqno(2.1)
$$
where $\widetilde{Q_j}$ denotes the cubes obtained by translation of
$Q_j$ and  centered at $0$. Thus $f\in BL^{p,s}_{|x|^{\alpha}}.$

\par Next, let $\textbf{S}$ be  the set of all  cube-type simple
functions with cube's sides which are parallel to the axes. We claim
that: $\textbf{S}$ is dense in $ L^{p}_{|x|^{\alpha}}.$ In fact, let
$f=\sum_{j=1}^{N}c_j\chi_{E_j}$, where  ${E_j} $ are bounded closed
sets which are mutually disjoint and $|c_j|\leq M $. It is easy to
see that  the set of all these type functions is dense in $
L^{p}_{|x|^{\alpha}}$. And there exist open sets
$\widetilde{O_j}\supset E_j,$ and $\widetilde{O_j}$ are mutually
disjoint. At the same time, for any $\varepsilon >0,$ since $E_j$
are closed set, there exist open sets
$\widetilde{\widetilde{O_j}}\supset E_j$ and
$\mu_\alpha(\widetilde{\widetilde{O_j}} \setminus E_j)<
\frac{1}{2}\varepsilon M^{-p}N^{-1}.$ Let $O_j= \widetilde{O_j}
\bigcap \widetilde{\widetilde{O_j}} ,$
 then the open sets $O_j$ are  mutually disjoint, $O_j\supset E_j$,
 and $\mu_\alpha(O_j \setminus E_j)<\frac{1}{2} \varepsilon
M^{-p}N^{-1},$ where $\mu_\alpha(x)=|x|^\alpha$. Since each $O_j$ is
an open set, then, by Whitney decomposition theorem,
$O_j=\bigcup_{i=1}^{\infty}Q_i^{(j)}, $ where $Q_i^{(j)}$ are cubes
whose interiors are mutually disjoint and whose sides are parallel
to the axes (see \cite {St61970}), and
$\sum_{i=1}^{\infty}\mu_\alpha(Q_i^{(j)})=\mu_\alpha(\bigcup_{i=1}^{\infty}Q_i^{(j)})
=\mu_\alpha(O_j)=\mu_\alpha(E_j)+\mu_\alpha(O_j\setminus E_j)<
\mu_\alpha(E_j)+\frac{1}{2} \varepsilon M^{-p}N^{-1}<\infty.$ Then
there exists $k_j$ such that
$\sum_{i=k_j}^{\infty}\mu_\alpha(Q_i^{(j)})< \frac{1}{2} \varepsilon
M^{-p}N^{-1}.$ Thus
$\mu_\alpha(\bigcup_{i=1}^{k_j-1}Q_i^{(j)}\setminus E_j)\leq
\mu_\alpha(O_j \setminus E_j)<\frac{1}{2} \varepsilon M^{-p}N^{-1}$
and $ \mu_\alpha(E_j\setminus \bigcup_{i=1}^{k_j-1}Q_i^{(j)})\leq
\mu_\alpha(\bigcup_{i=k_j}^{\infty}Q_i^{(j)})<\frac{1}{2}
\varepsilon M^{-p}N^{-1} $.
 Therefore,
\begin{eqnarray*}
\|\sum_{j=1}^{N}c_j\left(\chi_{E_j}-\sum_{i=1}^{k_j-1}\chi_{Q_i^{(j)}}\right)\|^p_{L^{p}_{|x|^{\alpha}}}
= \|\sum_{j=1}^{N}c_j\left(\chi_{E_j}-\chi_{\bigcup_{i=1}^{k_j-1}
Q_i^{(j)}}\right)\|^p_{L^{p}_{|x|^{\alpha}}} \\\leq
\sum_{j=1}^{N}|c_j|^p\left(
\mu_\alpha(\bigcup_{i=1}^{k_j-1}Q_i^{(j)}\setminus
E_j)+\mu_\alpha(E_j\setminus
\bigcup_{i=1}^{k_j-1}Q_i^{(j)})\right)<\varepsilon.
\end{eqnarray*}
Thus, the cube type simple function
$\sum_{j=1}^{N}\sum_{i=1}^{k_j-1}c_j\chi_{Q_i^{(j)}}$ approaches to
$f$ in $L^{p}_{|x|^{\alpha}}$. We finish the claim.

\par Now let $f\in L^{p}_{|x|^{\alpha}}$, and $ \{f_n\}\subset \textbf{S}$
approaches to $f$ in $L^{p}_{|x|^{\alpha}}.$ For any
$\varepsilon>0,$ choose $\{n_j\}$ such that
$$
\|f_{n_0}\|^p_{L^{p}_{|x|^{\alpha}}}\leq
\|f\|^p_{L^{p}_{|x|^{\alpha}}} +\frac{\varepsilon}{2},
\|f_{n_k}-f_{n_{k-1}}\|^p_{L^{p}_{|x|^{\alpha}}}\leq
\frac{\varepsilon}{2^k}.
$$
It is easy to see that any cube whose sidelength  may not be equal
can be expressed as a union of at most countable cubes with equal
sidelength. Then, a
 difference of the two functions in $\textbf{S}$ is a linear combination
 of  the characteristic functions of
 at most
 countable cubes, whose interiors are mutually disjoint. So
$$f_{n_k}-f_{n_{k-1}} =\sum_{j=1 }^{\infty }c_j^{(k)}\chi_{Q_j^{(k)}}
=\sum_{j=1 }^{ \infty}m_j^{(k)}b_j^{(k)}, k=0,1,\cdots,
$$
where $f_{n_{-1}}=0, m_j^{(k)}=c_j^{(k)}|Q_j^{(k)}|^{1/p+\alpha/np}
$ and $  b_j^{(k)}=|Q_j^{(k)}|^{-1/p-\alpha/np}\chi_{Q_j^{(k)}}. $
By (2.1), we have
$$\sum_{j=1 }^{ \infty}|m_j^{(0)}|^p \leq \|f_{n_0}\|^p_{L^{p}_{|x|^{\alpha}}}\leq
\|f\|^p_{L^{p}_{|x|^{\alpha}}} +\frac{\varepsilon}{2},$$

$$\sum_{j=1 }^{ \infty}|m_j^{(k)}|^p \leq \|f_{n_k}-f_{n_{k-1}}\|^p_{L^{p}_{|x|^{\alpha}}}\leq
\frac{\varepsilon}{2^k}.$$ Thus,
$$
f=\sum_{k=1 }^{ \infty} (f_{n_k}-f_{n_{k-1}}) + f_{n_0}= \sum_{k=0
}^{ \infty} \sum_{j=1 }^{ \infty} m_j^{(k)}b_j^{(k)}
$$
in $L^{p}_{|x|^{\alpha}}$-norm, and
$$ \sum_ {k=0
}^{ \infty} \sum_{j=1 }^{ \infty}|m_j^{(k)}|^p \leq
\|f\|^p_{L^{p}_{|x|^{\alpha}}} +\frac{\varepsilon}{2} + \sum_{k=1
}^{ \infty} \frac{\varepsilon}{2^k}= \|f\|^p_{L^{p}_{|x|^{\alpha}}}
+\frac{3}{2} \varepsilon,
$$
then  $\|f\|_{BL^{p,s}_{|x|^{\alpha}}}\leq \|f\|_{ {L}^{p
}_{|x|^{\alpha}}}$ follows.

\par{\bf Remark 2.1} For $\alpha=0$ Theorem 2.2 is known, the case
$p=1$ can be found in \cite{S3}, and $0<p\leq 1$ in \cite{Lo}.

\par Theorems 2.1 and   2.2 give

\par {\bf Corollary 2.1}  Let $0<p\leq 1 $ and $ p< s \leq \infty$. We have
$$ L^{p}=BL^{p,s}_{|x|^{0}}. $$
\\
{\bf B.   Between $BL^{p,s}_{|x|^{\alpha}}$   and $L^{p}_{}$}

\par {\bf Theorem 2.3} Let $0<p\leq s\leq\infty$. We have,
\par A) if $ n(p/s-1) \leq \alpha\leq n(p-1)$, or $0<p\leq 1,n(p-1)\leq \alpha \leq
0$,
$$BL^{p,s}_{|x|^{\alpha}}\subset L^{\frac{np}{n+\alpha}},  $$

 \par B) if $-n<\alpha \leq n(p/s-1)$,
$$ L^s\subset BL^{p,s}_{|x|^{\alpha}}, $$

   \par C) if $ \alpha = n(p/s-1)$,
$$ L^s = BL^{p,s}_{|x|^{\alpha}} . $$

\par {\bf Proof }
 A). Let $f\in BL^{p,s}_{|x|^{\alpha}},$ then, $f(x)=\sum
\lambda_ja_j(x)$, where each $a_j$ is a $(p,s,\alpha)$-block with
supp $a_j\subseteq Q_j$ (ball or cube), and $\sum
|\lambda_j|^{\overline{p}}<\infty.$ Let $0< r \leq s\leq\infty. $
\par When $1\leq r \leq s \leq \infty,$
by Minkowski inequality, we have
\begin{eqnarray*}
\|f\|_{L^{r}} &\leq  &\sum |\lambda_j|\| a_j\|_{L^{r}} \leq  \sum
|\lambda_j|\| a_j\|_{L^{s}} |Q_j|^{1/r-1/s}
\\&\leq &C  \sum
|\lambda_j| |Q_j|^{1/r-1/s} |Q_j|^{-\alpha/pn-1/p+1/s}
\\ &=&C  \sum
|\lambda_j|  |Q_j|^{-\alpha/pn-1/p+1/r}.
\end{eqnarray*}
For $n(p/s-1)\leq \alpha \leq n(p-1),$ take $r=\frac{np}{\alpha+n},$
then  $ 1\leq r\leq s .$ Thus
\begin{eqnarray*}
\|f\|_{L^{\frac{np}{\alpha+n}}} \leq  \sum |\lambda_j| \leq C ( \sum
|\lambda_j|^{\overline{p} } )^{1/\overline{p}}<\infty.
\end{eqnarray*}
\par When $0<   r\leq 1,$ and $r\leq s$, by H\"{o}lder inequality,   we have
\begin{eqnarray*}
\|f\|^r_{L^{r}} &\leq  &\sum |\lambda_j^r|\| a_j\|^r_{L^{r}} \leq
\sum |\lambda_j|^r\| a_j\|^r_{L^{s}} |Q_j|^{(1/r-1/s)r}
\\&\leq &C  \sum
|\lambda_j|^r |Q_j|^{(1/r-1/s)r} |Q_j|^{(-\alpha/pn-1/p+1/s)r}
\\ &=&C  \sum
|\lambda_j|^r  |Q_j|^{(-\alpha/pn-1/p+1/r)r}.
\end{eqnarray*}
For $   n(p-1)\leq \alpha \leq 0,$ take $r=\frac{np}{\alpha+n},$
then $ 0< p\leq r\leq 1,$ and $r\leq s.$ We then have
\begin{eqnarray*}
\|f\|_{L^{\frac{np}{\alpha+n}}} \leq  C ( \sum |\lambda_j|^r )^{1/r}
\leq C ( \sum |\lambda_j|^p )^{1/p}<\infty.
\end{eqnarray*}
\par B). Assume first that $\alpha <n(p/s-1).$   Let $f\in L^{s},$ write
$$f=\sum_{k=1}^{\infty} |B_k|^{\frac{\alpha}{pn}+\frac{1}{p}-\frac{1}{s}}\|f\|_{L^{s}}b_k,$$
where
$$b_1=\frac{f\chi_{B_1}}{\|f\|_{L^{s}}}, ~~
b_k=\frac{f\chi_{C_k}}{|B_k|^{\frac{\alpha}{pn}+\frac{1}{p}-\frac{1}{s}}\|f\|_{L^{s}}},~~k=2,3,\cdots,$$
here $B_k=B(0,2^k)$ and $ C_k=B_k\setminus B_{k-1}, k=1,2,\cdots.$
We see that all $b_k$ are $(p,s,\alpha)-$ blocks,  and
\begin{eqnarray*}
\|f\|_{ BL^{p,s}_{|x|^{\alpha}}} \leq  \left(\sum _{k=1}^{\infty}
|B_k|^{(\frac{\alpha}{np}+\frac{1}{p}-\frac{1}{s})\bar{p}}\|f\|^{\bar{p}}_{L^{s}}\right)^{1/\bar{p}}\leq
\left(\sum _{k=1}^{\infty}
2^{nk(\frac{\alpha}{np}+\frac{1}{p}-\frac{1}{s})\bar{p}}\right)^{1/\bar{p}}\|f\|_{L^{s}}=C\|f\|_{L^{s}}
\end{eqnarray*}
since $\frac{\alpha}{np}+\frac{1}{p}-\frac{1}{s}<0.$
\par Now assume  that $\alpha =n(p/s-1)$  (i.e.
  $s=np/(\alpha+n)$),
and $f\in L^s$.

\par  Let us first consider the case of $1\leq s <\infty$. Since $f\chi_{B_k}$
is a Cauchy sequence in $L^s$ then there exists
$\{n_k\}_{k=1}^{\infty}$ such that
$$\|f\chi_{B_{n_{k+1}}}-f\chi_{B_{n_{k}}}\|_{L^s}<2^{-k}\|f\|_{L^s}.\eqno(2.2)$$
Let
$$f=f\chi_{B_{n_{1}}}+\sum_{k=1}^{\infty}
(f\chi_{B_{n_{k+1}}}-f\chi_{B_{n_{k}}})=
\sum_{k=0}^{\infty}2^{-k}\|f\|_{L^s}b_k,$$ where
$$b_0=\frac{f\chi_{B_{n_1}}}{\|f\|_{L^{s}}}, ~~
b_k=\frac{f(\chi_{B_{n_{k+1}}}-\chi_{B_{n_{k}}})}{2^{-k}\|f\|_{L^{s}}},~~k=1,2,3,\cdots.$$
Noticing that $\frac{\alpha}{pn}+\frac{1}{p}-\frac{1}{s}=0$ and
(2.2), we see that $b_k$ are $(p,s,\alpha)-$blocks. This shows that
$f\in BL^{p,s}_{|x|^{\alpha}}$  and $\|f\|_{
BL^{p,s}_{|x|^{\alpha}}} \leq C\|f\|_{L^{s}}.$

\par For the case of $0<s<1$, since  $f\in L^s$, we have
$|f|^s\in L^1$, and it follows that  $|f|^s\chi_{B_k}$ is a Cauchy
sequence in $L^1$. Then we can find $\{n_k\}_{k=1}^{\infty}$, and
replace (2.2) with
$$\||f|^s\chi_{B_{n_{k+1}}}-|f|^s\chi_{B_{n_{k}}}\|_{L^1}<2^{-sk}\||f|^s\|_{L^1}.$$
And  the remain statement is the same as  the case of $1\leq s
<\infty$.

\par Then we can obtain C) from A) and B). Thus, we finish the proof of
Theorem 2.3.
\\ {\bf C.   Between $BL^{p,s}_{|x|^{\alpha}}$ and $H^{p}$}
\par {\bf Theorem 2.4}
  Let $0< s\leq \infty, 0<p<\infty$ and $n(p-1)\leq \alpha < \infty $.
Then
$$H^{\frac{np}{n+\alpha}}\subset BL^{p,s}_{|x|^{\alpha}}
. $$
\par {\bf Proof} This can be seen easily from the atom
decomposition theory of Hardy space $H^{\frac{np}{n+\alpha}}$ (see
for example \cite{CW}).
 \\{\bf D.   Between
$BL^{p,s}_{|x|^{\alpha}}$  and $BL^{p,s}_{|x|^{\alpha}}$ }

\par It is easy to check that

\par {\bf Theorem 2.5}
  Let $0<p<\infty, 0< s_1\leq s_2\leq \infty$ and $- \infty< \alpha < \infty $.
Then
$$BL^{p,s_2}_{|x|^{\alpha}} \subset
BL^{p,s_1}_{|x|^{\alpha}} .
$$
\\{\bf 2.2.2  Completeness}
\par Let $0<s\leq \infty,0<p\leq \infty$ and $-\infty <\alpha< \infty$. It
is easy to see that,  (i) if $f=0,$ then
$\|f\|_{BL^{p,s}_{|x|^{\alpha}}}=0,$   (ii) $\|\cdot
\|_{BL^{p,s}_{|x|^{\alpha}}}$
 are positive homogeneous,
and (iii) $\|\cdot\|_{BL^{p,s}_{|x|^{\alpha}}}$ is subadditive. Now
we prove (iii).
 Let $f, g \in BL^{p,s}_{|x|^{\alpha}}.$ For
any $\varepsilon >0$, there exist the block -decompositions of $f$
and $g$: $f=\sum m_kb_k $ and $g=\sum n_jc_j $ such that
$$\|f\|_{BL^{p,s}_{|x|^{\alpha}}}^{\bar{p}} \geq
 \sum_k|m_k|^{\bar{p}}  -\varepsilon, ~~~~{\rm and}
~~~~\|g\|_{BL^{p,s}_{|x|^{\alpha}}}^{\bar{p}} \geq
 \sum_j|n_j|^{\bar{p}}   -\varepsilon. $$ Then
$$\|f+g\|_{BL^{p,s}_{|x|^{\alpha}}}^{\bar{p}} \leq
 \sum_k|m_k|^{\bar{p}} + \sum_j|n_j|^{\bar{p}} \leq
 \|f\|_{BL^{p,s}_{|x|^{\alpha}}}^{\bar{p}}+\|g\|_{BL^{p,s}_{|x|^{\alpha}}}^{\bar{p}}+ 2\varepsilon. $$
Since $\varepsilon$ is arbitrary, (iii) follows.

\par For $0<s\leq \infty,0<p< s,-n(1-p/s) <\alpha\leq 0,$
by Theorem 2.1, we see that
\par (iv) if   $\|f\|_{BL^{p,s}_{|x|^{\alpha}}}=0,$
then $f=0$ $\mu_\alpha$-a.e..

\par So $
BL^{p,s}_{|x|^{\alpha}} $  are normed-spaces for  $0<s\leq
\infty,0<p< s$ and $-n (1-p/s)<\alpha\leq 0,$ and they are also
quasi-normed spaces in the sense of (i)$\sim$ (iii) for $0<s\leq
\infty,0<p\leq \infty$ and $-\infty <\alpha< \infty$.

\par {\bf Theorem 2.6 }  Let  $0<s\leq \infty,0<p< s$ and $-n(1-p/s)
<\alpha \leq 0$. Then, $ BL^{p,s}_{|x|^{\alpha}} $   are complete.
\par {\bf Proof}  Let $\{u_n\}$ be a Cauchy sequence in $ BL^{p,s}_{|x|^{\alpha}}. $
 For any $\varepsilon >0$, there exists a subsequence $\{u_{n_j}\} $
  of $\{u_n\}$ such that
  $$\|u_{n_{j+1}}-u_{n_j}\|_{BL^{p,s}_{|x|^{\alpha}}} \leq \frac{\varepsilon}{2^j}, ~~~~j=1,2,\cdots, $$
and
$$\|u_{n }-u_{n_1}\|_{BL^{p,s}_{|x|^{\alpha}}} \leq
\varepsilon $$ for $n>n_1.$
 Set
$\bar{u}=\sum_{j=1}^{\infty}(u_{n_{j+1}}-u_{n_j}).$ $\bar{u}$  can
be represented as a linear combination of $(p,s,\alpha)-$ blocks
since $u_{n_{j+1}}-u_{n_j}\in BL^{p,s}_{|x|^{\alpha}}$. And, by
Theorem 2.1, we see that
$$\|\bar{u}\|_{L^{p}_{|x|^{\alpha}}}\leq \|\bar{u}\|_{BL^{p,s}_{|x|^{\alpha}}}
\leq
\sum_{j=1}^{\infty}\|(u_{n_{j+1}}-u_{n_j})\|_{BL^{p,s}_{|x|^{\alpha}}}
\leq \sum_{j=1}^{\infty} \frac{\varepsilon}{2^j}=\varepsilon,$$ for
 $0<s\leq \infty, 0<p< s$ and $-n(1-p/s)
<\alpha  \leq 0$. Then
$\bar{u}=\sum_{j=1}^{\infty}(u_{n_{j+1}}-u_{n_j})$ converges
$\mu_\alpha$-a.e., and $\bar{u}\in BL^{p,s}_{|x|^{\alpha}}.$ Let
$u=\bar{u}+u_{n_1}$, then $u\in BL^{p,s}_{|x|^{\alpha}},$ and
 $$\|u_{n }-u_{ }\|_{BL^{p,s}_{|x|^{\alpha}}} \leq \|u_{n }-u_{n_1 }\|_{BL^{p,s}_{|x|^{\alpha}}} +
\|\bar{u} \|_{BL^{p,s}_{|x|^{\alpha}}}
 \leq 2
\varepsilon $$ for $n>n_1.$ Thus, we finish the proof for the
homogeneous case. And the proof for the non-homogeneous case is
similar.
 \\{\bf 2.2.3 Density}

 \par We denote by $C^m, 0\leq m\leq \infty,  $ the set of $m$-times
 continuously differentiable functions on ${\bf R}^n$,  by  $C_c^m $  the set of all functions in
 $C^m$ having a compact support, and by $\mathcal{S}$ Schwartz spaces.
   Let
\begin{eqnarray*}
\varphi(x)=\left\{\begin{array}{ll}Ce^{1/{(|x|^2-1)}} &\textrm{if
$|x|<1$},\\
0 &\textrm{if $|x|\geq 1$}
\end{array}\right.
 \end{eqnarray*}
such that $\int_{{\bf R}^n}\varphi(x)dx=1$ and
$\varphi_t(x)=\frac{1}{t^n}\varphi\left(\frac{x}{t}\right). $

\par {\bf Theorem 2.7} Let   $0<p<\infty,1\leq s<\infty$ and $ -n<\alpha<\infty.$ Then $C^0_c$ is
dense in $BL^{p,s}_{|x|^{\alpha}}.$

\par{\bf Proof} Let $f\in BL^{p,s}_{|x|^{\alpha}}$, i.e. $f(x)=\sum_{k=1}^{\infty} \lambda_ka_k(x)
$ where each $a_k$ is a $(p,s,\alpha)$-block with supp $
a_k\subseteq Q_k$ and $ \sum_{k=1}^{\infty}
 |\lambda_k|^{\bar{p}}<\infty $. Then for any $\varepsilon>0,$ we can find
 an
 $i_0$ such that
$$
 \sum_{k=i_0+1}^{\infty}
 |\lambda_k|^{\bar{p}}<\varepsilon^{\bar{p}}. \eqno(2.3)
 $$
Let $f_{i_0}(x)=\sum_{k=1}^{i_0} \lambda_ka_k(x) $, supp
$f_{i_0}\subseteq Q (\supseteq\bigcup_{k=1}^{i_0}Q_k)$ and
$f_{i_0}\in L^s(Q).$ We can find $g_{i_0}\in C^0(Q)$ such that
$$
\|f_{i_0}-g_{i_0}\|_{L^s(Q)}<\frac{\varepsilon^{}}{2|2Q|^{\alpha/np+1/p-1/s}}.\eqno(2.4)
$$
\par Let
\begin{eqnarray*}
\tilde{g}_{i_0}(x)=\left\{\begin{array}{ll}g_{i_0}(x) &\textrm{if
$x\in Q $},\\
0 &\textrm{if $x\in Q^c$},
\end{array}\right.
 \end{eqnarray*}
and
$$\tilde{g}_{i_0}^t=\tilde{g}_{i_0}\ast\varphi_t.$$
Then, for $t$ small enough, we have
$$
\|\tilde{g}_{i_0}^t-\tilde{g}_{i_0}\|_{L^s}<\frac{\varepsilon^{}}{2|2Q|^{\alpha/np+1/p-1/s}},\eqno(2.5)
 $$
and supp $\tilde{g}_{i_0}^t\subseteq $ supp
$2\tilde{g}_{i_0}\subseteq2Q$. Thus, by (2.4) and (2.5),
\begin{eqnarray*}
\|f_{i_0}-\tilde{g}_{i_0}^t\|_{L^s}\leq
\|f_{i_0}-{\tilde{g}}_{i_0}\|_{L^s}+\|\tilde{g}_{i_0}^t-\tilde{g}_{i_0}\|_{L^s}
<\frac{\varepsilon^{}}{|2Q|^{\alpha/np+1/p-1/s}},
 \end{eqnarray*}
and supp $(f_{i_0}-\tilde{g}_{i_0}^t)\subseteq  2Q$. Let
$h_{i_0}(x)=\frac{1}{\varepsilon}(f_{i_0}(x)-\tilde{g}_{i_0}^t(x))$.
Then $h_{i_0}(x)$ is a $(p,s,\alpha)$-block, and it follows that
$f_{i_0}(x)-\tilde{g}_{i_0}^t(x)=\varepsilon h_{i_0}(x)\in
BL^{p,s}_{|x|^{\alpha}}$, and
$$
\|f_{i_0}-\tilde{g}_{i_0}^t\|_{BL^{p,s}_{|x|^{\alpha}}}\leq
\varepsilon.\eqno(2.6)
 $$
By (2.3) and (2.6), and the subadditive of
$\|\cdot\|_{BL^{p,s}_{|x|^{\alpha}}}$, we have
\begin{eqnarray*}
\|f-\tilde{g}_{i_0}^t\|_{BL^{p,s}_{|x|^{\alpha}}} &\leq &
\|f_{i_0}-\tilde{g}_{i_0}^t\|_{BL^{p,s}_{|x|^{\alpha}}}+
\|\sum_{k=i_0+1}^{\infty}
\lambda_ka_k(x)\|_{BL^{p,s}_{|x|^{\alpha}}}
 \\ &\leq&
 \varepsilon+\left(\sum_{k=i_0+1}^{\infty}
 |\lambda_k|^{\bar{p}}\right)^{1/\bar{p}}
 \leq 2\varepsilon^p.
 \end{eqnarray*}

\par {\bf Theorem 2.8} Let   $0<p<\infty,1\leq s<\infty$ and $ -n<\alpha<\infty.$
If $u(x)$ is a continuous function on ${\bf R}^n$ with compact
support, then $u_t(x)=u\ast\varphi_t(x)\rightarrow u(x)$  in
$BL^{p,s}_{|x|^{\alpha}}$   as $t\rightarrow 0.$

\par {\bf Proof} Let supp $u\subseteq K$ which is a compact set in ${\bf
R}^n$. Then  supp $(u_t-u)\subseteq K_1\subset B(x_0,r)=B, r>1$. And
for any $\varepsilon>0$ and $t$  small enough, $\max_{x\in {\bf R}^n
}|u_t(x)-u(x)|<\varepsilon$, we have
\begin{eqnarray*}
\|u_t-u\|_{L^s}\leq \varepsilon|B|^{1/s}=\varepsilon|B|^
{\alpha/np+1/p}|B|^ {-\alpha/np-1/p+1/s}.
 \end{eqnarray*}
So, $a(x)=\frac{1}{\varepsilon|B|^ {\alpha/np+1/p}}(u_t(x)-u(t))$ is
a $(p,s,\alpha)$-block, and  it follows that
$u_t(x)-u(t)=\varepsilon|B|^ {\alpha/np+1/p}a(x)\in
BL^{p,s}_{|x|^{\alpha}}$ and
\begin{eqnarray*}
\|u_t-u\|_{BL^{p,s}_{|x|^{\alpha}}}\leq  \varepsilon|B|^
{\alpha/np+1/p} \rightarrow 0, {\rm ~as~}t\rightarrow 0.
 \end{eqnarray*}

\par {\bf Theorem 2.9} Let   $0<p\leq 1,1< s<\infty$ and $
-n(1-p/s)<\alpha<n(p-1),$ or let $1<p<s<\infty$ and $
-n(1-p/s)<\alpha \leq 0.$
 If $u(x)\in BL^{p,s}_{|x|^{\alpha}} $,
then $u_t(x)=u\ast\varphi_t(x)\rightarrow u(x)$  in
$BL^{p,s}_{|x|^{\alpha}}$  as $t\rightarrow 0.$
\par {\bf Proof}  Let $u\in BL^{p,s}_{|x|^{\alpha}}$. Since $|\varphi(x)|<C,$ we
have
 $ u_t(x)\leq CMu(x)$. For any $\varepsilon>0$, by Theorem 2.10 above, we can find a $v(x)\in
 C_c^0$ such that
$$
\|u-v\|_{BL^{p,s}_{|x|^{\alpha}}}\leq  \varepsilon.\eqno(2.7)
 $$
For $t$ small enough, using Theorem 2.8 and (2.7), and the
boundedness of $M$ on $BL^{p,s}_{|x|^{\alpha}}$ (Theorem 4.4 below),
we have
\begin{eqnarray*}
\|u_t-u\|_{BL^{p,s}_{|x|^{\alpha}}} &\leq &
\|u_t-v_t\|_{BL^{p,s}_{|x|^{\alpha}}}+
\|v_t-v\|_{BL^{p,s}_{|x|^{\alpha}}}+\|u-v\|_{BL^{p,s}_{|x|^{\alpha}}}.
 \\&\leq &
\|M(u-v)\|_{BL^{p,s}_{|x|^{\alpha}}}+2\varepsilon
\\&\leq &
C\|u-v\|_{BL^{p,s}_{|x|^{\alpha}}}+2\varepsilon\leq C\varepsilon,
 \end{eqnarray*}
where $M$ is Hardy-Littlewood maximal operator.
\par Theorem 2.9 implies that

\par {\bf Theorem 2.10} Let   $0<p\leq 1,1< s<\infty$ and $
-n(1-p/s)<\alpha<n(p-1),$ or let  $1<p<s<\infty$ and $
-n(1-p/s)<\alpha \leq 0.$ Then
  $C^\infty $ is dense in $
BL^{p,s}_{|x|^{\alpha}} $.

\par In fact, we have
\par {\bf Theorem 2.11} Let   $0<p<\infty,1\leq s<\infty$ and $
-n<\alpha<\infty.$

\par 1) $C_c^\infty $ is dense in $
BL^{p,s}_{|x|^{\alpha}} $,
\par 2) $\mathcal{S} $ is dense in $
BL^{p,s}_{|x|^{\alpha}} $.

\par {\bf Proof}  Since $
C_c^\infty \subset \mathcal{S},$ 1) implies 2). We just need to
prove 1). For $f\in BL^{p,s}_{|x|^{\alpha}}$ and any
$\varepsilon>0$, by Theorem 2.7 above, we can find a $g\in
 C_c^0$ such that
\begin{eqnarray*}
\|f-g\|_{BL^{p,s}_{|x|^{\alpha}}}< \varepsilon.
 \end{eqnarray*}
  As $t$ is small enough, by Theorem 2.8, we have
\begin{eqnarray*}
\|g_t-g\|_{BL^{p,s}_{|x|^{\alpha}}}<  \varepsilon.
 \end{eqnarray*}
So,
\begin{eqnarray*}
\|f-g_t\|_{BL^{p,s}_{|x|^{\alpha}}}\leq
\|f-g\|_{BL^{p,s}_{|x|^{\alpha}}}+\|g_t-g\|_{BL^{p,s}_{|x|^{\alpha}}}<2\varepsilon.
 \end{eqnarray*}

\par {\bf Theorem 2.12} Let   $0<p<+\infty,1\leq s<\infty$ and $ -n<\alpha<\infty.$
Then ${\mathcal{S}}\subset BL^{p,s}_{|x|^{\alpha}}$.
\par {\bf Proof}  Let $f\in  \mathcal{S}$, we have $|f(x)|<\frac{C_N}{(1+|x|^2)^N}$
for all $N\in \mathbf{N}$ and $x\in {\bf R}^n$.
Write
$$f(x)=\chi_{B_1}(x)f(x)+\sum_{k=1}^{\infty}\chi_{C_k}(x)f(x)=\lambda_1\frac{1}{\lambda_1}\chi_{B_1}(x)f(x)+
\sum_{k=2}^{\infty}\lambda_k\frac{1}{\lambda_k}\chi_{C_k}(x)f(x),$$
  here $B_k=B(0,2^k)$, and $ C_k=B_k\setminus B_{k-1}, k=1,2,\cdots.$ We have
\begin{eqnarray*}
\|\frac{1}{\lambda_k}\chi_{C_k} f \|_{L^{ s} } &\leq&
\frac{1}{\lambda_k}\frac{1}{(1+2^{2(k-1)})^N}|C_k|^{\frac{1}{s} }
\\ &\leq &
\frac{1}{\lambda_k}\frac{2^{k(\frac{\alpha}{pn}-\frac{1}{p})}}
{(1+2^{2(k-1)})^N}|C_k|^{-\frac{\alpha}{pn}-\frac{1}{p}+\frac{1}{s}
}
\\&<& |C_k|^{-\frac{\alpha}{pn}-\frac{1}{p}+\frac{1}{s} }
 \end{eqnarray*}
if
$$
\frac{1}{\lambda_k}\frac{2^{k(\frac{\alpha}{pn}+\frac{1}{p})}}{(1+2^{2(k-1)})^N}
 <1
 {\rm~~i.e.~~}\lambda_k>\frac{2^{k(\frac{\alpha}{pn}+\frac{1}{p})}}{(1+2^{2(k-1)})^N}.
 \eqno(2.8)
 $$
Let $N=\left[\frac{\alpha}{pn}+\frac{1}{p}\right]+2$ and
$\lambda_k=\frac{2^{k(\frac{\alpha}{pn}+\frac{1}{p})}}{\left(1+2^{2(k-1)}\right)^{
[\frac{\alpha}{pn}+\frac{1}{p} ]+1}}$, then (2.8) holds, and
\begin{eqnarray*}
 \lambda_k<\frac{2^{k(\frac{\alpha}{pn}+\frac{1}{p})}}{\left(1+|2|^{2(k-1)}\right)^{
\frac{\alpha}{pn}+\frac{1}{p}  }}\leq \frac{C}{  2 ^{
 k      (\frac{\alpha}{pn}+\frac{1}{p} ) }},
 \end{eqnarray*}
noticing $-n<\alpha.$ Thus, $\frac{1}{\lambda_k}\chi_{C_k}(x)f(x)$
are $(p,s,\alpha)$-blocks for $k=2,3,\cdots,$ and
$\sum_{k=2}^{\infty}|\lambda_k|^p\leq \sum_{k=2}^{\infty}\frac{C}{ 2
^{
 k      (\frac{\alpha}{n}+1) }}<\infty$. The   discuss above also applies to
 $\chi_{B_1}(x)f(x)$. Therefore, We obtain that
 $f\in BL^{p,s}_{|x|^\alpha}$.

\par \section*
 {\bf  3. A molecular theorem  }
\par The molecular theory for Hardy spaces $H^p$ was  established by
Coifman \cite{Coif}, Coifman and Weiss \cite{CW}, and Taibleson and
Weiss  \cite{TW2}. The weighted case can be found in \cite{LL}. In
this section, we will prove a molecular theorem for
$BL^{p,s}_{|x|^{\alpha}}$. Next, from Section 4 to Section 7 we will
make use of this theorem to obtain some estimates beyond endpoint
for some operators.

\par First, let us introduce a definition.
 \par {\bf Definition 3.1} Let $0 < s <+ \infty, 0< p <+\infty,- \infty < \alpha < +\infty
 , -\infty <A_0,B_0 <\infty, A_0-B_0
 =\frac{1}{s}-\frac{1}{p}-\frac{\alpha}{np}$ and $
  -\infty<\varepsilon < A_0. $ Set $
   a=A_0 -\varepsilon ,$ and $b=B_0-\varepsilon $.
  A function $M(x)\in L^s $
 is said to be a  $ (p, s, \alpha, \varepsilon )$-molecular (centerred at  $x_0$), if
\par (i)~~~~$M(x)|x-x_0|^{nb}\in L^s ,$
\par (ii)~~~~$\|M\|^{a/b}_{L^{s}   }\|M(x)|x-x_0|^{nb}\|^{1-a/b}_{L^{s}   }
\equiv \Re (M)<\infty.$

\par {\bf Theorem 3.1} ~Let $  p,s,\alpha  ,  \varepsilon  , a ,b$
as in Definition 3.1.
   Every  $ (p, s, \alpha, \varepsilon )$-molecular $M(x)$  centered
at any point is in $BL^{p,s}_{|x|^{\alpha}}$ and $\|M \|^{
}_{BL^{p,s}_{|x|^{\alpha}}}\leq C \Re(M)$. Here the constant $C$ is
independent of $M$.

\par {\bf Proof}   Without loss of generality, we can assume
$\Re(M)=1$. In fact, assume $\|M \|^{ }_{\dot{B}L^{p,s}
_{|x|^{\alpha}}}\leq C  $ holds whenever $\Re(M)=1$. Then, for
general $M,$ let $M'=M/\Re(M)$. We have  $\Re(M')=1$ and hence
$\|\Re(M)M' \|^{ }_{BL^{p,s}_{|x|^{\alpha}}}\leq \Re(M)\| M' \|^{
}_{BL^{p,s}_{|x|^{\alpha}}}\leq C \Re(M)$. And we can also assume
$M(x)$  has the center at $x_0=0$ by a translation transformation.
\par
Let $M$ be a $ (p, s, \alpha, \varepsilon )$-molecular centered at
$ 0$
 satisfying $\Re(M)=1$. Define $B(0,r)$ by setting
$\|M\|^{ }_{L^{s}   }= |B(0,r)|^{1/s-1/p-\alpha/np}$. Let
$2^{k_0-1}<r\leq2^{k_0}$,
 and consider the set $ E_0=B_{k_0}(0), E_k=C_{k_0+k}(0),$ for $ k=1, 2, ...
 $, here $B_{k_0}(0)=B(0,2^{k_0}), C_{k_0+k}(0)=B_{k_0+k}(0)\setminus B_{k_0+k-1}(0).$
 Set $M_k=M{\chi _{E_k}}$.
 \par By $\Re (M)=1$, we have
$$
 \|M(x)|x |^{nb}\|^{1-a/b}_{L^{s}   }=\|M\|^{-a/b}_{L^{s}   } .
$$
  Then,
$$
 \|M(x)|x |^{nb}\|^{ }_{L^{s}   }=\|M\|^{-\frac{a}{b}\frac{b}{b-a}}_{L^{s}   }
 =|B (0,r)|^{-\left(\frac{1}{s} -\frac{1}{p} -\frac{\alpha}{np} \right)
 \frac{a}{b-a}}=|B (0,r)|^a\leq C_{a,n}|B_{k_0}(0)|^a. \eqno(3.1)
$$
Thus, using (3.1), for $k=1,2,...$, we have
\begin{eqnarray*}
\int _{E_k} M_k^s(x)dx &=&\int _{E_k} M_k^s(x) |x|^{nbs}|x|^{-nbs}dx
\\
& \leq &
 \left(2^{k_0+k-1}\right)^{-nbs}\int _{E_k} M_k^s(x) |x|^{nbs} dx
\\
& \leq &
 C_{a,n}\left(2^{k_0+k-1}\right)^{ns(a-b)}2^{-knsa}
\\
& \leq &
 C_{a,b,s,n} 2^{-knsa} |B_{k_0+k}(0)|^{(1/s-1/p-\alpha/np)s},
\end{eqnarray*}
 when $b>0;$ and
\begin{eqnarray*}
\int _{E_k} M_k^s(x)dx &=&\int _{E_k} M_k^s(x) |x|^{nbs}|x|^{-nbs}dx
\\
& \leq &
 \left(2^{k_0+k+1}\right)^{-nbs}\int _{E_k} M_k^s(x) |x|^{nbs} dx
\\
& \leq &
 C_{a,n}\left(2^{k_0+k+1}\right)^{ns(a-b)}2^{-knsa}
\\
& \leq &
 C_{a,b,s,n} 2^{-knsa} |B_{k_0+k}(0)|^{(1/s-1/p-\alpha/np)s},
\end{eqnarray*}
when $b\leq 0.$ And
$$\|M_0\|^{ }_{L^{s}   }\leq  \|M\|^{ }_{L^{s}   }= |B(0,r)|^{1/s-1/p-\alpha/np}
\leq C_{\alpha,p,s,n}|B_{k_0}(0)|^{1/s-1/p-\alpha/np}.$$ Hence, for
$k=0, 1,2,...$ , we have
$$\|M_k\|^{ }_{L^{s}   } \leq C2^{-kna}|B_{k_0+k}(0)|^{1/s-1/p-\alpha/np},$$
where $C$ is an absolute constant. Let
$$
a_k(x)=2^{kna}M_k(x), ~~k=0,1,2,...,
$$
 then
$$
M(x)=\sum \limits_{k=0}^{\infty}M_k(x)=\sum \limits_{k=0}^{\infty}
2^{-kna}2^{kna}
 M_k(x)=\sum \limits_{k=0}^{\infty}2^{-kna}a_k(x),$$
and each $a_k$ is a $ (\alpha ,p, s )$-block centered at  $ 0$ with
${\rm supp}a_k \subset B_{k_0+k}(0),$ and $\sum
\limits_{k=0}^{\infty}2^{-kna\bar{p}}=C<\infty$ since $a>0$. That is
$\|M \|^{ }_{BL^{p,s}_{|x|^{\alpha}}}\leq C .$
 Thus, we finish the proof of Theorem 3.1.

\par \section*
 {\bf  4. Hardy-Littlewood maximal operator}

In this section, we hope to get some estimates of Hardy-Littlewood
maximal operator $M $ on     $ BL^{p,s}_{|x|^{\alpha}}$. In fact, we
prove that some sublinear operators satisfying (1.14) and (1.15) are
bounded on $BL^{p,s}_{|x|^{\alpha}}$.

\par {\bf Theorem 4.1} Let
 $    1< s < \infty,0< p \leq s$ and $ - n(1-p/s)\leq \alpha < n(p-1)$.
 Suppose that a
  sublinear operator $ H $  satisfies
(1.14) and  is bounded
 on ${L}^{s} $. Then
$$\|H h\|_{BL^{p,s}_{|x|^{\alpha}}}\leq
C  $$ for every $(  p,s,\alpha)-$block $h$, where constant $C$ is
independent of $h$.

\par
\par {\bf Proof  }
It  suffices to check that $Hh$ is a  $ (p, s, \alpha, \varepsilon
)$-molecular
 for every $(  p,s,\alpha)-$block $h$ centered at any
point and $\Re (Hh)\leq C $, where $C$ is independent of $h$.

\par  Since  $1
 -1/p-\alpha/np>0$ and $1-1/s>0$, we can choose $\varepsilon$ such that
$$0<\varepsilon <\min\{1 -1/p-\alpha/np,1-1/s\}.$$
 Set $ a=1-1/p-\alpha/np-\varepsilon $ and $ b=1-1/s-\varepsilon. $ We
see that $a>0,b>0$ and $b-a=-(1/s-1/p-\alpha/np)\geq 0$  by $-
n(1-p/s)\leq \alpha$. Given a $( p,s, \alpha)-$block $h$ with
 supp$h\subset B(x_0,r)\subset B_{k_0}(x_0)$ and
 $2^{k_0-1}<r\leq2^{k_0}$,
we have
\begin{eqnarray*}
\|Hh(x)|x-x_0 |^{nb}\|^{s}_{L^{s}   } &=&\int _{{\bf R}^n }
|Hh(x)|^s |x-x_0|^{snb} dx
\\
& \leq &
 \left(\int _{|x-x_0|<2^{k_0+ 1}}+\int _{|x-x_0|\geq 2^{k_0+ 1}} \right)
|Hh(x)|^s |x-x_0|^{snb} dx
 \\& =&  J_1+J_2.
\end{eqnarray*}

 \par For $J_1,$ noticing that $b>0$, and  by $L^s$ boundedness of $H$, we have
\begin{eqnarray*}
J_1&=& \int _{|x-x_0|< 2^{k_0+ 1}} |Hh(x)|^s |x-x_0|^{snb} dx
\\
& \leq & C 2^{sk_0nb}\|Hh \|^{s}_{L^{s}   }
 \\
& \leq & C 2^{sk_0nb}\|h \|^{s}_{L^{s}   }.
\end{eqnarray*}

For  $J_2,$    by (1.14) and H\"{o}lder inequality, and noticing
that   $b+1/s-1=-\varepsilon <0$,
 we have
\begin{eqnarray*}
J_2&=&  C \int _{|x-x_0|\geq 2^{k_0+ 1}} |Hh(x)|^s |x-x_0|^{snb} dx
\\
& \leq &    \int _{|x-x_0|\geq 2^{k_0+ 1}} \|h \|^{s}_{L^{1} }
|x-x_0|^{snb-sn} dx
\\
& = & C  2^{sk_0n(b-1+1/s)}\| h \|^{s}_{L^{1}   }
\\
& \leq & C 2^{sk_0nb}\|h \|^{s}_{L^{s}   }.
\end{eqnarray*}
Thus,
$$
\|Hh(x)|x-x_0 |^{nb}\|^{ }_{L^{s}   } \leq  C 2^{k_0nb}\|h
\|^{}_{L^{s}   }\leq  C 2^{k_0nb}2^{k_0n(1/s-1/p-\alpha/np)}.
$$
And by the boundedness of $H$ on $L^s$,
$$
\|Hh \|^{ }_{L^{s}   } \leq  C \| h  \|^{ }_{L^{s}   }
 \leq  C  2^{k_0n(1/s-1/p-\alpha/np)}.
$$
Then, noticing $a/b>0$ and $1-a/b\leq 0,$ we have
\begin{eqnarray*}
\Re (Hh)&\leq& C2^{k_0n(1/s-1/p-\alpha/np)a/b}
2^{k_0n(b+1/s-1/p-\alpha/np)(1-a/b)}
\\
&=&C2^{k_0n(b-a+1/s-1/p-\alpha/np)} =C2^0=C.
\end{eqnarray*}

Thus, we finish the proof of Theorem 4.1.

\par From Theorem 4.1 and the subadditivity of norm, it is easy to
get the following theorem.

\par {\bf Theorem 4.2}  Let
 $    1< s < \infty,0< p \leq s$ and $ - n(1-p/s)\leq \alpha < n(p-1)$.
  Suppose that a
  sublinear operator $ H $ is bounded
 on ${L}^{s} $, and
  satisfies
(1.14)  and (1.15).
 Then,
  $$H ~{\rm is~ of~ type}~  (BL^{p,s}_{|x|^{\alpha}},BL^{p,s}_{|x|^{\alpha}}).$$

\par From Theorems 2.1  and   4.2, we have

\par {\bf Theorem 4.3} Let
 $    1<s < \infty,0< p <s ,-n(1-p/s) < \alpha < n(p-1)$ and $\alpha
 \leq 0$.
   Suppose that a
  sublinear operator $ H $ is bounded
 on ${L}^{s} $, and
  satisfies
(1.14)  and  (1.15).
 Then
  $$H ~{\rm is~ of~ type~} (BL^{p,s}_{|x|^{\alpha}},{L}^{p}
  _{|x|^{\alpha}}).$$

\par Hardy-Littlewood maximal function satisfies the conditions in  Theorem 4.1 $\sim$ Theorem
4.3.  Define   Hardy-Littlewood maximal function as
$$Mf(x)=\sup_{x\in B}\frac{1}{|B|}\int_B|f(y)|dy$$
for $f\in L^1_{loc}$. Then
 $M$ satisfies the size conditions (1.14).
 In fact, let ${\rm supp} f \subset  B_k(x_0)$ and $ x\in
(B_{k+1}(x_0))^c $. For any $y\in B_k(x_0)=B(x_0,2^k)$ and $x\in
(B_{k+1}(x_0))^c $ for $k\in {\bf Z}$, we have $|x-y|\geq
|x-x_0|/2$. Then, the diameter of the ball $B$, which includes $x$
and intersects with $ B_k(x_0)$, is not less than  $|x-x_0|/2$, and
it follows that $|B|\geq C|x-x_0|^n.$ So
$$Mf(x)=\sup_{x\in B}\frac{1}{|B|}\int_B|f(y)|dy
=\sup_{x\in B}\frac{1}{|B|}\int_{B\cap B_k(x_0)}|f(y)|dy\leq
\frac{C\|f\|_{L^1}}{|x-x_0|^{n}}
$$
for $x\in (B_{k+1}(x_0))^c.$

\par   By Minkowski inequality, it
is easy to see that $M$ satisfies (1.15). It is known that  $M$ is a
bounded operator  on $L^s $ with $1<s<\infty$, so we have

\par {\bf Theorem 4.4} Let $ 1<s < \infty,
  0< p\leq s $ and $- n(1-p/s) \leq \alpha < n(p-1)$, then
  $$M~{\rm ~ is~ of~ type~
  }(BL^{p,s}_{|x|^{\alpha}},BL^{p,s}_{|x|^{\alpha}})
 .  $$
And let
 $    1<s < \infty,0< p <s ,-n(1-p/s) < \alpha < n(p-1)$ and $\alpha
 \leq 0$, then
  $$M {\rm ~is~ of~ type ~}(BL^{p,s}_{|x|^{\alpha}},{L}^{p}
  _{|x|^{\alpha}}) . \eqno(4.1)$$

\par (4.1) is sharp in the sense that $M$
   fails to be of type $(BL^{p,s}_{|x|^{\alpha}}, {L}^{p}
  _{|x|^{\alpha}}) $ for $s=1$ or $n(p-1)\leq \alpha
<\infty.$ In fact, we have

\par {\bf Theorem 4.5} A) Let $s=1, 0<p\leq 1$ and $ -\infty <\alpha \leq
n(p-1)$. Then there exists $f\in BL^{p,s}_{|x|^{\alpha}}$ such that
 $Mf$
   fails to be of type $(BL^{p,s}_{|x|^{\alpha}}, {L}^{p}
  _{|x|^{\alpha}}) $.
  \par B)
Let   $0<p\leq s\leq \infty. $ If $\alpha=n(p-1),$ or $0<p\leq
1,n(p-1)\leq \alpha \leq 0, $ then $Mf$
   fails to be of  type $(BL^{p,s}_{|x|^{\alpha}}, {L}^{p}
  _{|x|^{\alpha}}) $ for all $f\in BL^{p,s}_{|x|^{\alpha}}.$
\par C) If
$ 0<s\leq \infty$ and $ 0<p<\infty, n(p-1)\leq \alpha< \infty, $
  then there exists $f\in BL^{p,s}_{|x|^{\alpha}}$ such that
 $Mf$
   fails to be of type $(BL^{p,s}_{|x|^{\alpha}}, {L}^{p}
  _{|x|^{\alpha}}) $.

\par {\bf Proof} When  $s=1, 0<p\leq 1$ and $-\infty <\alpha \leq n(p-1),$
 by Theorem 2.4 B), we
have $L^1\subset BL^{p,s}_{|x|^{\alpha}}$.

\par  When $0<p\leq s\leq \infty,  \alpha=n(p-1),$ or $0<p\leq
1,n(p-1)\leq \alpha \leq 0, $ by Theorem 2.4 A), we see that $
BL^{p,s}_{|x|^{\alpha}}\subset L^\frac{np}{n+\alpha},$ and
$\frac{np}{n+\alpha}\leq 1$.
\par
\par When $ 0<s\leq
\infty, 0<p<\infty$ and $n(p-1)\leq \alpha< \infty, $ by Theorem
2.5, we have $H^\frac{np}{n+\alpha}\subset BL^{p,s}_{|x|^{\alpha}}$,
and $\frac{np}{n+\alpha}\leq 1$.
\par The remain of the proof follows from  the following

\par {\bf Claim}: Let $0<p<\infty, n(p-1)\leq \alpha <\infty,0<q\leq 1$ and
$ f\in L^q$. Then $Mf$ is never in ${L}^{p}
  _{|x|^{\alpha}}$ if $ f  \neq 0.$

\par See \cite{Longs2} for a proof of Claim. Thus, we finish the
proof of Theorem 4.5.

\par  The conditions in  Theorem 4.1 $\sim$ Theorem
4.3 are also satisfied by   Littlewood-Paley functions.  Suppose
that $\psi$ is integrable on ${\bf R}^n$ and
\par (i) $\int_{{\bf
R}^n}\psi (x)dx=0,$
\par (ii) $|\psi (x)|\leq c(1+|x|)^{-(n+\beta)},$ for some $\beta>0,$
\par (iii) $\int_{{\bf
R}^n}|\psi (x+y)-\psi (x)|dx\leq c|y|^\gamma,$ all $y\in {\bf R}^n$,
for some $\gamma>0.$
\par Let $\psi_t(x)=t^{-n}\psi (x/t)$ with
$t>0.$
 The Littlewood-Paley $g$-function of $f$ is defined by
$$g_\psi(f)(x)=\left(\int_0^\infty |f\ast\psi_t(x)|^2\frac{dt}{t}\right)^{1/2}.$$
\par The Lusin area function of $f$ is defined by
$$S_{\psi,a}(f)(x)=\left(\frac{1}{a^n|B_0|}\int_{\Gamma_a(x)} |f\ast\psi_t(y)|^2t^{-n}dy\frac{dt}{t}\right)^{1/2},$$
where $B_0$ is the unit ball of ${\bf R}^n$ and
$\Gamma_a(x)=\{(y,t)\in {\bf R}^{n+1}_+:|x-y|<at\}$.
\par The Littlewood-Paley $g_{\lambda}^*$-function of $f$ is defined by
$$g_{\psi,\lambda}^*(f)(x)=\left(\int_0^\infty\int_{{\bf
R}^n}
\frac{|f\ast\psi_t(y)|^2}{(1+\frac{|x-y|}{t})^{2\lambda}}t^{-n}dy\frac{dt}{t}\right)^{1/2}.$$

 \par As in \cite{LY} we can check that  $g_\psi(f), S_{\psi,a}(f)$ and $ g_{\psi,\lambda}^*(f)$ satisfy the size
conditions (1.14).  By Minkowski inequality, it is easy to see that
$g_\psi(f), S_{\psi,a}(f)$ and $ g_{\psi,\lambda}^*(f)$ satisfy
(1.15). By the boundedness of these operators  on $L^s $ with
$1<s<\infty$, we have

\par {\bf Theorem 4.6}  Let $\psi$ satisfies (i),(ii) and (iii) above. For
$g_\psi(f), S_{\psi,a}(f)$ and $ g_{\psi,\lambda}^*(f)$, the same
conclusions hold as those stated in Theorem 4.4.

\par The conditions in Theorem 4.1 $\sim$ Theorem 4.3 are also satisfied by many other  operators arising
in harmonic analysis, and we will discuss this in section 7.

\par \section*
 {\bf  5. Bochner-Riesz means }

\par $B_R^\lambda$ can be written as a convolution operator
$$
B_R^\lambda f(x)=(f\ast K_R^\lambda )(x)
$$
for $f\in \mathcal{S}$, where $K^\lambda_R(x) =
[(1-|\xi/R|^2)^\lambda_+]~\breve{}~(x),$ and $\check{g}$ denotes the
inverse Fourier transform of $g$.

 We will prove that

\par {\bf Theorem 5.1} Let $
\frac{n-1}{2(n+1)}< \lambda\leq \frac{n-1}{2}, p'_\lambda
<s<p_\lambda, 0<p\leq s $ and $ n(\frac{p}{s}-1)\leq
\alpha<n(\frac{p}{p'_\lambda}-1)$. Then (1.3) holds.

\par {\bf Theorem 5.2}  Let $0<R<\infty,$ and let $
\frac{n-1}{2(n+1)}< \lambda\leq \frac{n-1}{2}, p'_\lambda
<s<p_\lambda, 0<p< s , n(\frac{p}{s}-1)<
\alpha<n(\frac{p}{p'_\lambda}-1)$ and $\alpha\leq 0$. Then
  (1.4) holds.

\par Since Theorem 2.5, we can extend this to $p_\lambda \leq s <
\infty$.
\par {\bf Corollary 5.1}  Let $0<R<\infty,$ and let $
\frac{n-1}{2(n+1)}< \lambda\leq \frac{n-1}{2},  s=p_\lambda, 0<p< s
, n(\frac{p}{p_\lambda}-1)< \alpha<n(\frac{p}{p'_\lambda}-1)$ and
$\alpha\leq 0$. Then
  (1.4) holds.

\par When $
  \lambda= \frac{n-1}{2}$, these results can extend to the maximal
  Bochner-Riesz operators. Noticing that $ p'_{\frac{n-1}{2}}=1,   p_{\frac{n-1}{2}}=\infty,
  $ we have

\par {\bf Theorem 5.3} Let $
  \lambda= \frac{n-1}{2},1
<s<\infty, 0<p\leq s $ and $ n(\frac{p}{s}-1)\leq \alpha<n(p-1)$.
Then (1.6) holds.

\par {\bf Theorem 5.4} Let $
  \lambda= \frac{n-1}{2},1
<s<\infty, 0<p< s ,  n(\frac{p}{s}-1)<\alpha<n(p-1)$ and $\alpha\leq
0$.  Then
 (1.7) holds.

\par When $
  \lambda>\frac{n-1}{2}$, we have also

\par {\bf Theorem 5.5} Let $\lambda  > \frac{n-1}{2}, 1<s<\infty, 0<p\leq s $ and $
n(\frac{p}{s}-1)\leq \alpha<n(p-1)$.  Then (1.6) holds.

\par {\bf Theorem 5.6} Let $
  \lambda> \frac{n-1}{2},1
<s<\infty, 0<p< s ,  n(\frac{p}{s}-1)<\alpha<n(p-1)$ and $\alpha\leq
0$.  Then
 (1.7) holds.

\par We have the following convergence results.

\par {\bf Theorem 5.7} (i) Let $\lambda   , s , p $ and $
  \alpha $ under the conditions of Theorem 5.2, or
  Theorem 5.4, or Theorem 5.6.
     Then  (1.5) holds.

\par It is  known that the boundedness of $B^\lambda_*$ on $L^p$ and
the maximal principle (see Stein \cite{S1}) imply that  $B_R^\lambda
f$ converges to $f$ almost everywhere
 as $R$ tends to $\infty$, for $f\in L^p$. Similarly, we can also get a $\mu_\alpha$-a.e.
 convergence result for $B_R^\lambda
f$ by Theorem 5.3 (or Theorem 5.5) and a similar  maximal
 principle, for $f\in  BL^{p,s}
_{|x|^{\alpha}}$   with $\lambda , s , p $ and $
  \alpha $ under the conditions of Theorem 5.3 (or  Theorem 5.5). But
  this does not imply new information since $BL^{p,s}_{|x|^{\alpha}}\subset
  L^{\frac{np}{n+\alpha}}$   and  $\frac{np}{n+\alpha}>1$ (see Theorem
  2.4).

\par Next, We prove the Theorems above. First, let us give an outline of the proof of   Theorem 5.1.
In order to prove that  $ B_1^\lambda $ is  of type $
(BL^{p,s}_{|x|^{\alpha}},BL^{p,s}_{|x|^{\alpha}})$, by the equality
$B_1^\lambda f= Hf +Tf$ (i.e.(5.6) below, where $Hf$ is controlled
by Hardy-Littlewood maximal function $Mf$),
 we need only to prove that this  holds for $T$   since
Theorem 4.4.  This can be  divided into three steps.  Step 1, we
prove that
$$
\|B_1^\lambda h\|_{BL^{p,s}_{|x|^{\alpha}}} \leq C \eqno (5.1)
$$
for all $(p,s,\alpha)$-blocks $h$ of restrict II-type, where $C$ is
independent of $h$, and it follows that
$$
\|T h\|_{BL^{p,s}_{|x|^{\alpha}}} \leq C \eqno (5.2)
$$
 for all $(p,s,\alpha)$-blocks $h$ of restrict II-type, since Theorem 4.4.
 Step 2, we prove that (5.2) holds  for all $(p,s,\alpha)$-blocks $h$ of restrict
 I-type. Thus, (5.2) holds  for all $(p,s,\alpha)$-blocks
 $h$. Step 3,  we
prove that
$$
Tf=\sum_{j}\mu_jT a_j \eqno(5.3)
$$
for $f=\sum_{j}\mu_j a_j \in BL^{p,s}_{|x|^{\alpha}}$ with $s,p $
and $\alpha$ under the conditions of Theorem 5.1. From these, it is
easy to get that  $T$  is  of type $ (
BL^{p,s}_{|x|^{\alpha}},BL^{p,s}_{|x|^{\alpha}})$, then $
B_1^\lambda $ follows since (5.6) below and Theorem 4.4.

\par {\bf Proof of Theorem 5.1}
\par  {\bf Step 1}.
\par  we consider the case of $(p,s,\alpha)$-block  of restrict
II-type.

\par {\bf Proposition 5.1} Let $\lambda\leq \frac{n-1}{2}, p'_\lambda <s<\infty, 0<p\leq s $ and $
n(\frac{p}{s}-1)\leq \alpha<n(\frac{p}{p'_\lambda}-1)$. If
$B_1^\lambda$ is bounded on $L^s$, then (5.1) holds for all
$(p,s,\alpha)$-blocks $h$ of restrict II-type.

\par In fact, we can prove the following
\par {\bf Proposition  5.2} Let $  1<s<\infty, 0<p\leq s,
n(\frac{p}{s}-1)\leq \alpha<n(p-1)$ and $ \frac{n+\alpha}{p}<\delta
\leq n$. Suppose that a sublinear operator $K$ satisfies the size
condition as following
$$
|Kf(x)| \leq C \|f\|_{L^{1}}/|x-x_0|^\delta \eqno(5.4)
$$
when supp $f\subseteq B(x_0,2^k)$ and $|x-x_0| \geq 2^{k+1}$ with
$k\in {\bf Z} $. If $K$ is bounded on $L^s$, then
$$
\|K h\|_{BL^{p,s}_{|x|^{\alpha}}} \leq C
$$
for all $(p,s,\alpha)$-blocks $h$ of restrict II-type, where $C$ is
independent of $h$.

\par We can check that  the following size condition implies the
condition (5.4):
$$|Kf(x)|\leq C \int_{{\bf R}^n}\frac{|f(y)|}{|x-y|^\delta}dy, ~~~~ x\notin
{\rm supp}f, \eqno{(5.5)}$$ for any integral function with compact
support.  In fact, assume that $K$ satisfies (5.5). Let $y\in
B_k(x_0)=B(x_0,2^k) $ and $x\in (B_{k+1}(x_0))^c $ (i.e.
$|x-x_0|\geq 2^{k+1}$). For $k\in {\bf Z}$, we have  $|x-y|\geq
2^k\geq |y-x_0| ,$ and it follows that $|x-x_0|\leq |x-y| +
|y-x_0|\leq 2|x-y|$. Then $K$ satisfies (5.4).

 \par It is known that
$$|B_1^\lambda f(x)|\leq C \int_{{\bf R}^n}\frac{|f(y)|}{|x-y|^{\frac{n+1}{2}+\lambda}}dy   $$
for $f\in L^s$ with $1<s<\infty$.

\par Let $\delta=\frac{n+1}{2}+\lambda$, it is easy to check that $\lambda,s,p$ and $
\alpha$ satisfy the conditions of Proposition 5.1 if they satisfy
Proposition 5.2. So, Proposition 5.1 is a special case of
Proposition 5.2 when $\delta=\frac{n+1}{2}+\lambda$.

\par {\bf Proof of Proposition 5.2 }  By Theorem 3.1,
it   suffices to check that $K h$ is a $ (p, s, \alpha, \varepsilon
)$-molecular
 for every $(  p,s,\alpha)-$block $h$ of restrict II-type
 centered at any point $x_0$ and $\Re (Kh)\leq C $ with $C$ independent of $h$.

\par  Since  $n(\frac{p}{s}-1)\leq \alpha $ and
$\frac{n+\alpha}{p}<\delta$,
 we have $0<\frac{\delta}{n}-\frac{1}{p}-\frac{\alpha}{np} \leq \frac{\delta}{n}-\frac{1}{s}
$. Then, we can choose $A_0 $ and $ \varepsilon$ such that
$$0<\varepsilon < A_0 =  \frac{\delta}{n}-\frac{1}{p}-\frac{\alpha}{np} \leq  \frac{\delta}{n}-\frac{1}{s}
.$$
 Set $ a=A_0-\varepsilon $ and $ b= \frac{\delta}{n}-\frac{1}{s}-\varepsilon. $ We
see that $a>0,b>0 $ and $b-a
=-\frac{1}{s}+\frac{1}{p}+\frac{\alpha}{np}\geq 0$ since
$n(\frac{p}{s}-1)\leq \alpha $. Given a $( p,s, \alpha)-$block $h$
of restrict II-type with supp $a\subset B(x_0,r)\subset
B_{k_0}(x_0)$, $2^{k_0-1}<r\leq2^{k_0}$, and $ k_0\leq 0$, we see
that
\begin{eqnarray*}
\|K h(x)|x-x_0 |^{nb}\|^{s}_{L^{s}   }  & =&
 \left(\int _{|x-x_0|<2^{k_0+ 1}}+\int _{|x-x_0|\geq 2^{k_0+ 1}} \right)
|K h(x)|^s |x-x_0|^{snb} dx
\end{eqnarray*}
$$
 =  J_1+J_2.~~~~~~~~~~~~~~~~~~~~~~~~~~~~~$$
 For $J_1,$ by the $L^{s}   $-boundedness of $K$,   we have
\begin{eqnarray*}
J_1   \leq  C 2^{k_0nsb}\|h \|^{s}_{L^{s}   }  \leq  C
2^{k_0nsb}2^{k_0ns(1/s-1/p-\alpha/np)}
  .
\end{eqnarray*}
 For $J_2,$ by (5.4),  noticing that $(nb-\delta)s+n=-n\varepsilon <0,   \delta\leq n$ and $k_0\leq 0$,
 and using H\"{o}lder inequality,  we
have
\begin{eqnarray*}
J_2&=& \int _{|x-x_0|\geq 2^{k_0+ 1}} |K h(x)|^s |x-x_0|^{nbs} dx
\\
& \leq & C \int _{|x-x_0|\geq 2^{k_0+ 1}} \|h \|^{s}_{L^{1}   }
|x-x_0|^{nbs-\delta s} dx
\\
& = & C2^{k_0((nb-n)s+n  )} 2^{k_0 (n -\delta )s}\| h \|^{s}_{L^{1}}
\\
& \leq & C2^{k_0ns(b-1+1/s)}\| h \|^{s}_{L^{s}   } |B_{k_0}|^{s/s'}
\\
& \leq & C 2^{k_0nsb}2^{k_0ns(1/s-1/p-\alpha/np)}.
\end{eqnarray*}
Thus,
$$
\|Kh(x)|x-x_0 |^{nb}\|^{ }_{L^{s}   } \leq  C
2^{k_0nb}2^{k_0n(1/s-1/p-\alpha/np)}.
$$
And by the $L^s$-boundedness of $K$, we have
$$
\|Kh \|^{ }_{L^{s}   } \leq  C \| h \|^{ }_{L^{s}   }
 \leq  C  2^{k_0n(1/s-1/p-\alpha/np)}.
$$
Then, noticing $a/b>0$ and $1-a/b\geq 0,$   $b-a +\frac{1}{s}
-\frac{1}{p}-\frac{\alpha}{np} = 0,$ we have
\begin{eqnarray*}
\Re (K h)&\leq& C2^{k_0n(1/s-1/p-\alpha/np)a/b}
2^{k_0n(b+1/s-1/p-\alpha/np)(1-a/b)}
\\
&=&C2^{k_0n(b-a+1/s-1/p-\alpha/np)} =C,
\end{eqnarray*}
since $k_0\leq 0$. Thus, we finish the proof of Proposition 5.2.

\par {\bf Step 2}.

\par  We consider the case of $(p,s,\alpha)$-block  of restrict
I-type.

  \par Because $B_1^\lambda = f\ast K_1^\lambda$, we can decompose
  $B_1^\lambda$ into two parts as $B_1^\lambda = f\ast (\chi_{B_{2\pi}} K_1^\lambda )+f \ast (\chi_{B_{2\pi}^c}
  K_1^\lambda)$, where $B_{2\pi}=\{x:|x|\leq 2\pi\}$. It is easy to see that the first part can be
  controlled by Hardy-Littlewood maximal function $Mf$ since
  $K_1^\lambda$ is bounded. For the second part, we can use the
  asymptotic expansion of $K_1^\lambda$
$$
K_1^\lambda(x)\sim |x|^{-(n+1)/2-\lambda}\left[ e^{2\pi i |x|}\sum
_{j=0}^\infty \alpha_j |x|^{-j}+e^{-2\pi i |x|}\sum _{j=0}^\infty
\beta_j |x|^{-j}\right]
$$
as $|x|\rightarrow \infty$, for suitable constants $\alpha_j$ and
$\beta_j$, and express it  as a finite sum as follows, as Stein did
in \cite{S3}. First there are finitely many terms, given by some
constant multiples of
$$T^{j \pm }f(x)=\int_{|y|\geq {2\pi}} e^{\pm 2\pi i |y|} f(x-y)|y|^{-(n+1)/2-\lambda-j}dy$$
with $j\geq 0$. Next, there is an error term, corresponding to
convolution with a kernel belonging to $L^1$. In fact, the error
term can be controlled by  $Mf$, and this can be seen more clearly
in the following Lemma 5.1 (see \cite{An}), and  in  \cite{DC} for
$\lambda>\frac{n-3}{2}$.

\par {\bf
Lemma 5.1}(\cite{An})  For every $k\geq -\frac{1}{2}$, there are
real constants
 $\{c_{k,j}\}^\infty _{ j =0}$
 and $c_k$ such that, for every $N \in {\bf Z}^+
 $, the following holds:
 $$ J_k(x) = R_{k,N}\left(\frac{1
}{x^{2N+5/2}}\right) +\sum ^N _{ j=0 }\left(  \frac{c_{k,2j}}{
x^{2j+1/2}} \cos(x - c_k) - \frac{c_{k,2j+1}}{ x^{2j+3/2}} \sin(x -
c_k)\right)$$
 where $R_{k,N}$
 is controlled by
$$
\left|R_{k,N}\left(\frac{1}{x^{2N+5/2}}\right)\right| \leq
\frac{C_{k,N} }{x^{2N+5/2}}
$$
 for every $x \geq 2 \pi $ and some constants $C_{k,N}$.

 \par Thus, we need only to take $N_0$ in Lemma 5.1 such that $ (n+1)/2+\lambda +2N_0+5/2>
 n$. Then
 we can write
$$B_1^\lambda f= Hf +Tf\eqno (5.6)$$
where $Hf$ is
  controlled by Hardy-Littlewood maximal function $Mf$, and
   $$Tf=\sum_{j=0}^{N_0}[
 (c^+_{2j+\frac{1}{2}} T^{(2j+\frac{1}{2})+}f+c^-_{2j+\frac{1}{2}} T^{(2j+\frac{1}{2})-}f)
 +(c^+_{2j+\frac{3}{2}} T^{(2j+\frac{3}{2})+}f+c^-_{2j+\frac{3}{2}} T^{(2j+\frac{3}{2})-}f)]$$
 for  some constants $c^{\pm}_j$, where
$$T^{(2j+\frac{1}{2}) \pm }f(x)=\int_{|y|\geq 2\pi}
 e^{\pm 2\pi i |y|} f(x-y)|y|^{-\frac{n+1}{2}-\lambda-(2j+\frac{1}{2})}dy$$
$j=0,1,2, \cdots, N_0. T^{(2j+\frac{3}{2}) \pm }f $ are similar.

\par {\bf Proposition 5.3}  Let $
\frac{n-1}{2(n+1)}< \lambda\leq \frac{n-1}{2}, p'_\lambda
<s<p_\lambda, 0<p\leq s$ and $n(\frac{p}{s}-1)\leq
\alpha<n(\frac{p}{p'_\lambda}-1)$. Then,   (5.2)   holds
 for all $(p,s,\alpha)$-blocks $h$ of restrict I-type.

\par {\bf Proof  }  We prove  (5.2)   by  Theorem 3.1.

\par It is easy to see that $T^{(2j+\frac{1}{2}) \pm }f(x)$ is
controlled pointwise by Hardy-Littlewood maximal function $Mf$ when
$\frac{n+1}{2}+\lambda +2j+\frac{1}{2}>n$, and $T^{(2j+\frac{3}{2})
\pm }f(x)$ is controlled by   $Mf$ when $\frac{n+1}{2}+\lambda
+2j+\frac{3}{2}>n$. It is clear that we need only to consider
$$T^{0 + }f(x)=\int_{|y|\geq {2\pi}} e^{ 2\pi i |y|} f(x-y)|y|^{-(n+1)/2-\lambda}dy.$$
 In fact, once we have proved that
$$
\|T^{0+} h\|_{BL^{p,s}_{|x|^{\alpha}}} \leq C \eqno(5.7)
$$
for all $(p,s,\alpha)$-blocks $h$ of restrict I-type, under
conditions of Proposition 5.3, then, from this and Theorem 4.4, it
is easy to see that $T^{(2j+\frac{1}{2}) \pm }$ and
$T^{(2j+\frac{3}{2}) \pm }$, $j=0,1,2, \cdots, N_0,$   also satisfy
(5.7). Then, $T$ follows.

\par Let us prove (5.7). By Theorem 3.1, it  suffices to check that $T^{0+} h$ is a  $ (p,
s, \alpha, \varepsilon )$-molecular
 for every $(  p,s,\alpha)-$block $h$ of restrict I-type
 centered at origin $x_0=0$ and $\Re (T^{0+} h)\leq C $ with $C$ independent of $h$.

\par  Since $n(\frac{p}{s}-1)\leq \alpha<n(\frac{p}{p'_\lambda}-1),$  we
have
$$   \frac{n+1}{2n}+\frac{\lambda}{n}-\frac{1}{s} \geq
\frac{n+1}{2n}+\frac{\lambda}{n}-\frac{1}{p} -\frac{\alpha}{np}>0.$$
Then, we can choose $A_0 $ and $ \varepsilon$ such that

$$ 0<\varepsilon <\frac{n+1}{2n}+\frac{\lambda}{n}-\frac{1}{p}
-\frac{\alpha}{np}
 \leq \frac{n+1}{2n}+\frac{\lambda}{n}-\frac{1}{s}
.$$

 Set $ a=\frac{n+1}{2n}+\frac{\lambda}{n}-\frac{1}{p}
-\frac{\alpha}{np}-\varepsilon $ and $
b=\frac{n+1}{2n}+\frac{\lambda}{n}-\frac{1}{s}-\varepsilon. $ We see
that $a>0,b>0 $ and $b-a
=-\frac{1}{s}+\frac{1}{p}+\frac{\alpha}{np}\geq 0$.

\par As in \cite{S3}, we write
$$
|y|^{-(n+1)/2-\lambda}=\sum_{k=0}^\infty
2^{-[(n+1)/2+\lambda]k}\psi(y/2^k), ~~~~~~|y|\geq 1
$$
where $\psi$ is a smooth function, supported in $1/2\leq |y| \leq
2$. Then
$$T^{0+}=\sum_{k=0}^\infty T_k \eqno(5.8)$$
 with
$T_kf(x)=2^{-[(n+1)/2+\lambda]k} \int e^{ 2\pi i |y|}
f(x-y)\psi(y/2^k)dy.$

\par Given a $( p,s,
\alpha)-$block $h$ with supp $h\subset B(0,r)\subset B_{k_0}(0)$,
$2^{k_0-1}<r\leq2^{k_0}$ and $ k_0>0$, by (5.8) and Minkowski
inequality, we have
$$\|T^{0+} h(x)|x |^{nb}\|_{L^{s}   }\leq \sum_{k=0}^\infty \|T_k h(x)|x |^{nb}\| _{L^{s}   }.$$
 Let us estimate $\|T_k h(x)|x |^{nb}\| _{L^{s}   }$. We see
 that
\begin{eqnarray*}
\|T_k h(x)|x |^{nb}\|^{s}_{L^{s}   } & \leq &
 \left(\int _{|x|<2^{k_0+ 1}}+\int _{|x|\geq 2^{k_0+ 1}} \right)
|T_kh(x)|^s |x|^{snb} dx
\\
& = & J_1^k+J_2^k.
 \end{eqnarray*}
For $J_1^k,$ by
(34) in (\cite{S3}, p394), we have
\begin{eqnarray*}
J_1^k&=& \int _{|x|<2^{k_0+ 1}} |T_kh(x)|^s |x|^{nbs} dx
\\
& \leq & C 2^{k_0nsb}\|T_kh \|^{s}_{L^{s}   }
\\
& \leq & C 2^{k_0nsb} 2^{-[(n+1)/2+\lambda]ks+kn} \| h \|^{s}_{L^{s}
} .
\end{eqnarray*}
Noticing that $-\frac{n+1}{2}- \lambda +\frac{n}{s}  <0$ (since
$p'_\lambda<s)$,  we have

$$\sum_{k=0 }^\infty (J_1^k)^{1/s}
 \leq C 2^{k_0nb } \| h \| _{L^{s} }  $$
 for $k_0=    1,   2,\cdots$.

For $J_2^k,$ let $x\in \{x:|x|\geq 2^{k_0+1}\}, y\in {\rm supp} \psi
(\cdot/2^k)=\{y:2^{k-1}<|y|\leq 2^{k+1}\}$ and $ x-y\in
B_{k_0}=\{x-y:|x-y|\leq 2^{k_0}\}$, then $2^{k_0+1}\leq |x|\leq
|y|+|x-y|\leq 2^{k_0}+ 2^{k+1}$ and $k_0\leq k+1$. Thus, noticing
that  $T_kh(x)=2^{-[(n+1)/2+\lambda]k} \int e^{ 2\pi i |y|}
h(x-y)\psi(y/2^k)dy,$
  we have $J_2^k=0$ when $k+1<k_0$, while when  $k+1\geq k_0$,
\begin{eqnarray*}
J_2^k&=& \int _{|x|\geq 2^{k_0+ 1}} |T_kh(x)|^s |x|^{nbs} dx
\\
& \leq & C  (2\times 2^{k+1})^{nsb} \int _{|x|\geq 2^{k_0+ 1}}
|T_kh(x)|^s
 dx
\\
& \leq & C 2^{knsb} 2^{-[(n+1)/2+\lambda]ks+kn} \| h \|^{s}_{L^{s} }
\\
& = & C 2^{ks(nb -(n+1)/2-\lambda+n/s)} \| h \|^{s}_{L^{s} } ,
\end{eqnarray*}
since (34) in (\cite{S3}, p394). Noticing that $nb-\frac{n+1}{2}-
\lambda +\frac{n}{s}=-n \varepsilon <0$ and  $-\frac{n+1}{2}-
\lambda +\frac{n}{s}  <0,$ we have
$$\sum_{k=0 }^\infty (J_2^k)^{1/s} =\sum_{k+1\geq  k_0 }(J_2^k)^{1/s}
\leq C 2^{k_0(nb -(n+1)/2-\lambda+n/s)} \| h \|_{L^{s} } \leq C
2^{k_0nb } \| h \|_{L^{s} }
$$
 for $k_0> 0$.
Thus,
\begin{eqnarray*}
\|T^{0+}h(x)|x |^{nb}\|^{ }_{L^{s}   } &\leq & \sum_{k=0 }^\infty
[(J_1^k)^{1/s}+  (J_2^k)^{1/s}]
 \\ & = & \sum_{k=0 }^\infty
(J_1^k)^{1/s}+ \sum_{k=0 }^\infty (J_2^k)^{1/s}
\leq  C
2^{k_0nb}2^{k_0n(1/s-1/p-\alpha/np)}.
\end{eqnarray*}
And by the $L^s$ boundedness of $T^{0+}$, see \cite{S3}, we have
$$
\|T^{0+}h( \|^{ }_{L^{s}   } \leq  C \| h \|^{ }_{L^{s}   }
 \leq  C  2^{k_0n(1/s-1/p-\alpha/np)}.
$$
Then, noticing $a/b>0$, $1-a/b>0$ and $b-a +\frac{1}{s}
-\frac{1}{p}-\frac{\alpha}{np}=0,$ we have
\begin{eqnarray*}
\Re (T^{0+}h)&\leq& C2^{k_0n(1/s-1/p-\alpha/np)a/b}
2^{k_0n(b+1/s-1/p-\alpha/np)(1-a/b)}
\\
&=&C2^{k_0n(b-a+1/s-1/p-\alpha/np)} = C,
\end{eqnarray*}
where $k_0>0$. Thus, (5.7) holds, and it follows that (5.2) holds
for all $(p,s,\alpha)$-blocks $h$ of restrict I-type.

\par For  $(p,s,\alpha)$-blocks, by (5.6) and Theorem 4.4, we see that for
$T$ in (5.6), the same conclusions hold as those stated in
Proposition 5.1.
 Then, by Proposition 5.3, we have

\par {\bf Proposition 5.4} Let $
  \lambda,s,p $ and $ \alpha $ under the conditions of Theorem 5.1. Then
  (5.2) holds
 for all $(p,s,\alpha)$-blocks $h$.

\par{\bf Step 3}.
\par We prove (5.3). From (5.3) the conclusion required can be proved easily.
\par {\bf Proposition 5.5}    Let $
  \lambda,s,p,\alpha $ under the conditions of Theorem 5.1,  $T$ as in
  (5.6). Then
$$
T ~{\rm is ~ of ~type ~} (
BL^{p,s}_{|x|^{\alpha}},BL^{p,s}_{|x|^{\alpha}})  .
$$

\par {\bf Proof}   We need to
prove  (5.3).
  For $f\in BL^{p,s}_{|x|^{\alpha}}$ and any
$0<\varepsilon \leq \|f\|_{BL^{p,s}_{|x|^{\alpha}}}$, there exists a
decomposition of $f$: $f=\sum_{j}\mu_j a_j $
 such that
$(\sum_{j}|\mu_j|^{\bar{p}})^{1/\bar{p}}\leq \|f\|
_{BL^{p,s}_{|x|^{\alpha}}}+\varepsilon, $
 where all $a_j$ are $(p,s,\alpha)$-blocks.
Once (5.3) is proved, then, by Proposition 5.4,  we have
\begin{eqnarray*}
\|  T  f\|_{BL^{p,s}_{|x|^{\alpha}}} & \leq& \sum_{j }|\mu_j|\| T
a_j\|_{BL^{p,s}_{|x|^{\alpha}}} \leq C\sum_{j }|\mu_j|
 \\ &\leq& C
\left(\sum_{j }|\mu_j|^{\bar{p}}\right)^{1/\bar{p}} \leq
C\|f\|_{BL^{p,s}_{|x|^{\alpha}}}+\varepsilon \leq
C\|f\|_{BL^{p,s}_{|x|^{\alpha}}}.
\end{eqnarray*}

 \par  Next, let us prove (5.3). It is clear that we need only to
 prove (5.3) for $T^{\frac{1}{2}+}$ since the cases for $T^{(2j+\frac{1}{2}) \pm }$ and
$T^{(2j+\frac{3}{2}) \pm }$, $j=0,1,2, \cdots, N_0$, are similar. As
(5.8), we can write

  $$T^{\frac{1}{2}+}=\sum_{k=0}^\infty T_k^{\frac{1}{2}+}\eqno(5.9)$$
 with
$T_k^{\frac{1}{2}+}f(x)=2^{-[(n+1)/2+\lambda+1/2]k} \int e^{ 2\pi i
|y|} f(x-y)\psi(y/2^k)dy.$

 \par Let us first consider $T_k^{\frac{1}{2}+}$. Fix $k$ and $x\in{\bf R}^n$.  Let
$f=\sum_{j=1}^{\infty}\mu_j a_j $ where each $a_j$ is a
$(p,s,\alpha)$-block with supp $a_j\subset B_j$ and $\|a_j\|_{L^s}
\leq |B_j|^{-\frac{1}{p}-\frac{\alpha}{np}+\frac{1}{s}}$, and
$(\sum_{j=1}^{\infty}|\mu_j|^{\bar{p}})^{1/\bar{p}}<\infty $.

\par When $a_j$ is a
$(p,s,\alpha)$-block of restrict II-type,   by H\"{o}lder
inequality, we have
$$
\int |a_j(x-y)||\psi(y/2^k)|dy\leq \|\psi\|_{L^\infty}\|a_j\|_{L^1}
\leq \|\psi\|_{L^\infty} |B_j|^{1-\frac{1}{p}-\frac{\alpha}{np}}
\leq \|\psi\|_{L^\infty},
$$
since $|B_j|\leq 1 $ and $ 1-\frac{1}{p}-\frac{\alpha}{np}>0.$

\par  When $a_j$ is a $(p,s,\alpha)$-block of restrict I-type,
  by H\"{o}lder inequality, we have
$$
\int |a_j(x-y)||\psi(y/2^k)|dy\leq
\|2^{nk/s'}\psi\|_{L^{s'}}\|a_j\|_{L^s} \leq 2^{nk/s'}
\|\psi\|_{L^{s'}} |B_j|^{\frac{1}{s}-\frac{1}{p}-\frac{\alpha}{np}}
\leq 2^{nk/s'}\|\psi\|_{L^{s'}},
$$
since $|B_j|> 1$ and $ \frac{1}{s}-\frac{1}{p}-\frac{\alpha}{np}\leq
0.$

Let $F(y)=\sum_{j=1}^{\infty}|\mu_j| |a_j(x-y)||\psi(y/2^k)|$ for
$x\in {\bf R}^n $, by Minkowski inequality and  the two inequalities
above,
 we have
$$
\| F\|_{L^1}\leq
(\|\psi\|_{L^\infty}+2^{nk/s'}\|\psi\|_{L^{s'}})\sum_{j=1}^{\infty}|\mu_j|\leq
(\|\psi\|_{L^\infty}+2^{nk/s'}\|\psi\|_{L^{s'}})(\sum_{j=1}^{\infty}|\mu_j|^{\bar{p}})^{1/\bar{p}}<\infty.
$$
 Let $F_N(y)=\sum_{j=1}^{N}\mu_j  e^{ 2\pi i |y|} a_j(x-y)\psi(y/2^k)$ for
$x\in {\bf R}^n $. It is clear that $|F_N(y)|\leq |F(y)|$. Then the
Lebesgue dominated convergence theorem gives that
$
 \int \lim _{N\rightarrow \infty }F_N(y)dy=  \lim _{N\rightarrow \infty }\int F_N(y)dy
$,
  and it follows that
 $$
T_k^{\frac{1}{2}+}f(x)=\sum_{j=0}^\infty \mu_jT_k^{\frac{1}{2}+}
a_j(x). \eqno(5.10) $$
 We see also that
$$|T_k^{\frac{1}{2}+}a(x)|\leq 2^{-[ \frac{n+1}{2} +\lambda +\frac{1}{2}]k}
(\|\psi\|_{L^\infty}+2^{\frac{nk}{s'}}\|\psi\|_{L^{s'}}) \eqno(5.11)
$$
for all $(p,s,\alpha)$-blocks with $1<s<\infty, 0<p\leq s$ and $
n(\frac{p}{s}-1)\leq \alpha <n(p-1)$. Noticing (5.10) and (5.9), we
have
$$T^{\frac{1}{2}+}f
=\sum_{k=0}^\infty  \sum_{j=0}^\infty \mu_jT_k^{\frac{1}{2}+} a_j
=\sum_{j=0}^\infty  \sum_{k=0}^\infty \mu_jT_k^{\frac{1}{2}+} a_j
=\sum_{j=0}^\infty \mu_j \sum_{k=0}^\infty T_k^{\frac{1}{2}+} a_j
=\sum_{j=0}^\infty \mu_j T^{\frac{1}{2}+} a_j,
$$
the second "=" above holds since
\begin{eqnarray*}
\sum_{j=0}^\infty  \sum_{k=0}^\infty |\mu_j||T_k^{\frac{1}{2}+} a_j|
 & \leq &
(\|\psi\|_{L^\infty}+\|\psi\|_{L^{s'}}) \sum_{j=0}^\infty
\sum_{k=0}^\infty 2^{-[\frac{n+1}{2} +\lambda
+\frac{1}{2}-\frac{n}{s'}]k} |\mu_j|
\\
&\leq & (\|\psi\|_{L^\infty}+\|\psi\|_{L^{s'}})(\sum_{k=0}^\infty
2^{-[\frac{n+1}{2} +\lambda
+\frac{1}{2}-\frac{n}{s'}]k})(\sum_{j=1}^{\infty}|\mu_j|^{\bar{p}})^{1/\bar{p}}<\infty,
\end{eqnarray*}
these inequalities  hold because of  (5.11) and
$(n+1)/2+\lambda+1/2-n/s'>(n+1)/2+\lambda-n/s'>0$  which follows
from $s<p_\lambda$.
 Thus, we have proved that $T$ is of
type $(BL^{p,s}_{|x|^{\alpha}},BL^{p,s}_{|x|^{\alpha}})$.
  Thus, we finish the proof of Proposition 5.5.
\par  The proof of Theorem 5.1 has been finished.

\par {\bf Proof of Corollary 5.1}  Let   $
 n(\frac{p}{p_\lambda}-1)< \alpha <n(\frac{p}{p'_\lambda}-1)$, then there
 exists $s$ such that $ p'_\lambda
<s<p_\lambda$ and $n(\frac{p}{p_\lambda}-1)<n(\frac{p}{s}-1)<
\alpha<n(\frac{p}{p'_\lambda}-1)$. Then, by Theorem 5.2, we have
that $ B_R^\lambda ~{\rm is ~ of ~type ~}
(BL^{p,s}_{|x|^{\alpha}},L^{p} _{|x|^{\alpha}})   $.  By Theorem
2.5, we know $BL^{p,p_\lambda}_{|x|^{\alpha}}\subset
BL^{p,s}_{|x|^{\alpha}}$
 for $p,s, p_\lambda$ and $\alpha$. So, $ B_R^\lambda ~{\rm is ~ of ~type ~}
(BL^{p, p_\lambda}_{|x|^{\alpha}},L^{p} _{|x|^{\alpha}})   $.

\par {\bf Proof of Theorem 5.2  } Since
\begin{eqnarray*}
B_R^\lambda f(x)=\int_{|\xi|<R}\hat{f}(\xi)\left(1-\left|
\xi/R\right|^2\right)^\lambda e^{2\pi ix\cdot\xi}d\xi,
\end{eqnarray*}
and let $g(x)=f(x/R),$ we have
\begin{eqnarray*}
B_R^\lambda
f(x)=\int_{|\eta|<1}\hat{g}(\eta)\left(1-|\eta|^2\right)^\lambda
e^{2\pi iRx\cdot\eta}d\eta.
\end{eqnarray*}

\par Let $f \in BL^{p,s}_{|x|^{\alpha}}$ and $f\neq 0$, then there exists a decomposition of $f(x)$:
 $f(x)=\sum_{i=1}^\infty \mu _ia_i(x)$ with supp $a_i\subseteq
B_i$, $\|a_i\|_{L^s}\leq
|B_i|^{-\frac{1}{p}-\frac{\alpha}{pn}+\frac{1}{s}}$ and
$(\sum_{j=1}^{\infty}|\mu_j|^{\bar{p}})^{1/\bar{p}}\leq 2
\|f\|^p_{BL^{p,s}_{|x|^{\alpha}}}.$
 Then
$g(x)=\sum_{i=1}^\infty \mu _ia_i(x/R)=\sum_{i=1}^\infty \mu
_ib_i(x)$ with supp $b_i =$ supp $a_i(\cdot/R) \subseteq RB_i$ and
$\|b_i\|_{L^s}\leq
R^{\frac{n+\alpha}{p}}|RB_i|^{-\frac{1}{p}-\frac{\alpha}{pn}+\frac{1}{s}}$,
and this shows that each $R^{-\frac{n+\alpha}{p}} b_i$ is a
$(p,s,\alpha)$-block. So we have
$$g(x) =\sum_{i=1}^\infty \mu
_iR^{\frac{n+\alpha}{p}}R^{-\frac{n+\alpha}{p}}b_i(x),$$
 and
 \begin{eqnarray*}
\|g\|^p_{BL^{p,s}_{|x|^{\alpha}}}  \leq  \left(\sum_{i=1}^\infty
\left(|\mu
_i|R^{\frac{n+\alpha}{p}}\right)^{\bar{p}}\right)^{1/\bar{p}} \leq 2
R^{\frac{n+\alpha}{p}}\|f\|^p_{BL^{p,s}_{|x|^{\alpha}}}.
\end{eqnarray*}

\par By Theorems 5.1 and   2.1, we see that
  $$
B_1^\lambda ~{\rm is ~ of ~type ~} ( BL^{p,s}_{|x|^{\alpha}},L^{p}
_{|x|^{\alpha}})
$$
for $    p'_\lambda <s<p_\lambda, 0<p< s , n(\frac{p}{s}-1)<
\alpha<n(\frac{p}{p'_\lambda}-1)$ and $\alpha\leq 0$. Using these,
we have
\begin{eqnarray*}
\int_{{\bf R}^n} |B_R^\lambda f(x)|^p|x|^\alpha dx &=&\int_{{\bf
R}^n}\left|\int_{|\eta|<1}\hat{g}(\eta)\left(1-|\eta|^2\right)^\lambda
e^{2\pi iRx\cdot\eta}d\eta\right|^p|x|^\alpha dx
\\
&=& \frac{1}{R^{\alpha+n}}\int_{{\bf
R}^n}\left|\int_{|\eta|<1}\hat{g}(\eta)\left(1-|\eta|^2\right)^\lambda
e^{2\pi ix\cdot\eta}d\eta\right|^p|x|^\alpha dx
\\ &\leq &
\frac{C}{R^{\alpha+n}}\|g\|^p_{BL^{p,s}_{|x|^{\alpha}}}
 \leq  C
\|f\|^p_{BL^{p,s}_{|x|^{\alpha}}}.
\end{eqnarray*}
Thus, we  finish the proof of Theorem 5.3.

\par In particular, taking  $\lambda=\frac{n-1}{2}$ in Theorem 5.1, we have

\par {\bf Corollary  5.2} Let  $ 0 <s<\infty, 0<p\leq s $ and $
n(\frac{p}{s}-1)\leq \alpha<n(p-1)$.  Then
$$
B_1^{\frac{n-1}{2}} {~\rm is ~ of ~ type ~ }(
BL^{p,s}_{|x|^{\alpha}},  BL^{p,s}_{|x|^{\alpha}})
   .
$$

\par When $p>1$,  taking $\alpha=0$ in Theorem 5.2, we have

\par {\bf Corollary 5.3}   (i) Let $
\frac{n-1}{2(n+1)}< \lambda\leq \frac{n-1}{2},
p_\lambda'<s<p_\lambda , 1<p< s $ and $
   \frac{n}{p}<\frac{1}{p_\lambda'} $. Then
 $$
 B_1^\lambda {~\rm is ~ of ~ type ~ }( BL^{p,s}_{|x|^{0}},L^{p}).
$$

\par When $p>1$  and  taking $\alpha=0$ in Theorem 5.3, we have

\par {\bf Corollary 5.4}  Let $0<R<\infty.$ For
$B_R^\lambda$, the same conclusions hold as those stated in
Corollary 5.3.

\par {\bf Remark 5.2} Since the conditions $p_\lambda'
<s$ and $n(\frac{p}{s}-1)< \alpha$ in theorems above imply that $
\frac{1}{s}-\frac{1}{p} -\frac{\alpha}{np} \leq 0<
\frac{1}{p_\lambda'}-\frac{1}{p} -\frac{\alpha}{np}
  $, those theorems above maybe fail when
  $s=p_\lambda'$.

\par {\bf Proof of Theorem 5.3 and Theorem 5.4} It is easy to see that
   Theorem 5.3 follows from the following
 two propositions, and Theorem 5.4 follows from Theorems 5.3 and
   2.1.

\par {\bf Proposition 5.6} Let
 $    1< s < \infty,0< p \leq s$ and $ - n(1-p/s)\leq \alpha < n(p-1)$.
We have
$$
\| B_*^{\frac{n-1}{2}} h\|_{BL^{p,s}_{|x|^{\alpha}}}\leq C
\eqno(5.12)
 $$
 for every $(  p,s,\alpha)-$block $h$,
  where the constant $C$ is  independent of $h$.

 \par {\bf Proof } Let $0<R<\infty$ and  $0<\lambda<\infty$. We know that
$ B_R^\lambda f(x)=(f\ast K_R^\lambda )(x) $ where
\begin{eqnarray*}
K^\lambda_R(x)
&=&
[(1-|\xi/R|^2)^\lambda_+]~\breve{}~(x)
\\&=&
\int_{|\xi|<R} \left(1-\left| \xi/R\right|^2\right)^\lambda e^{2\pi
ix\cdot\xi}d\xi
\\&=&
R^n K^\lambda_1(Rx)
\end{eqnarray*}
$$~~~~~~= R^n |Rx|^{-\frac{n}{2}-\lambda}
J_{\frac{n}{2}+\lambda}(2\pi|Rx|). \eqno(5.13)
$$
Here $J_{\frac{n}{2}+\lambda}$ is the Bessel function, and it
  has the asymptotic properties as follows:
$$
J_m(t)=\sqrt{\frac{2}{\pi}}~t^{-\frac{1}{2}}\cos (t-\frac{m\pi}{2}-\frac{\pi}{4})+O(t^{-\frac{3}{2}}),
 ~~~~{\rm ~ as ~} t\rightarrow \infty, \eqno(5.14)
 $$
$$
J_m(t)= \frac{t^m}{2^m\Gamma(m+1)}+O(t^{m+1}),
 ~~~~{\rm ~ as ~} t\rightarrow 0, \eqno(5.15)
 $$
where $m>-\frac{1}{2},$ see \cite {S3}. From (5.13) and (5.15), we
see that the radial function $K^\lambda_R(x)$ is right-continues at
$|x|=0$ for fixed $0<R<\infty.$ From (5.13) and (5.14), we see that
$$K^\lambda_R(x)=R^n O(|Rx|^{-\frac{n+1}{2}-\lambda})~~~~{\rm ~ as ~} t\rightarrow \infty. $$
Such two facts imply that
$$|K^\lambda_R(x)|\leq C R^n (1+|Rx|)^{-\frac{n+1}{2}-\lambda},~~~~  0<|x|< \infty,  \eqno(5.16)$$
where $C$ is independent of $|Rx|$. (5.16) shows that
$$|K^{\frac{n-1}{2}}_R(x)|\leq C    | x|)^{-n},~~~~  0<|x|< \infty, 0<R<\infty,   $$
 and therefore
$$|B_*^{\frac{n-1}{2}}f(x)|\leq C \int_{{\bf R}^n}\frac{|f(y)|}{|x-y|^n}dy.$$
As before (i.e. as (5.5)), it is easy to check that
$B_*^{\frac{n-1}{2}}$ satisfies (1.14). And it is known that
$B_*^{\frac{n-1}{2}}$ is bounded on $L^s$ with $1<s<\infty$, see
\cite{SG}.
  Then, by Theorem 4.1, (5.12)
follows.

\par {\bf Proposition 5.7} Let
 $    0<\lambda <\infty, 0<R<\infty, 1< s < p_\lambda,
 0< p \leq s$ and $ - n(1-p/s)\leq \alpha < n(p-1)$.
We have
$$
B_R^{\lambda}f(x)=\sum_{j=0}^\infty \mu_jB_R^{\lambda} a_j(x),
\eqno(5.17) $$ and
$$
B_*^{\lambda}f(x)\leq \sum_{j=0}^\infty |\mu_j|B_*^{\lambda} a_j(x)
\eqno(5.18) $$
 for   $f=\sum\lambda_ja_j \in \dot{B}L^{p,s}
_{|x|^{\alpha}}$, where each $a_j$ is a   $(  p,s, \alpha)-$block
and $\sum|\lambda_j|^{\bar{p}}<\infty$.

\par {\bf Proof } Fix $0<R<\infty$ and $x\in{\bf R}^n$.  Let
$f=\sum_{j=1}^{\infty}\mu_j a_j $ where each $a_j$ is a
$(p,s,\alpha)$-block with supp $a_j\subset B_j$ and $\|a_j\|_{L^s}
\leq |B_j|^{-\frac{1}{p}-\frac{\alpha}{np}+\frac{1}{s}}$, and
$(\sum_{j=1}^{\infty}|\mu_j|^{\bar{p}})^{1/\bar{p}}<\infty $.

\par When $a_j$ is a
$(p,s,\alpha)$-block of restrict II-type,   by (5.16) and H\"{o}lder
inequality, we have
$$
\int |a_j(x-y)||K^\lambda_R(y)|dy\leq CR^n\|a_j\|_{L^1} \leq CR^n
|B_j|^{1-\frac{1}{p}-\frac{\alpha}{np}} \leq CR^n \eqno(5.19)
$$
since $|B_j|\leq 1, 1-\frac{1}{p}-\frac{\alpha}{np}>0$ and
$\frac{n+1}{2}+\lambda>0$.

\par  When $a_j$ is a $(p,s,\alpha)$-block of restrict I-type,
  by H\"{o}lder inequality, we have
$$
\int |a_j(x-y)||K^\lambda_R(y)|dy\leq
\|K^\lambda_R\|_{L^{s'}}\|a_j\|_{L^s} \leq \|K^\lambda_R\|_{L^{s'}}
|B_j|^{\frac{1}{s}-\frac{1}{p}-\frac{\alpha}{np}} \leq
\|K^\lambda_R\|_{L^{s'}} \eqno(5.20)
$$
since $|B_j|> 1$ and $ \frac{1}{s}-\frac{1}{p}-\frac{\alpha}{np}\leq
0.$ By (5.16), we see that
\begin{eqnarray*}
\|K^\lambda_R\|^{s'}_{L^{s'}} &\leq & C \int _{{\bf R}^n}
\frac{R^{ns'}dx}{ (1+|Rx|)^{(\frac{n+1}{2}+\lambda)s'}}
\\
&= & C R^{n(s'-1)}\int_0^\infty  \frac{r^{n-1}dr}{
(1+r)^{(\frac{n+1}{2}+\lambda)s'}}
\\
&\leq & C R^{n(s'-1)}\int_0^\infty
(1+r)^{n-1-(\frac{n+1}{2}+\lambda)s'}dr
\end{eqnarray*}
$$\leq  C R^{n(s'-1)}~~~~~~~~~~~~~~~~~~~~~~~ \eqno(5.21)
$$
as $n-(\frac{n+1}{2}+\lambda)s'<0$, i.e. $s<p_\lambda $.

Let $F(y)=\sum_{j=1}^{\infty}|\mu_j| |a_j(x-y)||K^\lambda_R(y)|$ for
$x\in {\bf R}^n $, by Minkowski inequality and (5.19) $\sim$ (5.21),
 we have
$$
\| F\|_{L^1}\leq C( R^{n}+
R^{n(s'-1)})\sum_{j=1}^{\infty}|\mu_j|\leq C( R^{n}+
R^{n(s'-1)})(\sum_{j=1}^{\infty}|\mu_j|^{\bar{p}})^{1/\bar{p}}<\infty.
$$
 Let $F_N(y)=\sum_{j=1}^{N}\mu_j    a_j(x-y)K^\lambda_R(y)$ for
$x\in {\bf R}^n $. It is clear that $|F_N(y)|\leq |F(y)|$. Then the
Lebesgue dominated convergence theorem gives (5.17).

\par And it follows that
$$
|B_R^{\lambda}f(x)|\leq \sum_{j=0}^\infty |\mu_j|B_*^{\lambda}
a_j(x)
 $$
  for $0<R<\infty$. Thus we have  (5.18).

\par {\bf Proof of Theorem 5.5 and Theorem 5.6}  When $\lambda>\frac{n-1}{2}$, we know that
$ |B_*^\lambda f(x)|\leq C Mf(x) $ for all $f\in L^1_{loc}$ and an
absolutely constant $C$ (see, for example, \cite{SS}). Then Theorems
5.5 and  5.6 follow from  Theorem 4.4.

\par {\bf Proof of Theorem 5.7}
 Let
$f=\sum_{j=1}^{\infty}\mu_j a_j $ where each $a_j$ is a
$(p,s,\alpha)$-block  and
$(\sum_{j=1}^{\infty}|\mu_j|^{\bar{p}})^{1/\bar{p}}<\infty $. Then,
for any $\varepsilon>0$, there exists an $N_\varepsilon>0$ such that
$(\sum_{j=N_\varepsilon+1}^{\infty}|\mu_j|^{\bar{p}})^{1/\bar{p}}<\varepsilon
$.
 By (5.17), we have
\begin{eqnarray*}
& &\|B_R^\lambda f-f\|_{L^{p} _{|x|^{\alpha}}}
\\&=&\|\sum_{j=1}^{\infty}\mu_j B_R^\lambda
a_j-\sum_{j=1}^{\infty}\mu_j a_j\|_{L^{p} _{|x|^{\alpha}}}
\\
&\leq &
 \sum_{j=1}^{N_\varepsilon}|\mu_j|\| B_R^\lambda a_j-
a_j\|_{L^{p} _{|x|^{\alpha}}} +
 \sum_{j=N_\varepsilon+1}^{\infty}|\mu_j|\|B_R^\lambda a_j \|_{L^{p} _{|x|^{\alpha}}}
+
 \sum_{j=N_\varepsilon+1}^{\infty}|\mu_j| \|  a_j \|_{L^{p} _{|x|^{\alpha}}}
\\&=& I+II+III.
\end{eqnarray*}
By Theorem 5.2 we have $\|B_R^\lambda a_j \|_{L^{p}
_{|x|^{\alpha}}}\leq \| a_j \|_{BL^{p,s}_{|x|^{\alpha}}}\leq C$, and
by Theorem 2.1 we have $\|  a_j \|_{L^{p} _{|x|^{\alpha}}}\leq \|
a_j \|_{BL^{p,s}_{|x|^{\alpha}}}\leq C$. Then,
$$
II+III\leq C \sum_{j=N_\varepsilon+1}^{\infty}|\mu_j| \leq C
\left(\sum_{j=N_\varepsilon+1}^{\infty}|\mu_j|^{\bar{p}}\right)^{1/\bar{p}}<C\varepsilon,
\eqno(5.22)
$$
where $C$ is an absolutely constant. For $I$, since $\|B_R^\lambda
a_j \|_{L^{p} _{|x|^{\alpha}}} \leq C<\infty $ and   $\|  a_j
\|_{L^{p} _{|x|^{\alpha}}}  \leq C<\infty$ for $j=1,2,\cdots,
N_\varepsilon$,  there exists a ball $B_\varepsilon$ such that
$$
\int_{B^c_\varepsilon}|B_R^\lambda a_j(x)|^p|x|^{\alpha}dx
<N_\varepsilon^{-p}\varepsilon^p{\rm ~~and~~ }
\int_{B^c_\varepsilon}|
a_j(x)|^p|x|^{\alpha}dx<N_\varepsilon^{-p}\varepsilon^p \eqno(5.23)
$$
for all $j=1,2,\cdots, N_\varepsilon$. It is known that $
\|B_R^\lambda a_j-a_j\|_{L^s } \rightarrow 0$ as $ R\rightarrow
\infty $, $j=1,2,\cdots, N_\varepsilon$, then for the
$N_\varepsilon^{-1}|B_\varepsilon|^{-(\frac{1}{p}-\frac{1}{s}-\frac{\alpha}{np})}\varepsilon$
, there exists an $R_0$ such that $ \|B_R^\lambda a_j-a_j\|_{L^s}
<N_\varepsilon^{-1}|B_\varepsilon|^{-(\frac{1}{p}-\frac{1}{s}-\frac{\alpha}{np})}\varepsilon
$ for all $R>R_0$ and  $j=1,2,\cdots, N_\varepsilon$. Then, by
H\"{o}lder inequality,
\begin{eqnarray*}
\int_{B_\varepsilon}|B_R^\lambda a_j(x)-a_j(x)|^p|x|^{\alpha}dx \leq
\left(\int_{B_\varepsilon}|B_R^\lambda
a_j(x)-a_j(x)|^sdx\right)^{p/s}
\left(\int_{B_\varepsilon}|x|^{\alpha s/(s-p)}dx\right)^{(s-p)/s}
\end{eqnarray*}
$$~~~~~~~~~~~~
  \leq     C
|B_\varepsilon|^{1-\frac{p}{s}-\frac{\alpha}{n}}
 \|B_R^\lambda a_j-a_j\|^p_{L^s} < CN_\varepsilon^{-p}\varepsilon^p
 \eqno(5.24)
 $$
since $\alpha s/(s-p) +n\leq 0$, where $C$ is an absolutely
constant. Thus, (5.23) and (5.24) give that $\int_{}|B_R^\lambda
a_j(x)-a_j(x)|^p|x|^{\alpha}dx <CN_\varepsilon^{-p}\varepsilon^p $,
and it follows that
$$I<C\varepsilon\eqno(5.25)$$
for all $R>R_0$. Then, by (5.22) and (5.25), we have that
$\|B_R^\lambda f-f\|_{L^{p} _{|x|^{\alpha}}}<C\varepsilon$ for all
$R>R_0$. Thus, we finish the proof of Theorem 5.7.

\par \section*
 {\bf  6. Spherical means }
We prove the following
\par {\bf Theorem 6.1}Let
 $  n\geq 3 ,0< p \leq 2 $ and $ n(\frac{1}{2}p-1) \leq \alpha < n(\frac{n-1}{n}p-1)$.
 Then (1.10) holds.






\par From Theorems 2.1 and  6.1, we have

\par {\bf Theorem 6.2} Let
 $  n\geq 3 ,0< p <2 , n(\frac{1}{2}p-1) < \alpha < n(\frac{n-1}{n}p-1)$ and $\alpha\leq 0$. Then (1.11) holds.
\par  In particularly, we have
\par {\bf Corollary 6.1} Let
 $ n\geq 3,0<p\leq \frac{n}{n-1}$ and $ n(\frac{1}{2}p-1)< \alpha < n(\frac{n-1}{n}p-1)$. Then
 $${\mathcal{M}} {\rm~ is~ of~ type~} (BL^{p,2}_{|x|^{\alpha}}, {L}^{p}_{|x|^{\alpha}}  ).$$

\par {\bf Theorem  6.3}  Under the hypothesis of Theorem 6.2, we
have
$$
\|A_t f-f\|_{L^p_{|x|^\alpha}} \rightarrow 0,~~ t\rightarrow 0
$$
for $f\in BL^{p,2}_{|x|^{\alpha}}$.

\par As $B_R^\lambda
f$, we can also get that $A_t f$ converges to $f$,
$\mu_\alpha$-a.e.,
 as $t$ tends to $0$
 by Theorem 6.1 for $f\in \dot{B}L^{p,2}
_{|x|^{\alpha}}$  with $ p $ and $
  \alpha $ under the conditions of Theorem 6.1, but
  this does not imply new information since $BL^{p,2}_{|x|^{\alpha}}\subset
  L^{\frac{np}{n+\alpha}}$, and  $\frac{np}{n+\alpha}>\frac{n}{n-1}$.

 \par We will use the following known facts in the proof of Theorem 6.1.
  Let ${\mathcal{M}} (f)$ be  the the spherical
 maximal function defined in section 1,  we have then
$${\mathcal{M}} (f)\leq \sum_{j=0}^{\infty}{\mathcal{M}} _j(f),\eqno(6.1)$$
where
$${\mathcal{M}} _j(f)\leq C(n)2^j M(f), ~~j=0,1,2,\cdots, \eqno(6.2)$$
and
$$\|{\mathcal{M}} _j(f)\|_{L^2}\leq C(n)2^{(1-n/2)j} \|f\|_{L^2}, ~~j=0,1,2,\cdots, \eqno(6.3)$$
for all function $f$ in $L^2=L^2({\bf R}^n)$. Here $M(f)$ is
Hardy-Littlewood maximal function. See \cite{Grafakos}.

\par {\bf Proof of Theorem 6.1 } By Theorem 3.1, it  suffices to show that
$$\|{\mathcal{M}} h\|_{ BL^{p,2}_{|x|^{\alpha}}}\leq
C \eqno(6.4)$$
 for any $( p,2, \alpha)-$block $h$, where  constant $C$ is
independent of $h$.
\par In fact,  let $f\in  BL^{p,2}_{|x|^{\alpha}}$, then for any $\varepsilon >0,$ there exists
$f=\sum\lambda_ih_i$ where each $h_i$ is a $(  p,2, \alpha)-$block
such that
 $$\left( \sum
|\lambda_i|^{\bar{p}}\right)^{1/\bar{p}}\leq \|f\|_{
BL^{p,2}_{|x|^{\alpha}}}+\varepsilon. $$ By Minkowski inequality,
$|A_tf(x)|\leq \sum|\lambda_i|A_t(|h_i|)(x)\leq
\sum|\lambda_i|{\mathcal{M}}(|h_i|)(x)$, and it follows that
$${\mathcal{M}}f(x)|\leq \sum|\lambda_i|{\mathcal{M}}(|h_i|)(x).$$
 Then, once (6.4) is established, by  the subadditivity of the norm,
we have
 $$\|{\mathcal{M}} f\|_{BL^{p,2}_{|x|^{\alpha}}}\leq
\sum |\lambda_i|\|{\mathcal{M}} h_i\|_{BL^{p,2}_{|x|^{\alpha}}} \leq
C \sum |\lambda_i| \leq C \left( \sum
|\lambda_i|^{\bar{p}}\right)^{1/\bar{p}} \leq
C\|f\|_{BL^{p,2}_{|x|^{\alpha}}}+\varepsilon.
 $$
Since $\varepsilon$ is arbitrary, we have $\|{\mathcal{M}}
f\|_{BL^{p,2}_{|x|^{\alpha}}} \leq C\|f\|_{BL^{p,2}
 _{|x|^{\alpha}}} .
 $

 \par In order to obtain (6.4), Let us estimate ${\mathcal{M}}_j h, j=0,1,2,\cdots$.

\par
We check that ${\mathcal{M}}_j h$ is a $ (p, 2, \alpha, \varepsilon
)$-molecular
 for every $(  p,2,\alpha)-$block $h$ centered at any
point and estimate $\Re ({\mathcal{M}}_j h) $.

\par  By $\alpha < n(\frac{n-1}{n}p-1)$, we have that
$$
\frac{n-1}{n-2}-\frac{n}{n-2}\left(\frac{1}{p}+\frac{\alpha}{np}\right)>0.\eqno(6.5)$$
By (6.5) and  noticing that  $ 1
 -1/p-\alpha/np>0$, we can choose $\varepsilon$ such that
$$0<\varepsilon <\min\left\{\frac{1}{2},1 -\frac{1}{p}-\frac{\alpha}{np},
\frac{n-1}{n-2}-\frac{n}{n-2}\left(\frac{1}{p}+\frac{\alpha}{np}\right)\right\}.\eqno(6.6)$$
 Set $
a=1-1/p-\alpha/np-\varepsilon $ and $ b=1/2-\varepsilon. $ We see
that $a>0,b>0$ from (6.6),  $b-a=-(1/2-1/p-\alpha/np)\geq 0$ from
$n(p/2-1)\leq \alpha$, and
$$
 na-2b=n-1-n(1/p+\alpha/np)-(n-2)\varepsilon>0 \eqno(6.7)
 $$
from (6.6).

 Given a $(  p,2, \alpha)-$block $h$ with supp $h\subset
B(x_0,r)\subset B(x_0,2^{k_0}),$ and
$$ \|h\|^{ }_{L^{2} }\leq
|B(x_0,r)|^{1/2-1/p-\alpha/np}\leq C_{\alpha, p}
|B(x_0,2^{k_0})|^{1/2-1/p-\alpha/np},$$ where
$2^{k_0-1}<r\leq2^{k_0}$, then we have
\begin{eqnarray*}
\|{\mathcal{M}}_jh(x)|x-x_0 |^{nb}\|^{2}_{L^{2}   } &=&\int _{{\bf
R}^n } |{\mathcal{M}}_jh(x)|^2 |x-x_0|^{2nb} dx
\\
& \leq &
 \left(\int _{|x-x_0|<2^{k_0+ 1}}+\int _{|x-x_0|\geq 2^{k_0+ 1}} \right)
|{\mathcal{M}}_jh(x)|^2 |x-x_0|^{2nb} dx
\end{eqnarray*}
$$
 =  J_1+J_2.~~~~~~~~~~~~~~~~~~~~~~~~~~~~~\eqno(6.8)$$

\par For $J_1,$ by (6.3), we have
\begin{eqnarray*}
J_1&=& \int _{|x-x_0|<2^{k_0+ 1}} |{\mathcal{M}}_jh(x)|^2
|x-x_0|^{2nb} dx
\\
& \leq & C 2^{2k_0nb}\|{\mathcal{M}}_jh \|^{2}_{L^{2}   }
\\
& \leq & C2^{2(1-n/2)j} 2^{2k_0nb}\| h \|^{2}_{L^{2}   }
\\
& \leq & C2^{2(1-n/2)j} 2^{2k_0nb}2^{2k_0n(1/2 -1/p-\alpha/np)}
\end{eqnarray*}

For $J_2,$ by (6.2), and noticing that  Hardy-Littlewood maximal
function $M$ satisfies  (1.14),  and  that $b-1+1/2=-\varepsilon
<0,$ we have
\begin{eqnarray*}
J_2&=& \int _{|x-x_0|\geq 2^{k_0+ 1}} |{\mathcal{M}}_jh(x)|^2
|x-x_0|^{2nb} dx
\\
& \leq & C 2^{2j} \int _{|x-x_0|\geq 2^{k_0+ 1}} |Mh(x)|^2
|x-x_0|^{2nb} dx
\\
& \leq & C2^{2j}\int _{|x-x_0|\geq 2^{k_0+ 1}} \|h \|^{2}_{L^{1}   }
|x-x_0|^{2n(b-1)} dx
\\
& \leq & C2^{2j}2^{2k_0n(b-1+1/2)}\| h \|^{2}_{L^{2}   }
|B(x_0,2^{k_0})|^{}
\\
& \leq & C 2^{2j}2^{2k_0nb}2^{2k_0n(1/2-1/p-\alpha/np)}.
\end{eqnarray*}
Thus,
$$
\|{\mathcal{M}}_jh(x)|x-x_0 |^{nb}\|^{ }_{L^{2}   } \leq  C 2^j
2^{k_0nb}2^{k_0n(1/s-1/p-\alpha/np)}.
$$
And by (6.3),
$$
\|{\mathcal{M}}_jh(x) \|^{ }_{L^{2}   } \leq  C 2^{(1-n/2)j} \| h(x)
\|^{ }_{L^{2} }
 \leq  C 2^{(1-n/2)j} 2^{k_0n(1/2-1/p-\alpha/np)}.
$$
Then, noticing $a/b>0$ and $1-a/b\geq 0,$ we have
\begin{eqnarray*}
\Re _1({\mathcal{M}}_j h)&\leq& C 2^{(1-\frac{n}{2})j\times
\frac{a}{b}} 2^{j(1-\frac{a}{b})} 2^{k_0n(1/2-1/p-\alpha/np)a/b}
2^{k_0n(b+1/2-1/p-\alpha/np)(1-a/b)}
\\
&=&C 2^{-j\times \frac{na-2b}{2b}} 2^{k_0n(b-a+1/2-1/p-\alpha/np)}
=C2^{-j\times \frac{na-2b}{2b}}.
\end{eqnarray*}
Thus, by Theorem 3.1, we have
$$\|{\mathcal{M}}_j h\|_{BL^{p,2}_{|x|^{\alpha}}}\leq C2^{-j\times \frac{na-2b}{2b}}.$$
By (6.1), it follows that
$$\|{\mathcal{M}} h\|_{BL^{p,2}_{|x|^{\alpha}}}
\leq \sum_{j=0}^\infty \|{\mathcal{M}}_j h\|_{BL^{p,2}
_{|x|^{\alpha}}}\leq C\sum_{j=0}^\infty2^{-j\times
\frac{na-2b}{2b}}=C$$ noticing (6.7) and $b>0$.

\par
Thus, we finish the proof of Theorem 6.1.

\par In order to prove Theorem 6.3, let us state a uniform boundedness
principle related to $BL^{p,s}_{|x|^{\alpha}}$, which is similar to
that related to $L^{p} $.

\par {\bf Lemma 6.1 (Uniform boundedness principle)}
Let $1<s<\infty, 0<p\leq s$ and $ -n<\alpha<n(p-1)$.
 Let $\mathcal{D}$ be a dense subspace of $
BL^{p,s}_{|x|^{\alpha}}$ and suppose that $T_R$ is a sequence of
linear operators such that
$$
\|T_Rf-f\|_{L^p_{|x|^\alpha}} \rightarrow 0, R\rightarrow
\infty,\eqno(6.9)
$$
 for $f\in \mathcal{D}$.
 Then in order for $T_Rf\rightarrow f$ in
$L^p_{|x|^\alpha}$ norm  as $R\rightarrow \infty $ for all $f\in
BL^{p,s}_{|x|^{\alpha}}$, it is a necessary and sufficient condition
that we have the estimate
$$
 \|T_Rf\|_{L^p_{|x|^\alpha}}\leq C\|f\|_{BL^{p,s}_{|x|^{\alpha}}}\eqno(6.10)
$$
for all sufficiently large $R,$ where the constants $C$ are
independent of $R$.

\par {\bf Proof }
   By density, for $f\in
BL^{p,s}_{|x|^{\alpha}}$ and any $\varepsilon>0$, there exists a
$g\in \mathcal{D}$ such that $\|f-g\|_{BL^{p,s}_{|x|^{\alpha}}}\leq
\varepsilon$. Noticing (6.9) and (6.10), we have
\begin{eqnarray*}
\|T_Rf-f\|^{\bar{p}}_{L^p_{|x|^\alpha}}&\leq &
\|T_Rf-T_Rg\|^{\bar{p}}_{L^p_{|x|^\alpha}}+
\|T_Rg-g\|^{\bar{p}}_{L^p_{|x|^\alpha}}+\|g-f\|^{\bar{p}}_{L^p_{|x|^\alpha}}
\\ &\leq & C \|f-g\|^{\bar{p}}_{BL^{p,s}_{|x|^{\alpha}}}+\varepsilon^{\bar{p}}
+\|g-f\|^{\bar{p}}_{BL^{p,s}_{|x|^{\alpha}}} \leq
C\varepsilon^{\bar{p}}
\end{eqnarray*}
for $f\in BL^{p,s}_{|x|^{\alpha}}$.
\par On the other hands, for $f\in
BL^{p,s}_{|x|^{\alpha}}$, set
$\varepsilon=\|f\|_{BL^{p,s}_{|x|^{\alpha}}}$.
 As $R$ larger enough, we have
\begin{eqnarray*}
\|T_Rf\|^{\bar{p}}_{L^p_{|x|^\alpha}}\leq
\|T_Rf-f\|^{\bar{p}}_{L^p_{|x|^\alpha}} +\|
f\|^{\bar{p}}_{L^p_{|x|^\alpha}} \leq  C \|f
\|^{\bar{p}}_{BL^{p,s}_{|x|^{\alpha}}}.
 \end{eqnarray*}
  Thus, we finish the proof of Lemma 6.1.

\par {\bf Proof of Theorem 6.3 } It is easy to see that
$L^2 \bigcap  BL^{p,2}_{|x|^{\alpha}}$ is dense in $
BL^{p,2}_{|x|^{\alpha}}$. For   $ f\in L^2 \bigcap
 BL^{p,2}_{|x|^{\alpha}}$,
 we know that $A_t f$ converges to $f$ pointwise (see \cite{S3}),
 noticing that $|A_t f(x)|\leq {\mathcal{M}} f(x)$ for all
 $0<t<\infty$ and
 $
\|{\mathcal{M}} f\|_{L^{p} _{|x|^{\alpha}}} \leq C
\|f\|_{BL^{p,s}_{|x|^{\alpha}}}<\infty $ by Theorem 6.2, then the
Lebesgue dominated convergence theorem gives that
$$
\|A_t f-f\|_{L^p_{|x|^\alpha}} \rightarrow 0, ~~~~R\rightarrow
\infty
$$
for all $ f\in L^2\bigcap BL^{p,2}_{|x|^{\alpha}}$. Then
 Theorem 6.3 follows from Lemma 6.1 and Theorem 6.2.

\par \section*
 {\bf  7. Some  sublinear operators}

In this section, we  prove that a lot of  sublinear operators in
harmonic analysis are of some new estimate at or beyond endpoint,
some of them, such as Hilbert transform, Riesz transforms and the
regular singular integral operators,  and the corresponding
truncated operators and maximal operators
  are of   $(BL^{p,s}_{|x|^{\alpha}}, BL^{p,s}_{|x|^{\alpha}})$
   and $(BL^{p,s}_{|x|^{\alpha}}, {L}^{p} _{|x|^{\alpha}})$ estimates, while
   others, such as Calder\'{o}n-Zygmund  operator,
    the strongly singular integral operator and some
oscillatory singular integrals operators,  admit respectively an
extension which is bounded on $BL^{p,s}_{|x|^{\alpha}}$.

\par {\bf Theorem 7.1} Let
 $    1< s < \infty,0< p \leq s$ and $ - n(1-p/s)\leq \alpha < n(p-1)$.
 Suppose a  sublinear operator  $T$ satisfies the size
condition
$$|Tf(x)|\leq C \int_{{\bf R}^n}\frac{|f(y)|}{|x-y|^n}dy, ~~~~ x\notin
{\rm supp}f, \eqno{(7.1)}$$
for any integral function with compact support. Then, if $T$ is
bounded on $L^s$,  we have
 $$\|T h\|_{BL^{p,s}_{|x|^{\alpha}}}\leq
C  \eqno(7.2)$$
 for every $(  p,s,\alpha)-$block $h$, where $C$ is independent of
 $h$.

\par {\bf Proof} As (5.5), it is easy to check that the condition
(7.1) implies the conditions (1.14) and (1.16). Then, Theorem 7.1
follows from Theorem 4.1.

\par {\bf Theorem 7.2} Let
 $    1< s < \infty,0< p \leq s$ and $ - n(1-p/s)\leq \alpha < n(p-1)$.
 Let
  $K\in L^1_{loc}$ which satisfies
$$|K(x)|\leq C 1/|x|^n,~~x\neq 0,\eqno(7.3)$$
 $K_\varepsilon=K\chi_{\{|x|>\varepsilon\}}$ with $\varepsilon>0$. Define
$Tf(x)=  K \ast f(x), T^\varepsilon f(x)=  K_\varepsilon \ast f(x)$
and $T^* f(x)=\sup_\varepsilon |T^\varepsilon f(x)|.$ Then, if $
T^\varepsilon$ is bounded on $L^s$,  we have
 $$
T^\varepsilon ~{\rm is ~ of ~type ~} (BL^{p,s}
_{|x|^{\alpha}},BL^{p,s}_{|x|^{\alpha}})  ,
$$
and if $ T^*$ is bounded on $L^s$,  we have
 $$
T^*, T^\varepsilon {\rm ~and ~}T ~{\rm are ~ of ~type ~}
(BL^{p,s}_{|x|^{\alpha}},BL^{p,s}_{|x|^{\alpha}})  .
$$

\par {\bf Proof} Let us first consider $T^\varepsilon$. We pick a smooth radial function $\varphi$ with
support inside the ball $B(0,2)$ which is equal to 1 on the closed
unit ball $\overline{B(0,1)}$. Let
$\psi_k(x)=\varphi(2^{-k}x)-\varphi(2^{-k+1}x)$ supposed in the
annuli $2^{k-1}\leq |x|\leq 2^{k+1}, k=1,2,\cdots.$ Then, we can
write $T^\varepsilon f$ as
$$
T^\varepsilon f= \sum_{k=0}^{\infty}T^\varepsilon_kf,\eqno{(7.4)}
$$
where $T^\varepsilon_0f(x)=  (K_\varepsilon\varphi) \ast f(x)$ and
$T^\varepsilon_kf(x)= (K_\varepsilon \psi_k) \ast f(x),
k=1,2,\cdots.$

 \par Fix  $x\in{\bf R}^n$.  Let
$f=\sum_{j=1}^{\infty}\mu_j a_j $ where each $a_j$ is a
$(p,s,\alpha)$-block with supp $a_j\subset B_j$, and $\|a_j\|_{L^s}
\leq |B_j|^{-\frac{1}{p}-\frac{\alpha}{np}+\frac{1}{s}}$, and
$(\sum_{j=1}^{\infty}|\mu_j|^{\bar{p}})^{1/\bar{p}}<\infty $.

\par It is clear that we can let $0<\varepsilon <2$. As Proposition 5.5, we have
$$
T_k ^\varepsilon f(x)=\sum_{j=1}^\infty \mu_jT_k a_j(x),
~~~~k=1,2,\cdots,\eqno{(7.5)}
$$
and
$$|T_k ^\varepsilon a(x)|\leq
2^{-nk/s}(\|\psi\|_{L^\infty}+\|\psi\|_{L^{s'}}),
~~~~k=1,2,\cdots,\eqno{(7.6)}
$$
for all $(p,s,\alpha)$-blocks with $1<s<\infty, 0<p\leq s$ and $
n(\frac{p}{s}-1)\leq \alpha <n(p-1)$. Similarly,
$$
T_0 ^\varepsilon f(x)=\sum_{j=1}^\infty \mu_jT_0 ^\varepsilon
a_j(x)\eqno{(7.7)}
$$
and
$$|T_0 ^\varepsilon a(x)|\leq
\varepsilon^{-n}(\|\varphi\|_{L^\infty}+\|\varphi\|_{L^{s'}})\eqno{(7.8)}
 $$
for all $(p,s,\alpha)$-blocks. From (7.6) and (7.8) we see that
\begin{eqnarray*}
\sum_{j=1}^\infty  \sum_{k=0}^\infty |\mu_j||T_k ^\varepsilon a_j|
 & \leq &
(\varepsilon^{-n}(\|\varphi\|_{L^\infty}+\|\varphi\|_{L^{s'}})+\|\psi\|_{L^\infty}+\|\psi\|_{L^{s'}})
\sum_{j=1}^\infty \sum_{k=0}^\infty 2^{-nk/s} |\mu_j|
\\
&\leq & C_\varepsilon
(\sum_{j=1}^{\infty}|\mu_j|^{\bar{p}})^{1/\bar{p}}<\infty.
\end{eqnarray*}
Then, from (7.4), (7.5) and (7.7), we have
$$T^\varepsilon f
=\sum_{k=0}^\infty  \sum_{j=1}^\infty \mu_jT_k ^\varepsilon a_j
^\varepsilon a_j =\sum_{j=1}^\infty  \sum_{k=0}^\infty \mu_jT_k
^\varepsilon a_j =\sum_{j=1}^\infty \mu_j T^{\varepsilon}
a_j.\eqno(7.9)
$$
 From this and Theorem 7.1, it is easy to get that $
T^\varepsilon ~{\rm is ~ of ~type ~}
(BL^{p,s}_{|x|^{\alpha}},BL^{p,s}_{|x|^{\alpha}}) $  if $
T^\varepsilon$ is bounded on $L^s$.

\par From (7.9), we see that
$$|T^\varepsilon f(x)|
\leq \sum_{j=1}^\infty |\mu_j|| T^{\varepsilon} a_j(x)| \leq
\sum_{j=1}^\infty |\mu_j| T^{*} a_j(x)
$$
for $0<\varepsilon <\infty$,  and it follows that
$$T^* f(x)
\leq \sum_{j=1}^\infty |\mu_j| T^{*} a_j(x).
$$
 By Theorem 7.1,  $
T^* ~{\rm is ~ of ~type ~}
(BL^{p,s}_{|x|^{\alpha}},BL^{p,s}_{|x|^{\alpha}}) $  if $ T^*$ is
bounded on $L^s$, and $T$ and $T^\varepsilon$ follow.    Thus, we
finish the proof of Theorem 7.2.

\par By Theorems 7.2 and 2.1, we have
\par {\bf Theorem 7.3}  Let
 $    1< s < \infty,0< p < s,  - n(1-p/s)< \alpha < n(p-1)$ and $\alpha\leq 0$. Let
 $      K, K_\varepsilon, T, T^\varepsilon$ and $ T^*$ as in Theorem
 7.2.
  Then, if $
T^\varepsilon$ is bounded on $L^s$,  we have
 $$
T^\varepsilon ~{\rm is ~ of ~type ~} (BL^{p,s}_{|x|^{\alpha}},L^{p}
_{|x|^{\alpha}})  ,
$$
and if $ T^*$ is bounded on $L^s$,  we have
 $$
T^*, T^\varepsilon {\rm ~ and ~}T ~{\rm are ~ of ~type ~}
(BL^{p,s}_{|x|^{\alpha}},L^{p} _{|x|^{\alpha}})  .
$$

\par  The conditions  of Theorem 7.2 or Theorem 7.1 are  satisfied by many operators arising
in harmonic analysis, for example,

\par{\bf A.}
 Hilbert transform and the corresponding truncated
operators and maximal operator: $$Hf(x)= \lim
_{\varepsilon\rightarrow 0}H_{\varepsilon}f(x),
H_{\varepsilon}f(x)=\int_{|x-t|>\varepsilon}\frac{f(t)}{x-t}dt,~
{\rm and }~ H^*f(x)=\sup_{\varepsilon>0 } |H_{\varepsilon}f(x)|.$$

\par It is well known that
$H,H_{\varepsilon}$ and $H^*$
 are of type $(L^p,L^p)$ for $1<p<\infty $. Then we see that
$H,H_{\varepsilon}$ and $H^*$ satisfy the conditions of Theorem 7.2.

\par It is known that
  $H,H_{\varepsilon}$ and $H^*$ are of weak-type
 $(L^1,L^1)$, and   $Hf$ is of types $(H^p,L^p)$ and
 $(H^p,H^p)$ for $1/2<p\leq 1 $.

\par At the same time,
$Hf,H_{\varepsilon}f$ and $H^*f$
 are of type $(L^p_w,L^p_w)$ for $1<p<\infty $ and $w\in A_p$ and  weak-type
 $(L^1_w,L^1_w)$ for $w\in A_1$
.  $Hf$ is of types $(H^p_w,L^p_w)$ and $(H^p_w,H^p_w)$ for
$1/2<p\leq 1 $ and $w\in A_1$( see \cite{LL}).

\par{\bf B.}  Riesz transform and the corresponding truncated
operators and maximal operator:
$$R_j  f(x)=\lim_{\varepsilon \rightarrow 0} R_j^\varepsilon f(x), R_j^\varepsilon
f(x)=c_n\int_{|y|>\varepsilon}\frac{y_j}{|y|^{n+1}}f(x-y)dy
  ~{\rm
 and }~
R_j^*f(x)=\sup_{\varepsilon >0} | R_j  f(x)| ,$$ $ j=1,\cdots,n, $
with $c_n=\frac{\Gamma((n+1)/2)}{\pi ^{(n+1)/2}}$.

\par It is well-known that  $R_j, R_j^\varepsilon$ and $R_j^*$ are of type $(L^s,L^s)$ for $1<s<\infty
$. Then $R_j, R_j^\varepsilon$ and $R_j^*$ satisfy the conditions of
Theorem 7.2.

\par It is also known that
$R_j, R_j^\varepsilon$ and $R_j^*$ are of  weak-type $(L^1,L^1)$,
and $R_jf$ is of type $(H^s,H^s)$ for $0<s\leq 1 $ (see
\cite{SW60,FS72}).

\par
 $R_jf$ and  $R_j^*f$
 are of type $(L^p_w,L^p_w)$ for $1<p<\infty $ and $w\in A_p$ and  weak-type
 $(L^1_w,L^1_w)$ for $w\in A_1$.
$R_jf$ is of types $(H^p_w,L^p_w)$ and  $(H^p_w,H^p_w)$ for
$n/(n+1)<p\leq 1 $ and $w\in A_1$( see \cite{LL}).

\par{\bf C.} The regular
singular integral operators and the corresponding truncated
operators and maximal operators defined in \cite{GRubio}:
$$T_{RS}f(x)=  K \ast f(x), T_{RS}^\varepsilon f(x)= K_\varepsilon
\ast f(x)~{\rm and} ~ T_{RS}^* f(x)=\sup_\varepsilon |T^\varepsilon
f(x)|$$ where the kernel $K$ satisfies (7.3) and a regular condition
$$|K(x-y)-K(x)|\leq |y|/|x|^{n+1}, ~~|x|>2|y|>0, $$
and $K_\varepsilon=K\chi_{\{|x|>\varepsilon\}}$ with
$\varepsilon>0$.
\par It is known that $T_{RS} ,  T_{RS}^\varepsilon   $ and  $T_{RS}^*    $  are of type $(L^s,L^s)$ for
$1<s<\infty$ (see \cite{GRubio}). Then $T_{RS} , T_{RS}^\varepsilon
$
 and  $T_{RS}^*    $ satisfy the conditions of Theorem 7.2.

\par It is also known that
  $T_{RS} ,  T_{RS}^\varepsilon   $ and  $T_{RS}^*    $ are of weak-type
 $(L^1,L^1)$,  and  $T_{RS}$ is of type $(H^p,L^p)$  for $n/(n+1)<p\leq 1 $.

\par  And
 $T_{RS} ,  T_{RS}^\varepsilon   $ and  $T_{RS}^*    $
 are of type $(L^p_w,L^p_w)$ for $1<p<\infty $ and $w\in A_p$ and  weak-type
 $(L^1_w,L^1_w)$ for $w\in A_1$. $T_{RS}$ is of type $(H^p_w,L^p_w)$
 for $n/(n+1)<p\leq 1 $ and $w\in A_1$ (see \cite{GRubio}).

\par{\bf D.} The R.Fefferman type
singular integral operator and the corresponding truncated operators
and maximal operators:
$$T_s f(x)= \lim_{\varepsilon \rightarrow 0}
T_\varepsilon f(x), T_s^\varepsilon f(x)=\int_{|y|>\varepsilon}
\frac{\Omega(\frac{y}{|y|})h(|y|)}{|y|^n} f(x-y)dy   {\rm~ and}~
T_s^* f(x)= \sup_{\varepsilon > 0} |T_s^\varepsilon f(x)|,$$
   where $h\in L^\infty ([0,\infty))$, and $\Omega $ is a homogeneous
   function of degree $0$, $\Omega \in L^\infty (S^{n-1}),$ and $\int _{{\bf
S}^{n-1}}\Omega(u)d\sigma(u)=0,\varepsilon>0$.

\par   $T_s ,  T_s^\varepsilon   $ and  $T_s^*    $  are of type $(L^s,L^s)$ for
$1<s<\infty$ (see, for example, \cite{DR}). Then  $T_s ,
T_s^\varepsilon   $ and  $T_s^*    $ satisfy the conditions of
Theorem 7.2.

\par It is also known that these  operators are  of   type $(L^s
_{w},L^s _{w})$ for $1<s<\infty$ and  $ w\in A_s$ (see \cite{DR}).


\par{\bf E. } The Fourier integrals operator and the
Carleson operator:
 $$S_N f(x)= \int_{|\xi|\leq N }\hat{f}(\xi)
e^{2\pi ix\xi} d\xi ~{\rm and }~ Cf(x)= \sup _{N>0}|S_Nf(x)|. $$

\par  It is known that
$$S_Nf=\frac{i}{2}( {\rm Mod}_NH {\rm Mod}_{-N}f-{\rm Mod}_{-N}H {\rm Mod}_{N}f),$$
where $H$ is Hilbert transform and ${\rm Mod}_{N}f(x)=e^{2\pi i
Nx}f(x).$ Then the estimate of $Cf$ can be concluded from the
estimate of
$$\bar{C}f(x)=\sup_{R>0}\left| {\rm p.v.} \int_{\bf
R}\frac{e^{iRy}}{x-y}f(y)dy\right|.$$
 It is easy to see that $\bar{C}f(x)$ can be controlled   by  the maximal Carleson operator
$$
C^*f(x)=\sup_{R>0}\sup_{\varepsilon >0}\left|   \int_{|x-y|>\varepsilon|
}\frac{e^{-iR(x-y)}}{x-y}f(y)dy\right|.
$$
It is known that $S_N, C, \bar{C} $ and $ C^*$ are bounded on $L^p$,
see \cite{Hunt,GTT}, and therefore they satisfy the conditions of
Theorem 7.2.

 \par It is also known that
$S_Nf $ is of  weak-type $(L^1({\bf R}),L^1({\bf R}))$ and type
$(H^p({\bf R}),L^p({\bf R}))$
 for $0<p\leq 1$, but
fails to be  of
  type $(L^1({\bf R}),L^1({\bf R}))$ \cite{Ko1},
  and $Cf$ fails to be  of
  weak-type $(L^1({\bf R}),L^1({\bf R}))$ \cite{Ko1} and weak-type $(H^1({\bf R}),L^1({\bf R}))$
  (see \cite {Zy}).

\par $S_Nf , Cf$ and $C^*f$
 are also of type $(L^{p} _{w}({\bf R}), {L}^{p}
  _{w}({\bf R})) $ for $1<p<\infty$ and $w\in A_p$, see \cite {HY,Grafakos}.

\par As we see, for the truncated operators
and maximal operators above, though the weak-type $(L^1,L^1)$
estimates are true, but the $(H^p,H^p)$ and $(H^p,L^p)$ estimates
fail (or is not known) for $0<p\leq 1$. Though the $(H^p,H^p)$ and
$(H^p,L^p)$ estimates hold for Hilbert transform for $p\in (1/2,1]$
and for the regular singular integral operators for $p\in
(n/(n+1),\leq 1] $, but for $p$ in the the remain ranges of $(0,1]$
these estimates fail (or are not known). In fact, for these
operators, no strong estimates are found for $p$ in the the remain
ranges of $(0,1]$. Here, replacing the $(H^p,H^p)$ and $(H^p,L^p)$
estimates, we establish
$(BL^{p,s}_{|x|^{\alpha}},BL^{p,s}_{|x|^{\alpha}})$ and
$(BL^{p,s}_{|x|^{\alpha}},L^{p}_{|x|^{\alpha}})$ estimates for these
operators for $p$ in the the remain ranges.
  \par From Theorems 7.2 and 7.3, we have

\par {\bf Corollary 7.1} Let $s,p$ and $\alpha$ as in Theorem 7.2,
then
 $H, H_{\varepsilon} , H^*, R_j, R_j^{\varepsilon} , R_j^*, T_{RS},$ $ T_{RS}^{\varepsilon} , T_{RS}^*, T_s ,
T_s^\varepsilon  , T_s^* , S_N, C ,\bar{C}   {\rm ~and ~} C^* ~{\rm
are ~ of ~type ~} (BL^{p,s}_{|x|^{\alpha}},BL^{p,s}_{|x|^{\alpha}})
. $
 \par Let $s,p$ and $\alpha$ as in Theorem 7.3, then
 $H,H_{\varepsilon} , H^*, R_j,R_j^{\varepsilon} , R_j^*, T_{RS}, T_{RS}^{\varepsilon} , T_{RS}^*, T_s ,
T_s^\varepsilon  ,$ $ T_s^* , S_N, C ,\bar{C}   {\rm ~and ~} C^*
~{\rm are ~ of ~type ~}
(BL^{p,s}_{|x|^{\alpha}},L^{p}_{|x|^{\alpha}}) . $

\par
For $S_N$ we also have  a convergence result as follows, whose proof
is similar to that of Theorem 6.4.

\par {\bf Theorem 7.4} Let $  1<s<\infty, 0<p< s,
n(\frac{p}{s}-1)< \alpha<n(p-1)$ and $\alpha\leq 0$.  Then
$$
\|S_N f-f\|_{L^p_{|x|^\alpha}} \rightarrow 0, R\rightarrow \infty
$$
for all $ f\in BL^{p,s}_{|x|^{\alpha}}$.

\par The following operators satisfy (7.1), so, for them the same
conclusions hold as those stated in Theorem 7.1. From these, we can
extend these operators    to some bounded operators on
$BL^{p,s}_{|x|^{\alpha}}$.

\par {\bf F. } Calder\'{o}n-Zygmund  operator and the corresponding truncated operators
and maximal operators defined  in \cite{Mey}.
\par Let $ K(x, y)$
be a locally integrable function defined off the diagonal $ x = y $
in ${\bf R}^n \times {\bf R}^n,$ which satisfies the standard
estimate:
$$|K(x, y)|\leq\frac{ c}{ |x- y|^{n}}, $$
  and,
for some $\delta > 0$,
$$|K(x, y) -K(z, y)| + |K(y, x) -K(y, z)|\leq
c \frac{|x- z|^\delta}{ |x- y|^{n+\delta }} , $$
 whenever $2|x - z| <
|x - y|.$
\par A linear operator $T_C:C^\infty_0 \rightarrow L^1_{loc}$ is a Calder\'{o}n-Zygmund
operator if $${\rm it~ extends~ to~ a ~bounded~ operator~ on }~L^2,
$$
 and there is a kernel $K$ satisfying the standard estimate
such that
$$T_Cf (x) =\int_ {{\bf R}^n} K(x, y)f (y) dy  $$
for any $f\in C^\infty_0 $ and $x$ not in  supp$(f )$.
   The corresponding truncated operators
and maximal operators are respectively defined as
$$T_C^\varepsilon f (x) =\int_ {|x-y|\geq \varepsilon} K(x, y)f (y) dy
~{\rm and }~ T_C^* f (x) =\sup_{\varepsilon>0}\left|T_C^\varepsilon
f (x)\right| .$$

\par $T_Cf, T_C^\varepsilon f$ and $T_C^* f$ are of type $(L^p,L^p)$ for $1<p<\infty$
   (see, for example, \cite{Mey}), and clearly they satisfy (7.1).
   So, theses operators satify the conditions of Theorem 7.1.

\par{\bf G.} The strongly singular integral operator
$$\widehat{T_bf}(\xi)=\frac{e^{i|\xi|^b}}{|\xi|^{nb/2}}\varphi(|\xi|)\hat{f}(\xi),$$
where $0<b<1$ and $\varphi$ is a smooth radial cut-off function with
$\varphi\equiv1$ on $\{|\xi|\geq 1\}$ and $\varphi\equiv0$ on
$\{|\xi|\leq 1/2\}$. The convolution form of $T_b$ can be written as
$$T_bf(x)={\rm p.v.~}\int_{{\bf R}^n}
\frac{e^{i|x-y|^{-b'}}}{|x-y|^n} \chi_E(|x-y|)f(y)dy, $$ where
$b'=b/(1-b)$ and $\chi_E$ is the characteristic function of the unit
interval $E=(0,1)\subset {\bf R}.$
\par  $T_b$ is of type $(L^p,L^p)$ for
$1<p<\infty$ (see \cite{Hir,Wain}), then, $T_b$ satisfies the
conditions of Theorem 7.1.


\par {\bf  H.}    The
oscillatory singular integrals operators
$$T_of(x)={\rm p.v.} \int_{{\bf
R}^n}e^{\lambda\Phi(x,y)}K(x,y)\varphi(x,y)f(y)dy,$$ where $K(x,y)$
is a Calder\'{o}n-Zygmund kernel,  $\varphi \in C_0^\infty({\bf
R}^n\times {\bf R}^n),$  $\lambda \in {\bf R}$, and $ \Phi(x,y)$ is
real-valued.

\par  $T_o$ is of type $(L^p,L^p), 1<p<\infty$, for the real bilinear
form
 $\Phi(x,y)=(Bx,y), \varphi=1$ (see \cite{PhongStein}), for the
 polynomial $\Phi(x,y)=P(x,y), \varphi=1$ (see \cite{ RicciStein}), and
 for the real analytic function
 $\Phi(x,y)$ (see \cite{Pan}). Then, in these cases, $T_o$ satisfies the
conditions of Theorem 7.1.


\par From Theorem 7.1 we have

\par {\bf Corollary 7.2} Let $s,p$ and $\alpha$ as in Theorem 7.1.
Then (7.2) holds   for $ T_C , T_C^\varepsilon  , T_C^* , T_b$ and $
T_o $ stated as above.

\par Let $FBL^{p,s}_{|x|^{\alpha}}$ be the set of all finite linear
combination of $(p,s,\alpha)$-blocks. Let $f$ be an element of
$FBL^{p,s}_{|x|^{\alpha}}$ and pick a representation of
$f=\sum_{j=1}^N \lambda_ja_j$ such that
$$\|f\|_{BL^{p,s}_{|x|^{\alpha}}}\approx \left(\sum\limits_{j=1}^{N}|\lambda _j|^{\bar{p}}\right)^{1/{\bar{p}}}.$$
Using Corollary 7.2, Theorem 3.1 and the sublinearity of the
quantity $\|\cdot\|_{BL^{p,s}_{|x|^{\alpha}}}$, it follows that, for
these operators $T$ in Corollary 7.2,
$$\|Tf\|_{BL^{p,s}_{|x|^{\alpha}}}\leq C \|f\|_{BL^{p,s}_{|x|^{\alpha}}} $$
if $s,p$ and $\alpha$ under the conditions of  Theorem 7.1, where
$C$ is independent of $f$.

Since   $FBL^{p,s}_{|x|^{\alpha}}$ is dense in
$BL^{p,s}_{|x|^{\alpha}}$ and  $BL^{p,s}_{|x|^{\alpha}}$ is complete
for $0<s\leq \infty,0<p< s$ and $-n(1-p/s) <\alpha \leq 0$ (Theorem
2.6), then we can extend these operators  $ T $ to some bounded
operators on $BL^{p,s}_{|x|^{\alpha}}$, so we have

\par {\bf Theorem 7.5} Let
 $    1< s < \infty,0< p < s, - n(1-p/s)< \alpha < n(p-1)$
and  $\alpha \leq 0$. Then the operators $ T_C , T_C^\varepsilon  ,
T_C^* , T_b$ and $ T_o $ initially defined for $f\in L^s$ with
compact support admit respectively an extension which is bounded on
$BL^{p,s}_{|x|^{\alpha}}$.

\par \section*
{\bf 8. Remarks}
\par (1)  When $p>1$, the spaces $
BL^{p,s}_{|x|^{\alpha}}$ can not be replaced  by $
\tilde{B}L^{p,s}_{|x|^{\alpha}}$ defined below for certain problems,
such as the convergence of Bochner-Riesz means and the spherical
means. In other words,
 the quasinorm of a function $f(x)$
$$ \|f\|_{ BL^{p,s}_{|x|^{\alpha}}}=\inf_{f(x)=
\sum_{i=-\infty}^\infty\lambda_ia_i(x), a(x)~ {\rm are}
~(p,s,\alpha)-{\rm blocks}}
 \sum_{i=-\infty}^\infty|\lambda_i| $$
in Definition 2.2 can not be replaced  by
$$\|f\|_{\tilde{B}L^{p,s}_{|x|^{\alpha}}}=
 \inf_{f(x)=
\sum_{i=-\infty}^\infty\lambda_ia_i(x), a(x)~ {\rm are}
~(p,s,\alpha)-{\rm blocks}}
\left(\sum\limits_{k=-\infty}^{\infty}|\lambda _k|^p\right)^{1/p}.$$

\par Let
$0< s \leq \infty, 0<p<\infty$ and $ -\infty<  \alpha<\infty $.
Denote
\begin{eqnarray*}
 \tilde{B}L^{p,s}_{|x|^{\alpha}}=\{  f&: &
 f=\sum\limits_{k=-\infty}^{\infty} \lambda _ka_k ,
 \\&&\textrm { where
 each $a_k$ is a $(p,s, \alpha)$-block on ${\bf R}^n$,}
 \sum\limits_{k=-\infty}^{\infty} |\lambda _k|^p <+ \infty \}.
 \end{eqnarray*}
 Here the
 "convergence" means a.e. convergence.

\par
And denote
$$\tilde{B}L^{p }_{|x|^{\alpha}}= \{f: Mf \in {L}^{p }_{|x|^{\alpha}} \},$$
where $Mf$ is  Hardy-littlewood maximal function. Then we have

\par {\bf Theorem 8.1}  Let $0< s\leq \infty, 0<p<\infty$ and $
0\leq\alpha <\infty$.  Then
\begin{eqnarray*}
\tilde{B}L^{p}_{|x|^{\alpha}}\subset
\tilde{B}L^{p,s}_{|x|^{\alpha}}.
\end{eqnarray*}
\par {\bf Proof}.
 Let $f\in \tilde{B}L^{p}_{|x|^{\alpha}}$. For $k\in {\bf Z}$, let
$$E_k=\{x\in {\bf R}^n:Mf(x)>2^k\}=\bigcup_{j=1}^\infty Q_j^k,
 {\rm ~~the~ Whitney ~decomposition ~of~} E_k.$$
  Clearly, we have
$$\bigcup _{k=-\infty}^\infty E_k={\bf R}^n,
\bigcup_{j,k}Q_j^k\setminus E_{k+1}={\bf R}^n, {\rm~and ~}
E_{k+1}\subset E_k.$$ Thus
$$f(x)=\sum_{k,j}\chi_{Q_j^k\setminus
E_{k+1}}(x)f(x)=\sum_{k,j}\lambda _j^ka_j^k.$$
Here
$$\lambda_j^k=C2^k(\mu_\alpha(Q_j^k))^{1/p},$$
$$a_j^k(x)=\left(C2^k(\mu_\alpha(Q_j^k))^{1/p}\right)^{-1}\chi_{Q_j^k\setminus
E_{k+1}}(x)f(x),$$
 where
$\mu_\alpha(Q_j^k)=\int_{Q_j^k}|x|^{\alpha}dx.$ Then we have
\begin{eqnarray*}
\sum_{j,k}|\lambda_j^k|^p&=&C^p\sum_{j,k}2^{kp}\mu_\alpha(Q_j^k)=C^p\sum_{k}2^{kp}\mu_\alpha(E_k)
\\&=&C^p\sum_{k}2^{kp}\mu_\alpha(\{x\in
{\bf R}^n:Mf(x)>2^k\})
\\&\leq& C\int_{{\bf
R}^n}(Mf(x))^p|x|^{\alpha}dx=C\|f\|^p_{
\tilde{B}L^{p}_{|x|^{\alpha}}}.
\end{eqnarray*}
Next, let us show that each $a_j^k$ is a $(p,s,\alpha)$-block. In
fact, noticing that $2^k<Mf(x)\leq 2^{k+1}$ on $ Q_j^k\setminus
E_{k+1} $, it follows that $|f(x)|<Mf(x)\leq 2^{k+1}$. Then we have
\begin{eqnarray*}
\|a_j^k\|_{L^{s}} &=&C\left(
2^{k}\mu_\alpha(Q_j^k)^{1/p}\right)^{-1}\left(\int_{Q_j^k\setminus
E_{k+1}}|f(x)|^sdx\right)^{1/s}
\\&\leq& C2^{-k}\mu_\alpha(Q_j^k)^{-1/p}2^k|Q_j^k|^{1/s}
\\&\leq& C |Q_j^k|^{-\alpha/pn-1/p+1/s},
\end{eqnarray*}
since  $\mu_\alpha(Q_j^k)\geq |Q_j^k|^{\alpha/n+1}$ when $\alpha\geq
0$.  Thus, we finish the proof of Theorem 8.1.

\par By the boundedness of Hardy-Littlewood maximal function $Mf$ on
${L}^{p}_{|x|^{\alpha}}$ we know  that, when  $1<p<\infty$ and $
-n<\alpha<n(p-1) $,
$$\tilde{B}L^{p}_{|x|^{\alpha}}=
{L}^{p}_{|x|^{\alpha}}. $$
 This and Theorem 8.1 give
 that
\par {\bf Corollary 8.1}  Let $0< s \leq \infty, 1<p<\infty$ and $
0\leq  \alpha<n(p-1) $. We have
$$ L^{p}_{|x|^{\alpha}}\subset
\tilde{B}L^{p,s}_{|x|^{\alpha}}. $$ In particular,
$$ L^{p}\subset
\tilde{B}L^{p,s}_{|x|^{0}}. $$

 \par It is known that  Bochner-Riesz means $B_R^\lambda$ for
 $ \lambda \leq \frac{n-1}{2}$ and $ p\leq p_\lambda'$ and the spherical maximal operator $\mathcal{M}$
 for $p\leq \frac{n}{n-1}$ and
$n\geq 2$ fail to be of type  $(L^{p},L^{p})$, see \cite{S3}. And
it follows that
\par {\bf Corollary 8.2} Let $0< s \leq \infty$. Then $B_R^\lambda$ for
 $ \lambda \leq \frac{n-1}{2}$ and $ p\leq p_\lambda'$, and
$\mathcal{M}$ for $  1<p\leq \frac{n}{n-1}$ and $n\geq 2$ fail to be
of type $(\tilde{B}L^{p,s}_{|x|^{0}},L^{p})$.


\par (2) For Fourier integrals operator $S_N$ and   Carleson operator
$C$, the second conclusion of Corollary 7.1 is  sharp in the sense
that

\par {\bf Theorem 8.2}
A) If  $ s=1,p=1$ and $ \alpha=0, $ then there exists $f\in
BL^{p,s}_{|x|^{\alpha}}({\bf R})$ such that
 $S_Nf$
   fails to be of type $(BL^{p,s}_{|x|^{\alpha}}({\bf R}), {L}^{p}
  _{|x|^{\alpha}}({\bf R})) $.
\par B) If
 $  s=1, 0<p\leq 1$ and $ -1 <\alpha \leq p-1, $
  then there exists $f\in BL^{p,s}_{|x|^{\alpha}}({\bf R})$ such that
 $Cf$
   fails to be of type $(BL^{p,s}_{|x|^{\alpha}}({\bf R}), {L}^{p}
  _{|x|^{\alpha}}({\bf R})) $.

\par {\bf Proof} A) is clear since $\dot{B}L^{1,1}_{|x|^{0}} =L^1$
\cite{Ko1}. Let us prove B).  When $  s=1, 0<p\leq 1$ and $ -1
<\alpha \leq p-1, $ by Theorem 2.3B), we see that $L^1\subset
\dot{B}L ^{p,s} _{|x|^\alpha}.$ Then Kolmogorov's example shows that
there exists $f\in \dot{ B}   L ^{p,s} _{|x|^\alpha}$ such that
$$
\limsup_{N\rightarrow \infty}S_Nf(x)=\infty,~~{\rm a.e.}.
$$
It follows that
$$
\limsup_{N\rightarrow \infty}S_Nf(x)=\infty,~~\mu_\alpha{\rm -a.e.}.
$$
Then $Cf(x)=\infty,~~\mu_\alpha{\rm -a.e.}$. This is a contradiction
of $\|Cf\|_{L^{p} _{|x|^{\alpha}}({\bf R})} \leq C\|f\|_{\dot{B}
L^{p,s} _{|x|^{\alpha}}({\bf R})}$.

{ }
\par
Shunchao Long
\par Department of Mathematics,
\par Xiangtan University,
 \par Hunan, 411105  P.R.China
 \par E-mail: sclong@xtu.edu.cn


\begin{thebibliography}{s2}


\bibitem {Ande} K.F.Andersen, Weighted norm inequalities for Bochner-Riesz
spherical summation multipliers, Proc. Amer. Math. Soc. 103(1988),
165-170.

\bibitem {BRV}O.Blasco, A.Ruiz and  L.Vega,
Non interpolation in Morrey-Campanato and block spaces,
  Ann. Sc. Norm. Super. Pisa, Cl. Sci. IV. Ser. (1)28(1999), 31-40.


\bibitem {B}
S.Bochner, Summation of multiple Fourier series by spherical means,
 Trans. Amer. Math. Soc. 40(1936), 175-207.

\bibitem {B3}
J.Bourgain, Besicovitch type maximal operators and applications to
Fourier analysis, Geom. Funct. Anal. 1(1991), 147-187.


\bibitem {B4}
 J.Bourgain, ¡°On the restriction and multiplier problems in ${\bf R}^3$¡±
in Geometric Aspects of Functional Analysis: Israel Seminar (GAFA)
1989/90, Lecture Notes in Math. 1469, Springer, 1991, 179-191.


 \bibitem {B1}
 J.Bourgain, Averages in the plane over convex curves and maximal
operators, J. Anal. Math. 47(1986), 69-85.


\bibitem {B2}
J.Bourgain, Estimations de certaines functions maximales, C. R.
Acad. Sci. Paris 301(1985), 499-502.


\bibitem {Ca}
 C.Calder\'{o}n, Lacunary spherical means, Illinois J. Math. 23(1979), 476-484.



\bibitem {CRV}
A.Carbery, J.L.Rubio de Francia and L.Vega, Almost everywhere
summability of Fourier integrals. J. London Math. Soc. 38(1988),
513-524.

\bibitem {CaS}
L.Carleson  and P.Sj\"{o}lin,   Oscillatory Integrals and a
Multiplier Problem
 for the Disc, Studia Math. 44(1972), 287-299.


\bibitem {Christ1}
M.Christ, Weak type (1,1) bounds for rough operators, Ann. Math.
128(1988), 19-42.

\bibitem {Christ2}
M.Christ, Weak type endpoint bounds for Bochner-Riesz multipliers,
Rev. Mat. Iberoamericana, (1)3(1987), 25-31.

\bibitem {Coif}
R.Coifman, Characterizations of Fourier transforms of Hardy spaces,
Proc. Nat. Acad. Sci. U.S.A. 71(1971), 4133-4134.

\bibitem {CW}
R.Coifman and G.Weiss, Extensions of Hardy spaces and their use in
analysis, Bull.
 Amer. Math. Soc. 83(1977), 569-645.

\bibitem {CW2}
 R.R.Coifman and G.Weiss, Review of the book:Littlewood-Paley
and multiplier theory, Bull. Amer. Math. Soc. 84(1978), 242-250.

\bibitem {Cord1}
A.C\'{o}rdoba, A note on Bochner-Riesz operators, Duke Math. J. 46
(1979), 505-511.

\bibitem {DC}
 K.M.Davis, Y.C.Chang, Lectures on Bochner-Riesz means,
  London Mathematical Society Lecture Note Series. 114, Cambridge Univ. Press. Cambridge.
  1987.

 \bibitem {DV}
  J.Duoandikoetxea and L.Vega, Spherical
means and weighted inequalities, J. London Math. Soc. 53(1996),
343-353.


\bibitem {DS}
J.Duoandikoetxea and E.Seijo,Weighted inequalities for some
spherical maximal operators, Illinois J. Math.  (4)46(2002),
1299-1312.


\bibitem {DR}
J.Duoandikoetxea and L.J.Rubio de Francia,  Maximal and singular
integral operators via Fourier transform estimates, Invent. Math.
(3)84(1986), 541-561.


\bibitem {F1}
C.Fefferman, A note on spherical summation multipliers, Israel J.
Math. 15(1973), 44-52.

\bibitem {F2}
C.Fefferman, Inequalities for strongly singular convolution
operators, Acta Math. 124(1970), 9-36.

\bibitem {FS72}
C.Fefferman and E.M.Stein, $H^p$ spaces of several variables, Acta
Math. 129(1972), 137-194.


\bibitem {Gc}
J.Garc\'{\i}a-Cuerva, Weighted $H^p$ spaces. Dissertations Math.
162(1979), 1-63.


\bibitem {GRubio}
 J.Garc\'{\i}a-Cuerva and J.Rubio de Francia, Weighted norm
inequalities and related topics, North-Holland, Amsterdam, 1985.

\bibitem {Grafakos}
L.Grafakos, Classical and modern Fourier analysis, Pearson
Education, 2004.

\bibitem {GTT}
L.Grafakos, T.Tao and E.Terwilleger, $L^p$ bounds for a maximal
dyadic sum operator, Math. Z. 246(2004), 321-337.

\bibitem {H}
C.S.Herz, On the mean inversion of Fourier and Hankel transforms,
 Proc. Nat. Acad. Sci. U.S.A. 40(1954), 996-999.


\bibitem {Hir1}
  I.I.Hirschman, Multiplier transformations. II, Duke Math. J. 28(1961), 45-56.

\bibitem {Hir}
 I.I.Hirschmann, On multiplier transformations, Duke Math. J. 26
(1959), 221-242.

\bibitem {Hor}
L.H\"{o}rmander,  Oscillatory integrals and multipliers on $FL^p$,
Arkiv. Math. 11(1973), 1-11.

\bibitem {Hunt}
R.A.Hunt, On the convergence of Fourier series. In orthogonal
expansions and their continuous analogues, pages 235-255,
Carbondale, IL, 1968.

\bibitem {HY}
R.A.Hunt and W.S.Young, A weighted norm inequality for Fourier
series,. Bull. Amer. Math. Soc. 80(1974), 274-277.

\bibitem {Ko1}
A.N.Kolmogorov, Une s\'{e}rie de Fourier-Lebesgue divergente presque
partout, Fundamenta Math. 4(1923), 324-329.



\bibitem {Lee}
S.Lee, Improved bounds for Bochner-Riesz and maximal Bochner-Riesz
operators, Duke Math.J. (1)122(2004), 205-232.

\bibitem {LL}
M.Y.Lee and C.C.Lin, The molecular characterization of weighted
Hardy spaces, J. Funct. Anal. 188(2002), 442-460.


\bibitem {LY}
X.Li and D.Yang, Boundedness of some sublinear operator s on Herz
spaces,   Illinois J. of Math.  40(1996), 484-501.


\bibitem {Lo} R.Long, The spaces generated by blocks, Scientia Sinica A XXI II (1984), 16-26.


\bibitem {Longs}
S.Long, Convergence of Fourier series, Bochner-Riesz means and
restriction problems, Ph.D.Thesis, Xiangtan University, Xiangtan,
2008.

\bibitem {Longs2}
S.Long, Convergence of Fourier series at or beyond endpoint,
preprint.


\bibitem {LTW1}
S.Z.Lu, M.Taibleson and G.Weiss, Spaces generated by blocks, Beijing
Normal University Press, Beijing, 1989.

\bibitem {An}
Marco Annoni. Almost everywhere convergence for modified
Bochner-Riesz means at the critical index for $p \geq 2$,
Ph.D.Thesis, University of Missouri- Columbia,  2010.

\bibitem {Mey}
Y.Meyer, Wavelets and operators,  Cambridge Studies in Advanced
Math. vol. 37, Cambridge Univ. Press, Cambridge, 1992.


\bibitem {Pan6}
Y.Pan, Oscillatory singular integrals on $L^p$ and Hardy spaces,
Proc. Amer. Math. Soc. (9)124 (1996), 2821-2825.


\bibitem {Pan3}
Y.Pan, Hardy spaces and oscillatory singular integrals, Rev. Mat
Iberoamericana, 7(1991), 55-64.


\bibitem {Pan}
Y.Pan, Uniform estimates for oscillatory integral operators, J.
Funct. Anal. 100 (1991), 207-220.


\bibitem {PhongStein}
D.H.Phong and E.M.Stein, Hilbert integrals, singular integrals and
Radon transforms, I, Acta Math. 57(1987), 179-194.

\bibitem {RicciStein}
 F.Ricci and E.M.Stein, Harmonic analysis on nilpotent groups and singular
integrals, I, J. Funct. Anal. 73(1987), 179-194.


\bibitem {Rubio}
J.L.Rubio.de.Francia, Weighted norm inequalities and vector valued
inequalities, Lecture Notes in Math. 908, 86-101, 1982.




\bibitem {Sato}
 S.Sato,  Divergence of the Bochner-Riesz means in the weighted Hardy spaces,
 Studia Math. 118(1996), 261-275.






\bibitem {Seeg}
A.Seeger, Endpoint inequalities for Bochner-Riesz multipliers in the
plane, Pacific J. Math. (2)174(1996), 543-553.

\bibitem {STW1}
A.Seeger, T.Tao and J.Wright, Pointwise convergence of lacunary
spherical means, Harmonic Analysis at Mount Holyoke South Hadley,
MA, 2001, Contemp. Math. vol. 320, Amer. Math. Soc. Providence, RI
(2003),  341-351.


\bibitem {STW2}
 A.Seeger, T.Tao and  J.Wright, Singular maximal functions and
Radon transforms near $L^1$, Amer. J. Math. 126(2004), 607-647.

\bibitem {STW3}
 A.Seeger, T.Tao and J.Wright, Endpoint mapping properties of
spherical maximal operators, J. Inst. Math. Jussieu 2(2003),
109-144.


\bibitem {SWW1}
A.Seeger, S.Wainger and J.Wright, Pointwise convergence of spherical
means, Math. Proc. Camb. Phil. Soc. 118(1995), 115-124.


\bibitem {SWW2}
A.Seeger, S.Wainger and J.Wright, Spherical maximal operator on
radial functions, Math. Nachr. 187(1997), 241-265.



\bibitem {SS}
X.Shi and Q.Sun, Weighted norm inequalities for Bochner-Riesz
operators and singular integral operators, Proc. Amer. Math. Soc.
116(1992), 665-673.

\bibitem {Sjolin76}
 P.Sj\"{o}lin, $L^p$ estimates for strongly singular convolution
operators in ${\bf R}^n$, Ark. Mat. 14 (1976), 59-64.















\bibitem {S3}
E.M.Stein, Harmonic analysis: real-variable methods, orthogonality
and
 oscillatory integrals, Princeton Univ. Press, Princeton, N. J., 1993,

\bibitem {S2}
 E.M.Stein,   On limits of sequences of operators. Ann.
Math. 74(I)(1961), 140-171.

\bibitem {St61970}
E.M.Stein, Singular integrals and differentiability properties of
functions, Princeton Univ. Press, Princeton, N. J., 1970.

\bibitem {Stei}
E.M.Stein, An $H^1$ function with non-summable Fourier expansion,
Lecture Notes in Math. 992, 1983,  193-200.


\bibitem {S4}
E.M.Stein, Localization and suminability of multiple Fourier series,
Acta Math. 100(1958), 93-147.


\bibitem {S1}
E.M.Stein, Maximal functions: spherical means, Proc. Nat. Acad. Sci.
U.S.A. 73(1976), 2174-2175.



\bibitem {SW60}
 E.M.Stein and G.Weiss, On the theory of harmonic functions of
several variables, Acta Math. 103(1960), 25-62.



\bibitem {SG}
E.M.Stein and G.Weiss, Introduction to Fourier analysis on Euclidean
spaces, Princeton Univ. Press, Princeton, N. J., 1971.


\bibitem {STW}
E.M.Stein, M.Taibleson and G.Weiss, Weak type estimates for maximal
operators on certain $H^p$ spaces, Rend. Circ. Mat. Palermo (2)
Suppl. 1 (1981), 81-97.


\bibitem {ST}
J.-O.Str\"{o}mberg and A.Torchinsky, Weighted Hardy spaces, Lecture
 Notes in Math. 1381, Springer, 1989.



\bibitem {TW1}
M.Taibleson and G.Weiss, Certain function spaces associated with
a.e. convergence of Fourier series. In:Proc. Conf. on Harmonic
Analysis in honor of Zygmund, Woodsworth, vol.1, (1983), 95-113.

\bibitem {TW2}
M.H.Taibleson and G.Weiss, The molecular characterization of certain
Hardy spaces, Asterisque 77(1980), 67-149.




\bibitem {T1}
T.Tao, Weak-type endpoint bounds for Riesz means, Proc. Amer. Math.
Soc. (9)124(1996),  2797-2805.

\bibitem {T2}
T.Tao, The weak-type endpoint Bochner-Riesz conjecture and related
topics, Indiana J. of Math. 47(1998), 1097-1124.


\bibitem {TV1}
 T.Tao and A.Vargas, A bilinear approach to cone multipliers,
I: Restriction estimates, Geom. Funct. Anal. 10(2000), 185-215.


\bibitem {TV2}
T.Tao and A.Vargas, A bilinear approach to cone multipliers, II:
Application, Geom. Funct. Anal. 10(2000), 216-258.


\bibitem {TVV}
 T.Tao, A.Vargas and L.Vega, A bilinear approach to the
restriction and Kakeya conjectures, J. Amer. Math. Soc. 11(1998),
967-1000.


\bibitem {To}
P.Tomas, Restriction theorems for the Fourier transform, Proc. Symp.
Pure Math. XXXV, vol 35, Amer.Math.Soc.Providence, RI, 1979,
111-114.


\bibitem {Wain}
 S.Wainger, Special trigonometric series in $k$
dimensions, Mem. Amer. Math. Soc. 59 (1965).


\bibitem {W}
T.Wolff, An improved bound for Kakeya type maximal functions, Rev.
Mat. Iberoamericana, 11(1995), 651-674.



\bibitem {Zy}
A.Zygmund, Trigonometric series, vol. 1, 2nd, Cambridge Univ. Press,
New York, 1959.
























































































\end{thebibliography}
\end{document}